\documentclass[11pt,reqno]{amsproc}

\usepackage[utf8]{inputenc}
\usepackage[margin=1in]{geometry}
\usepackage{amsmath,amssymb,amsfonts,amsthm,mathrsfs,mathtools}
\usepackage{enumitem}
\usepackage{booktabs,tabularx,array}
\usepackage[dvipsnames]{xcolor}
\usepackage{graphicx}
\usepackage{hyperref}
\hypersetup{colorlinks=true,citecolor=magenta,linkcolor=blue,filecolor=teal,urlcolor=cyan}
\numberwithin{equation}{section}

\newtheorem{theorem}{Theorem}[section]
\newtheorem{proposition}[theorem]{Proposition}
\newtheorem{lemma}[theorem]{Lemma}
\newtheorem{corollary}[theorem]{Corollary}
\theoremstyle{definition}

\theoremstyle{remark}
\newtheorem{remark}[theorem]{Remark}

\newcommand{\RR}{\mathbb{R}}

\newcommand{\ZZ}{\mathbb{Z}}
\newcommand{\CC}{\mathbb{C}}
\newcommand{\cS}{\mathcal{S}}
\newcommand{\cF}{\mathcal{F}}
\newcommand{\cG}{\mathcal{G}}
\newcommand{\cV}{\mathcal{V}}

\newcommand{\dd}{\mathop{}\!\mathrm{d}}
\makeatletter
\renewcommand{\subsubsection}{%
  \@startsection{subsubsection}{3}%
    \z@{.5\linespacing\@plus.7\linespacing}{-.5em}%
    {\normalfont\bfseries}}
\makeatother
\allowdisplaybreaks

\usepackage{etoolbox}

\makeatletter
\patchcmd{\@maketitle}
  {\global\topskip8pc\relax}
  {\global\topskip6pc\relax}
  {}
  {\PackageError{title-spacing}
    {Could not patch the title spacing}
    {Check the document class.}}
\makeatother

\begin{document}

\title{Finite-time blowup in the inviscid GOY and Sabra models}
\author{Taehun Kim}
\date{}

\begin{abstract}
For the inviscid GOY and Sabra shell models with zero boundary conditions
at the lowest shells, we construct strong solutions with finitely 
supported initial data whose vorticity amplitudes blow up in the 
supremum norm in finite time, at the self-similar rate $(T-t)^{-1}$.
The blowup is stable under small perturbations of the initial data 
in the same norm, including complex-valued ones.

Guided by the work of Mailybaev, we use a travelling wave of the 
renormalized equations to obtain blowup in the original variables.
We construct such a wave for the equations extended to all integer 
shells, decaying in both directions. We then show that the dynamics 
near this wave contract once the directions generated by time-translation 
and phase symmetries are removed by adjusting the symmetry parameters 
at each step (a modulation argument). This yields blowup solutions 
on the positive shells satisfying the original boundary conditions.
Computer-assisted interval arithmetic is used to establish the 
existence of the travelling wave and the contraction estimate 
needed for the modulation argument.
\end{abstract}

\maketitle

\section{Introduction}
\label{sec:intro}

Shell models provide a simplified description of energy transfer 
in fluid flows; see \cite{Biferale2003} for an overview.
They divide Fourier space into shells of geometrically spaced 
frequencies, for example $2^{n} \le |\xi| < 2^{n+1}$ in the dyadic 
case, and associate a single real or complex velocity amplitude with each 
shell. These models retain quadratic interactions between different 
frequency scales with relatively few dynamical variables and formally 
conserve energy in the inviscid case.
For a formal derivation of a dyadic model based on a Littlewood-Paley
decomposition of the Euler equations, see
\cite[Section 2]{CheskidovFriedlanderPavlovic2007}.
Two of the most widely studied shell models are the GOY
(Gledzer-Ohkitani-Yamada) model \cite{Gledzer1973, OhkitaniYamada1989} 
and its modification, the Sabra model \cite{LvovEtAl1998}. 

In this paper, we study the inviscid, unforced GOY and Sabra models
with specific parameters.\footnote{We follow the conventions of 
\cite[Eqs. (1) and (21)]{LvovEtAl1998}, with viscosity and forcing set to zero, 
$k_{n} = 2^{n}$, and $(a, b, c) = (1, -1/2, -1/2)$.} The GOY model is given by
\begin{equation}\label{eq:main-GOY}
\begin{aligned}
\frac{\dd}{\dd t} U_{n} &= i \left ( 2^{n+1}U_{n+2}U_{n+1} - 2^{n-1}U_{n+1}U_{n-1}
-2^{n-2}U_{n-1}U_{n-2}\right )^{\ast} \quad (n \ge 1), \\
U_{-1}&=0, \quad U_{0}=0,
\end{aligned}
\end{equation}
while the Sabra model is given by
\begin{equation}\label{eq:main-Sabra}
\begin{aligned}
\frac{\dd}{\dd t}U_{n} &= i\left ( 2^{n+1}U_{n+2}U_{n+1}^{\ast}-2^{n-1}U_{n+1}U_{n-1}^{\ast}
+2^{n-2}U_{n-1}U_{n-2}\right ) \quad (n \ge 1), \\
U_{-1}&=0, \quad U_{0}=0.
\end{aligned}
\end{equation}
In both models, $(U_{n}(t))_{n\ge 1}$ is a complex-valued sequence, and the 
superscript $\ast$ denotes complex conjugation.

There is a vast mathematical literature on shell models.
For the inviscid GOY and Sabra models, 
\cite{ConstantinLevantTiti2007}
establishes local well-posedness for sufficiently regular initial
data, a continuation criterion, and
global existence of finite-energy weak solutions.
For the viscous models, \cite{ConstantinLevantTiti2006} establishes 
global regularity. Further quantitative estimates are obtained in
\cite{VincenziGibbon2021, GibbonVincenzi2025}.
Stochastic versions of GOY and Sabra models have also been studied
rigorously; see
\cite{BarbatoEtAl2006, BessaihFerrario2012, BarbatoMorandin2013,
BessaihHausenblasRazafimandimby2016, FernandoHausenblasRazafimandimby2016}.
For other shell models, rigorous results include energy
dissipation in the inviscid case
\cite{CheskidovFriedlanderPavlovic2010, BarbatoFlandoliMorandin2011},
vanishing-viscosity limits
\cite{CheskidovFriedlander2009, FriedlanderGlattHoltzVicol2016},
and uniqueness and nonuniqueness of Leray-Hopf solutions
\cite{Filonov2017, FilonovKhodunov2021,Palasek2025}.
For further results and references, see
\cite{CheskidovDaiFriedlander2023}.

Shell models also provide a setting for studying possible
mechanisms of finite-time blowup in the Euler equations,
particularly through energy transfer to higher frequencies.
Finite-time blowup in the Katz-Pavlovi\'c model contrasts
with global regularity for the Obukhov model, 
despite their common scaling and formal
energy conservation \cite{KatzPavlovic2005, KiselevZlatos2005}.
This contrast highlights the importance of the specific
nonlinear interactions. The blowup argument for the Katz-Pavlovi\'c
model exploits the preservation of nonnegative amplitudes and the 
resulting one-way transfer of energy to higher shells.
Blowup has also been established for certain mixtures of the 
Katz-Pavlovi\'c and Obukhov interactions by a contradiction 
argument, even though nonnegativity of the amplitudes is not 
preserved \cite{JeongLi2015}.

For the inviscid GOY and Sabra models, which lack this 
positivity-preserving structure, numerical studies provide
evidence for finite-time blowup; see \cite{GilsonDombre1997,
Mailybaev2013Blowup, CST25} 
and the references therein. These numerical studies motivate a 
rigorous construction of finite-time blowup solutions to 
\eqref{eq:main-GOY} and \eqref{eq:main-Sabra}, which is the main goal 
of this paper. In particular, for the Sabra model with the parameters 
in \eqref{eq:main-Sabra}, our construction gives an affirmative 
answer to the question, raised in 
\cite[Section VIII]{ConstantinLevantTiti2007}, of whether initially 
smooth solutions can lose regularity in finite time.

Another motivation for studying finite-time blowup comes from
spontaneous stochasticity. This phenomenon refers to a situation 
in which, for fixed deterministic initial data, the solution 
distribution at a fixed time converges to a limit that is not 
concentrated at a single state as viscosity and random perturbations 
vanish in an appropriate joint limit \cite{EyinkBandak2020}.
In fluid models, this provides a mechanism by which microscopic
uncertainties, such as thermal fluctuations, can affect
predictions at much larger scales.

In shell models, one mechanism for spontaneous stochasticity is
the amplification of uncertainty at high frequencies, followed
by its propagation to lower frequencies through nonlinear
interactions, an inverse cascade of uncertainty
\cite{BandakEtAl2024Thermal, Mailybaev2016Inviscid}.
When the shell amplitudes decay sufficiently slowly, the time scale
on which nonlinear interactions occur becomes shorter at higher 
frequencies, which can facilitate this transfer of uncertainty in 
finite time. For example, the numerical study 
\cite{BandakEtAl2024Thermal} investigates this process in the Sabra 
model with thermal noise, starting from initial states whose amplitudes 
decay slowly at high frequencies. The high-frequency activity needed 
for this mechanism need not be present initially; it may emerge 
dynamically from regular initial data as the solution approaches 
finite-time blowup, potentially allowing spontaneous stochasticity 
after the blowup time.

This scenario of spontaneous stochasticity after blowup is supported 
by analytical arguments and numerical simulations in a real subsystem 
common to the GOY and Sabra models with nonzero large-scale boundary 
data \cite{Mailybaev2016Velocity}, and in the complex Sabra model
\cite{Mailybaev2016Inviscid}. The role of 
finite-time singularities in the emergence of spontaneous stochasticity 
has also been studied beyond shell models; see
\cite{DrivasMailybaevRaibekas2024, GrottoRomitoViviani2024, 
OrtizCampolinaMailybaev2026}. We view the construction of blowup 
solutions for \eqref{eq:main-GOY} and \eqref{eq:main-Sabra} as a starting 
point for a rigorous investigation of whether spontaneous stochasticity 
can emerge in these models after blowup.

\subsection{Main result}

To present our main result concisely, we rewrite both
\eqref{eq:main-GOY} and \eqref{eq:main-Sabra} in their
vorticity forms. For \eqref{eq:main-GOY}, let 
$\omega_{n} = -i2^{n}U_{n}$.
The complex quantity $\omega_{n}$ plays the role of a vorticity
amplitude. In these variables, \eqref{eq:main-GOY} becomes
\begin{equation}\label{eq:main-GOY-renorm}
\begin{aligned}
\frac{\dd}{\dd t}\omega_{n} &= \left ( N^{+}_{\mathrm{GOY}}(\omega) \right )_{n} := 
\left ( 2\omega_{n-1}\omega_{n-2}+\frac{1}{2}\omega_{n+1}\omega_{n-1}-\frac{1}{4}
\omega_{n+2}\omega_{n+1}\right )^{\ast} \quad (n \ge 1), \\ 
\omega_{-1} &=0, \quad \omega_{0}=0.
\end{aligned}
\end{equation}
Similarly, for \eqref{eq:main-Sabra}, setting 
$\omega_{n} = i2^{n}U_{n}$ gives
\begin{equation}\label{eq:main-Sabra-renorm}
\begin{aligned}
\frac{\dd}{\dd t}\omega_{n} &= \left ( N^{+}_{\mathrm{Sabra}}(\omega)\right )_{n} 
:= 2\omega_{n-1}\omega_{n-2} + \frac{1}{2}\omega_{n+1}\omega_{n-1}^{\ast} 
-\frac{1}{4}\omega_{n+2}\omega_{n+1}^{\ast} \quad (n \ge 1), \\
\omega_{-1} &=0, \quad \omega_{0}=0.
\end{aligned}
\end{equation}
Both equations share the following real invariant subsystem,
obtained by restricting $\omega$ to real-valued sequences:
\begin{equation}\label{eq:main-renorm}
\begin{aligned}
\frac{\dd}{\dd t}\omega_{n} &= \left ( N^{+}_{\mathrm{real}}(\omega)\right )_{n}
:= 2\omega_{n-1}\omega_{n-2}+\frac{1}{2}\omega_{n+1}\omega_{n-1}-\frac{1}{4}
\omega_{n+2}\omega_{n+1} \quad (n \ge 1), \\
\omega_{-1} &= 0, \quad \omega_{0} =0.
\end{aligned}
\end{equation}

Throughout this introduction, we write
$\ell^{q}=\ell^{q}(\ZZ_{>0}; \CC)$ for the complex equations 
\eqref{eq:main-GOY-renorm} and \eqref{eq:main-Sabra-renorm}, and 
$\ell^{q}=\ell^{q}(\ZZ_{>0}; \RR)$ for the real subsystem
\eqref{eq:main-renorm}. For every $1 \le q \le \infty$, we have \[
\|N^{+}(\omega)\|_{\ell^{q}} \lesssim
\|\omega\|_{\ell^{q}}\|\omega\|_{\ell^{\infty}} \le \|\omega\|_{\ell^{q}}^{2},
\]
where $N^{+}$ denotes any one of the three nonlinearities above.
Thus the nonlinearities are continuous quadratic maps over
$\RR$, and the initial-value problems for
\eqref{eq:main-GOY-renorm}, \eqref{eq:main-Sabra-renorm},
and \eqref{eq:main-renorm} are locally well-posed in the
class of strong solutions, namely functions
$\omega \in C_{t}^{1}(I; \ell^{q})$ satisfying the corresponding
equation in $\ell^{q}$ on a time interval $I$.
Moreover, a strong solution on $[0, T)$, with $0<T<\infty$,
extends beyond $T$ in the same class provided that
\begin{equation}\label{eq:continuation-criterion}
\int_{0}^{T}\|\omega(t)\|_{\ell^{\infty}}\dd t < \infty.
\end{equation}

The main theorem of this paper is the following.
\begin{theorem}\label{thm:main}
For either \eqref{eq:main-GOY-renorm} or
\eqref{eq:main-Sabra-renorm}, there exist finitely supported
initial data $\omega(0)$ such that the strong solution
$\omega \in C_{t}^{1}([0, 1); \ell^{\infty})$ to that equation
satisfies \[
\lim_{t \to 1-}\|\omega(t)\|_{\ell^{\infty}}=\infty.
\]
The blowup is stable in the following sense. There exists
$\epsilon>0$ such that, for every
$\widetilde{\omega}(0)\in\ell^{\infty}(\ZZ_{>0}; \CC)$
satisfying \[
\|\omega(0)-\widetilde{\omega}(0)\|_{\ell^{\infty}}<\epsilon,
\]
the corresponding strong solution $\widetilde{\omega}$
has a finite maximal existence time $T>0$. Moreover, we have the 
blowup rate \[
0 < \liminf_{t\to T-}(T-t)\|\widetilde{\omega}(t)\|_{\ell^{\infty}}
\le \limsup_{t\to T-}(T-t)\|\widetilde{\omega}(t)\|_{\ell^{\infty}} < \infty.
\]
In particular, since 
$\widetilde\omega(0)$ ranges over a complex $\ell^{\infty}$-ball, 
blowup holds for an open set of complex-valued initial data, not 
only within the real subsystem \eqref{eq:main-renorm}.

If additionally $\widetilde{\omega}(0) \in \ell^{q}$ for some
$1 \le q < \infty$, then
$\widetilde{\omega}\in C_{t}^{1}([0, T); \ell^{q})$ and satisfies \[
0 < \liminf_{t\to T-}(T-t)\|\widetilde{\omega}(t)\|_{\ell^{q}}
\le\limsup_{t\to T-}(T-t)\|\widetilde{\omega}(t)\|_{\ell^{q}} < \infty,
\]
so that the same blowup rate holds in $\ell^{q}$.
\end{theorem}

\begin{remark}
One natural space for studying \eqref{eq:main-GOY} and
\eqref{eq:main-Sabra} is the space $H^{s}$, defined by \[
H^{s} := \Big \{ U=(U_{n})_{n \ge 1}: U_{n} \in \CC, \quad \|U\|_{H^{s}}^{2}:=
\sum_{n \ge 1}2^{2sn}|U_{n}|^{2} < \infty \Big \} \quad (s \ge 1).
\]
This Sobolev scale is used in \cite{ConstantinLevantTiti2007}.
The $H^{1}$ norm also appears in numerical studies of finite-time
blowup for the inviscid Sabra model \cite{Mailybaev2013Blowup, CST25}.

Note that strong $H^{1}$ solutions $U$ of either model
correspond to strong $\ell^{2}$ solutions $\omega$ of its
vorticity formulation, with
$\|U\|_{H^{1}}=\|\omega\|_{\ell^{2}}$.
Since the finitely supported initial data $\omega(0)$
in Theorem \ref{thm:main} belong to $\ell^{2}$,
the theorem gives strong $H^{1}$ solutions to both models
that blow up in finite time with the same rate.
The blowup is also stable under small $H^{1}$ perturbations
of the initial data, since $\ell^{2} \hookrightarrow \ell^{\infty}$.

Moreover, for $s \ge 1$, $H^{s}$ regularity is preserved
throughout the existence interval of the corresponding
$\ell^{\infty}$ vorticity solution; see
\cite{ConstantinLevantTiti2007}.
Since $\omega(0)$ is finitely supported, the corresponding
initial velocity $U(0)$ belongs to $H^{s}$ for all $s \ge 1$
in either model. Thus the theorem gives blowup solutions
to both models which are smooth in the sense that 
$U \in C_{t}^{1}([0, 1); H^{s})$ for every $s \ge 1$.
\end{remark}

\subsection{Renormalization} \label{subsec:renorm}

For the proof of Theorem \ref{thm:main}, we use the self-similar
variables described in \cite[Section 2.3]{CST25}. Each of the 
three vorticity equations \eqref{eq:main-GOY-renorm}, 
\eqref{eq:main-Sabra-renorm}, and \eqref{eq:main-renorm} can be written 
as\footnote{We write $\omega=(\omega_{n})_{n \ge 1}$, $N^{+}(\omega)=
((N^{+}(\omega))_{n})_{n \ge 1}$, and use analogous notation for other 
sequences. Equalities between sequences are understood componentwise.}
\begin{equation}\label{eq:main-renorm-common}
\partial_{t}\omega = N^{+}(\omega),
\end{equation}
where $N^{+}$ denotes its nonlinearity. We rescale time and amplitude 
by setting
\begin{equation} \label{eq:ss-variables}
\tau = -\log(1-t), \quad W(\tau)=(1-t)\omega(t).
\end{equation}
In these variables, equation \eqref{eq:main-renorm-common} becomes
\begin{equation} \label{eq:main-ss}
\frac{\dd}{\dd\tau}W = -W + N^{+}(W).
\end{equation}

Under this change of variables, a strong $\ell^{\infty}$ solution of 
\eqref{eq:main-renorm-common} on $[0, 1)$ corresponds to a solution 
$W \in C_{\tau}^{1}([0, \infty); \ell^{\infty})$ satisfying 
\eqref{eq:main-ss}, and the $(1-t)^{-1}$ blowup rate in $\ell^{\infty}$ 
is equivalent to \begin{equation} \label{eq:condition-nontriviality}
0 < \liminf_{\tau \to \infty} \|W(\tau)\|_{\ell^{\infty}}
\le \limsup_{\tau \to \infty} \|W(\tau)\|_{\ell^{\infty}} < \infty.
\end{equation}

\begin{remark} \label{rmk:analogy-with-Euler}
To understand what form we expect $W$ to take, recall the usual
self-similar ansatz for the vorticity of the Euler equations:
\begin{equation} \label{eq:ss-ansatz-Euler}
\omega_{\mathrm{Euler}}(t, x) 
= \frac{1}{1-t}\Omega_{\mathrm{Euler}}\left (\frac{x}{(1-t)^{\gamma}}\right ),
\end{equation}
for some constant $\gamma>0$.
In the analogy with Euler, $\omega_{n}$ plays the role of 
a vorticity amplitude at frequency $2^{n}$. Therefore the amplitude growth 
$(1-t)^{-1}$ in \eqref{eq:ss-ansatz-Euler} corresponds to the 
$(1-t)^{-1}$ blowup rate for the sequence $(\omega_{n})$, and the 
$(1-t)^{\gamma}$ spatial scaling in \eqref{eq:ss-ansatz-Euler} corresponds 
to the shifting of a sequence $(\omega_{n})$ by $-\gamma\log_{2}(1-t)$.
Therefore, a vorticity solution of the Euler equations having 
the self-similar form \eqref{eq:ss-ansatz-Euler} is formally analogous 
to a sequence $(\omega_{n})$ having the form \[
\omega_{n}(t) = (1-t)^{-1}\Omega(n+\gamma\log_{2}(1-t)),
\]
for some time-independent shell profile $\Omega: \RR \to \CC$. 
In the self-similar variables \eqref{eq:ss-variables}, this becomes \[
W_{n}(\tau) = \Omega\left (n-\frac{\gamma\tau}{\log 2} \right ).
\]
Thus the profile moves toward higher shell indices at speed
$\gamma/\log 2$, advancing by one shell over a time interval
of length $(\log 2)/\gamma$.
This correspondence between self-similar blowup and travelling waves 
is discussed for the Sabra model on $\ZZ$ in \cite[Section 2.3]{CST25}.
\end{remark}

\subsection{Proof sketch}

In this subsection we outline the proof of Theorem \ref{thm:main}.

The analogy in Remark \ref{rmk:analogy-with-Euler} suggests finding 
a travelling wave solution for \eqref{eq:main-ss}. Note that in the 
equation \eqref{eq:main-ss}, the shell indices are restricted
to $\ZZ_{>0}$.
The main idea of our proof is to first construct a travelling 
wave on all of $\ZZ$ with a decaying left tail, rather than
seek one whose left tail is identically zero.
We then use this wave to construct a nearby solution on
the positive shells satisfying the original boundary conditions.

To this end, we first consider the same evolution law as in
\eqref{eq:main-ss}, but for all shell indices $n\in\ZZ$
and without boundary conditions:
\begin{equation} \label{eq:main-ss-Z}
\partial_{\tau}Z_{n} = -Z_{n} + \left (N(Z)\right )_{n} \quad (n \in \ZZ),
\end{equation}
where, as before, the nonlinearities are defined by
\begin{equation*}
\begin{aligned}
\left (N_{\mathrm{GOY}}(Z)\right )_{n} &:= 
\left ( 2Z_{n-1}Z_{n-2} + \frac{1}{2}Z_{n+1}Z_{n-1} 
-\frac{1}{4}Z_{n+2}Z_{n+1}\right )^{\ast}, \\
\left (N_{\mathrm{Sabra}}(Z)\right )_{n} &:= 
2Z_{n-1}Z_{n-2} + \frac{1}{2}Z_{n+1}Z_{n-1}^{\ast}
-\frac{1}{4}Z_{n+2}Z_{n+1}^{\ast}, \\
\left (N_{\mathrm{real}}(Z)\right )_{n} &:= 
2Z_{n-1}Z_{n-2}+\frac{1}{2}Z_{n+1}Z_{n-1} 
-\frac{1}{4}Z_{n+2}Z_{n+1}.
\end{aligned}
\end{equation*}
Since the GOY and Sabra nonlinearities agree with $N_{\mathrm{real}}$
on real-valued sequences, it is enough to construct a real-valued 
travelling wave solving \eqref{eq:main-ss-Z} with $N=N_{\mathrm{real}}$.
Its linearized dynamics, however, must be studied under complex 
perturbations, where the two models differ through their phase 
symmetries; this is done separately for each model in Section 
\ref{sec:monodromy}.

We seek a nonzero, spatially decaying solution 
$Z^{(\ast)}$ and a time $\delta^{(\ast)}>0$ such that
\begin{equation} \label{eq:orbit-shift-relation}
Z^{(\ast)}(\tau+\delta^{(\ast)})=RZ^{(\ast)}(\tau) \quad (\tau\ge 0),
\end{equation}
where $R$ is the one-shell right shift defined by $(RZ)_{n}=Z_{n-1}$.
The relation \eqref{eq:orbit-shift-relation} describes a bilateral 
travelling wave with relative period $\delta^{(\ast)}$: over each 
such time interval, the profile advances one shell to the right.
In Section \ref{sec:bilateral}, we construct such a wave together 
with its relative period. Following ideas of Mailybaev
\cite{Mailybaev2012, Mailybaev2013Bifurcations} and
Campolina, Simonnet, and Thalabard \cite{CST25}, we compute an
approximation to $(Z^{(\ast)}, \delta^{(\ast)})$ by tracking 
the weighted center of a numerically evolving profile; see Appendix 
\ref{app:numerical-approximation}.
To construct an exact solution near this numerical approximation, we use 
a Banach fixed-point argument, reducing the existence proof to a 
collection of quantitative estimates. These estimates initially 
involve infinite sums. Analytic bounds on the tails reduce their 
verification to finite computations, which are carried out using 
rigorous interval arithmetic.

At this stage in the proof, we have constructed a travelling wave 
solution to the bilateral analogue \eqref{eq:main-ss-Z} instead of 
\eqref{eq:main-ss}. Simply restricting the 
solution to positive shells does not yield a solution to \eqref{eq:main-ss}, 
as the shells with nonpositive indices affect those with positive indices
through the nonlinearity. So instead of simply truncating the sequence, 
we add a boundary correction $B$ which cancels this influence of 
nonpositive shells on positive shells. Define the boundary 
correction $B$ by
\begin{equation*}
\begin{aligned}
(B_{\mathrm{GOY}}(Z))_{1} = -2Z_{0}^{\ast}Z_{-1}^{\ast} 
-\frac{1}{2} Z_{2}^{\ast}Z_{0}^{\ast}, \quad 
(B_{\mathrm{GOY}}(Z))_{2} = -2Z_{1}^{\ast}Z_{0}^{\ast}, 
\quad (B_{\mathrm{GOY}}(Z))_{n} =0 \quad (n \neq 1, 2), \\
(B_{\mathrm{Sabra}}(Z))_{1} = -2Z_{0}Z_{-1} 
-\frac{1}{2} Z_{2}Z_{0}^{\ast}, \quad 
(B_{\mathrm{Sabra}}(Z))_{2} = -2Z_{1}Z_{0}, 
\quad (B_{\mathrm{Sabra}}(Z))_{n} =0 \quad (n \neq 1, 2). \\
\end{aligned}
\end{equation*}
If $\widehat{Z}$ solves the boundary-corrected equation
\begin{equation} \label{eq:main-ss-Z-modified}
\partial_{\tau}\widehat{Z} = -\widehat{Z} + N(\widehat{Z})
+ B(\widehat{Z}),
\end{equation}
then restricting $\widehat{Z}$ to the positive indices gives a 
solution to \eqref{eq:main-ss}. We therefore look for a
solution $\widehat{Z}$ of \eqref{eq:main-ss-Z-modified}.

Recall that $Z^{(\ast)}$ advances by one shell over each relative 
period $\delta^{(\ast)}$. As $Z^{(\ast)}$ decays spatially, 
the boundary correction $B(Z^{(\ast)}(\tau))$ tends to zero as the self-similar 
time $\tau$ goes to infinity. So for large $J$, after translating the 
travelling wave $J$ shells to the right, $R^{J}Z^{(\ast)}$ has small 
boundary correction and thus approximately solves 
\eqref{eq:main-ss-Z-modified}. This is why it is reasonable to search 
for a solution $\widehat{Z}$ to \eqref{eq:main-ss-Z-modified} near the 
travelling wave solution $R^{J}Z^{(\ast)}$ of \eqref{eq:main-ss-Z}.

A first approach is to look for a solution $\widehat{Z}$ that approaches 
$R^{J}Z^{(\ast)}$. Since $R^{J}Z^{(\ast)}$ is a travelling wave,
after suitable recentering, we can compare the recentered profiles 
of $\widehat{Z}$ with a fixed profile $Z^{(\ast)}(0)$. To be more 
precise, write 
$\widehat{Z}((j-J)\delta^{(\ast)}) = R^{j}\widehat{Z}^{[j]}$ for 
$j \ge J$. Since 
$R^{J}Z^{(\ast)}((j-J)\delta^{(\ast)})=R^{j}Z^{(\ast)}(0)$, we want 
$\widehat{Z}^{[j]}$ to approach $Z^{(\ast)}(0)$.

Note that the map from 
$\widehat{Z}^{[j]}$ to $\widehat{Z}^{[j+1]}$ is given by evolving the 
boundary-corrected equation for one relative period $\delta^{(\ast)}$ 
and then shifting one shell to the left. (In these recentered coordinates, 
the boundary correction acts at indices $1-j$ and $2-j$.) As we discussed 
earlier, the boundary correction is small near the travelling wave 
$R^{J}Z^{(\ast)}$, so this map from $\widehat{Z}^{[j]}$ to 
$\widehat{Z}^{[j+1]}$ can be approximated by the bilateral relative 
return map \[
F := R^{-1}\Phi^{\delta^{(\ast)}},
\]
where $\Phi^{s}$ denotes the evolution for self-similar time $s$ under 
\eqref{eq:main-ss-Z}. Note that the travelling-wave relation gives
$F(Z^{(\ast)}(0))=Z^{(\ast)}(0)$. Thus, writing 
$\widehat{\xi}^{[j]}=\widehat{Z}^{[j]}-Z^{(\ast)}(0)$ and neglecting the 
boundary correction, we obtain the formal linearized iteration 
$\widehat{\xi}^{[j+1]}=M\widehat{\xi}^{[j]}$, where the linearized bilateral relative 
return map $M$ is defined by \[
M := DF(Z^{(\ast)}(0)) = R^{-1}D\Phi^{\delta^{(\ast)}}(Z^{(\ast)}(0)).
\]

However, we do not expect $\{\widehat{Z}^{[j]}\}_{j \ge J}$ to 
converge to $Z^{(\ast)}(0)$ for all nearby initial profiles 
$R^{J}\widehat{Z}^{[J]}$, as $M$ has nonstable modes coming from 
time-translation symmetries and phase symmetries inherent to the 
equation \eqref{eq:main-ss-Z}. One could try to prove the stability 
statement in Theorem \ref{thm:main} by constructing solutions 
approaching $R^{J}Z^{(\ast)}$ and then enlarging this family 
using the admissible symmetries. An inverse function theorem argument 
could then show that the resulting set of initial data contains an 
open ball. For the GOY model, the relevant symmetry transformations
depend smoothly on their parameters in $\ell^{\infty}(\ZZ; \CC)$, 
so this approach is available.

But in the Sabra model, phase symmetry transformations do not 
define a $C^{1}$ action on the space $\ell^{\infty}(\ZZ; \CC)$, 
which is an obstacle to this inverse function theorem argument. The 
inverse function theorem approach would choose a single transformation
of the initial data so that the resulting solution approaches the
reference wave $R^{J}Z^{(\ast)}$. We instead use a modulation 
argument, where we choose the parameters step by step, by imposing 
conditions on finitely many shell coordinates at each step. On these 
coordinates the phase action is smooth, so determining the parameters 
does not require a $C^{1}$ phase action on the full 
$\ell^{\infty}(\ZZ; \CC)$ space. The use of modulation in stability 
arguments is classical \cite{Weinstein1985, MerleRaphael2005}. For 
applications to self-similar singularities in fluid equations, see 
\cite{ElgindiGhoulMasmoudi2021Euler, 
BuckmasterShkollerVicol2023Point}.

We briefly explain the modulation argument used in our proof.
Let $V$ be the span of the nonstable modes of $M$ coming from 
symmetries and choose a complement $Y$. Then, after each iteration 
of $F$, we can apply a suitable symmetry transformation to bring the 
profile into $Z^{(\ast)}(0)+Y$. Formally, in terms of the difference
$\widehat{\xi}^{[j]}=\widehat{Z}^{[j]}-Z^{(\ast)}(0)$, we first apply the 
linearized bilateral relative return map $M$, and then apply $I-\Pi$, 
where $\Pi$ is a projection onto $V$ with kernel $Y$. This is why we 
study the modulated operator \[
\widetilde{M} = (I-\Pi)M|_{Y} : Y \to Y
\]
in Section \ref{sec:monodromy}, where we prove that a power of 
$\widetilde{M}$ is 
contractive. Using this contraction estimate, we show that the 
modulated iterates converge to $Z^{(\ast)}(0)$. We then undo the 
symmetry transformations and return to the original variables to
obtain solution segments, which we glue together to form a blowup 
solution. This construction works for all initial data in an open 
ball. More details can be found in Section \ref{sec:boundary}.

Interval arithmetic is used at exactly two points in 
the proof: in Subsection \ref{subsec:orbit-finite-verification}, to 
verify the bounds guaranteeing an exact travelling wave $Z^{(\ast)}$ 
near the numerical approximation, and in Subsection 
\ref{subsec:monodromy-verification}, to verify the bounds showing 
that a power of $\widetilde{M}$ is contractive. The numerical 
approximation of Appendix \ref{app:numerical-approximation} need not 
itself be rigorous; it only serves as the center of the fixed-point 
argument. All remaining steps are analytic. The verification uses 
Arb through python-flint \cite{pythonflint, Johansson2017arb}; all 
code is available in the supplementary material and in the 
GitHub repository \cite{CAParchive}.

This use of interval arithmetic is well established in the study of
differential equations. Examples include existence proofs for elliptic
equations \cite{Nakao1988, Plum1992}, blowup results for the
harmonic map heat flow
\cite{BiernatDonninger2018, BiernatDonningerSchorkhuber2017}
and the De Gregorio model \cite{ChenHouHuang2021}, and stable nearly 
self-similar blowup for 2D Boussinesq and 3D axisymmetric Euler 
equations with boundary \cite{ChenHou2022, ChenHou2025}. For a survey 
of computer-assisted proofs in PDE, see \cite{GomezSerrano2019}.

\subsection*{Acknowledgments}

I am grateful to Vlad Vicol for introducing me to this problem, 
which he had previously discussed with Theodore Drivas. I also thank 
Giorgio Cialdea and Gerard Castro-López for helpful discussions. I 
acknowledge support from the Simons Foundation through a Simons 
Investigator Award to Vlad Vicol.

\subsection*{Use of AI systems}

Claude Fable 5.1 and GPT-6 Astra were used in the development of the 
proofs, the preparation of the manuscript, and the development of the 
accompanying code. Following Vlad Vicol's suggestion to seek a fixed 
point for the bilateral problem, Claude Fable 5.1 produced an initial proof 
draft. While studying and revising this draft, I considered the 
stability of the blowup solution and the extension of the stability 
result to complex-valued perturbations. At my request, GPT-6 Astra 
produced draft proofs of these extensions.

I subsequently revised the mathematical framework and the arguments. 
In particular, I reformulated the main result in $\ell^{\infty}$ 
rather than $\ell^{2}$, replaced an earlier Lyapunov-Perron argument 
with modulation, and used an unweighted $\ell^{1}$ space for the 
travelling-wave construction, with spatial decay of the travelling 
wave established separately. I drafted and revised the manuscript 
with assistance from GPT-6 Astra in the English exposition and in 
literature searches. GPT-6 Astra also generated code, which I 
reviewed, checked, and partially rewrote. I take full responsibility 
for the manuscript and the accompanying code.

\section{Existence of a bilateral travelling wave}
\label{sec:bilateral}

Consider the bilateral versions of the renormalized equation 
\eqref{eq:main-ss} given in the introduction, \eqref{eq:main-ss-Z}. The 
goal of this section is to construct a spatially decaying travelling wave for 
these equations. Since the bilateral renormalized GOY and Sabra 
equations both reduce on real-valued sequences to the real 
subsystem \eqref{eq:main-ss-Z} with $N=N_{\mathrm{real}}$, we focus on 
constructing a travelling wave for this real subsystem.
Throughout this section, we write $N = N_{\mathrm{real}}$.

It will be convenient to define the symmetric bilinear map $Q$ by \[
\left (Q(Z, \widetilde{Z})\right )_{n} := Z_{n-1}\widetilde{Z}_{n-2}
+\widetilde{Z}_{n-1}Z_{n-2}+\frac{1}{4}Z_{n+1}\widetilde{Z}_{n-1}
+\frac{1}{4}\widetilde{Z}_{n+1}Z_{n-1} 
-\frac{1}{8}Z_{n+2}\widetilde{Z}_{n+1}
-\frac{1}{8}\widetilde{Z}_{n+2}Z_{n+1},
\] so that $N(Z)=Q(Z, Z)$.

We will find a travelling wave solution to \eqref{eq:main-ss-Z} in 
$\ell^{1}(\ZZ; \RR)$. The choice of this function space is to require 
decay to avoid trivial solutions like $Z_{n}=\frac{4}{9}$ for all $n \in \ZZ$, 
and to simplify the computation of induced operator norms. The bilateral 
travelling wave we find will turn out to have exponential spatial decay, 
which will be proved in Subsection \ref{subsec:orbit-spatial-decay}.

Note that we have 
\begin{equation} \label{eq:orbit-bilinear-bound}
\|Q(Z, \widetilde{Z})\|_{\ell^{1}(\ZZ; \RR)} 
\le C_{Q} \|Z\|_{\ell^{1}(\ZZ; \RR)}
\|\widetilde{Z}\|_{\ell^{\infty}(\ZZ; \RR)} 
\le C_{Q} \|Z\|_{\ell^{1}(\ZZ; \RR)}
\|\widetilde{Z}\|_{\ell^{1}(\ZZ; \RR)},
\end{equation}
where $C_{Q} = \frac{11}{4}$. In particular $N$ is locally Lipschitz 
on $\ell^{1}(\ZZ; \RR)$, so that we have local existence and uniqueness for 
\eqref{eq:main-ss-Z} on $\ell^{1}(\ZZ; \RR)$.

Let $R$ denote the one-shell right shift, \[
(RZ)_{n} = Z_{n-1} \quad (n \in \ZZ).
\]
The goal of this section is the following.

\begin{theorem} \label{thm:orbit-existence}
There exist $\delta^{(\ast)} > 0$ and a nonzero solution
$Z^{(\ast)} \in C_{\tau}^{1}([0, \infty); \ell^{1}(\ZZ; \RR))$ 
to \eqref{eq:main-ss-Z} such that 
$Z^{(\ast)}(\tau+\delta^{(\ast)})=RZ^{(\ast)}(\tau)$ for all $\tau \ge 0$.

Moreover, this $Z^{(\ast)}$ satisfies the following spatial decay. For 
every $1 < \nu_{-} < e^{\delta^{(\ast)}}$ and $\nu_{+}>1$, define \[
\ell_{\nu_{-}, \nu_{+}}^{1} := \left \{ z=(z_{n})_{n \in \ZZ}: 
\|z\|_{\ell_{\nu_{-}, \nu_{+}}^{1}} := 
\sum_{n \ge 0} \nu_{+}^{|n|}|z_{n}| + \sum_{n < 0} \nu_{-}^{|n|}|z_{n}| 
< \infty \right \}.
\] Then, we have $Z^{(\ast)} \in C_{\tau}^{1}([0, \infty); 
\ell_{\nu_{-}, \nu_{+}}^{1})$.
\end{theorem}

\begin{figure}[htbp]
\centering
\includegraphics[width=0.76\linewidth]{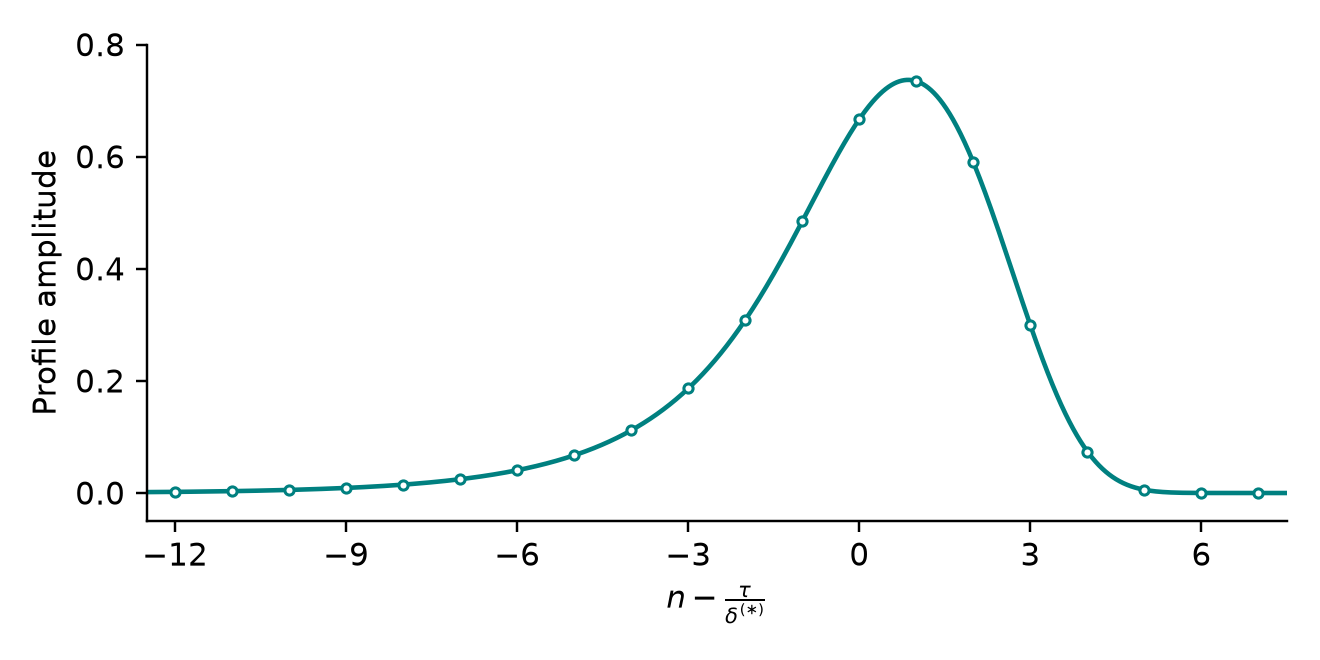}
\caption{Approximate values of $(Z^{(\ast)}(\tau))_{n}$ plotted against
$n-\frac{\tau}{\delta^{(\ast)}}$. The wave advances one shell to the right 
per relative period $\delta^{(\ast)}$, so its profile is 
stationary in this coordinate. The circles mark the values at 
$\tau=0$.}
\label{fig:bilateral-profile}
\end{figure}

Figure \ref{fig:bilateral-profile} shows the numerical approximation
used to construct this travelling wave. The plotting script 
\texttt{plot\_orbit.py} and the data used to generate the figure 
are available in the supplementary material and in the GitHub 
repository \cite{CAParchive}. Appendix
\ref{app:numerical-approximation} describes how we generated the 
data.

The remainder of this section is devoted to the proof of Theorem
\ref{thm:orbit-existence}. We first prove that any bilateral travelling 
wave in $\ell^{1}(\ZZ; \RR)$ has the required decay. To construct a travelling 
wave, we begin by formulating the problem as an infinite-dimensional 
nonlinear equation. Given an approximate solution to this equation, we set 
up a Banach fixed-point argument to find a true solution nearby. This 
reduces the problem to verifying explicit bounds. We control the 
infinite-dimensional tails in these bounds analytically, leaving only 
finite computations. Finally, we specify the approximate solution and 
certify the remaining finite estimates using interval arithmetic.

Throughout this section, $\ell^{1}$ denotes $\ell^{1}(\ZZ; \RR)$.

\subsection{Spatial decay of the travelling wave}
\label{subsec:orbit-spatial-decay}

We first assume the existence of a travelling wave and prove the
required decay in Theorem \ref{thm:orbit-existence}.
Let $Z^{(\ast)}\in C_{\tau}^{1}([0, \infty); \ell^{1})$ be a
solution of \eqref{eq:main-ss-Z} satisfying
\begin{equation}\label{eq:orbit-decay-shift}
Z^{(\ast)}(\tau+\delta^{(\ast)}) = RZ^{(\ast)}(\tau),
\end{equation}
where $\delta^{(\ast)}>0$. Define \[
R_{j} := \sup_{0 \le \tau \le \delta^{(\ast)}} \sup_{n \ge j} 
|(Z^{(\ast)}(\tau))_{n}|, \quad 
L_{j} := \sup_{0 \le \tau \le \delta^{(\ast)}} \sup_{n \le -j}
|(Z^{(\ast)}(\tau))_{n}| \quad (j \ge 0).
\]

We claim that both $R_{j}$ and $L_{j}$ tend to zero as $j \to \infty$.
Take any $\epsilon>0$. Since 
$Z^{(\ast)}\in C([0, \delta^{(\ast)}]; \ell^{1})$, there exist 
finitely many times $\tau_{1}, \ldots, \tau_{k}$ such that for 
any $0 \le \tau \le \delta^{(\ast)}$ there exists some $\tau_{i}$
such that 
\begin{equation} \label{eq:orbit-time-cover}
\|Z^{(\ast)}(\tau)-Z^{(\ast)}(\tau_{i})\|_{\ell^{1}} < \frac{1}{2}\epsilon.
\end{equation}
Now choose $J$ large enough so that \[
\sum_{|n| \ge J} |(Z^{(\ast)}(\tau_{i}))_{n}| < \frac{1}{2}\epsilon 
\quad (i=1, 2, \ldots, k).
\]
Then for any $0 \le \tau \le \delta^{(\ast)}$ and $j \ge J$, we have \[
\sum_{|n| \ge j} |(Z^{(\ast)}(\tau))_{n}|  
\le \sum_{|n| \ge J} |(Z^{(\ast)}(\tau))_{n}| 
\le \sum_{|n| \ge J}|(Z^{(\ast)}(\tau_{i}))_{n}| +
\|Z^{(\ast)}(\tau)-Z^{(\ast)}(\tau_{i})\|_{\ell^{1}} 
< \epsilon,
\]
for some $\tau_{i}$ satisfying \eqref{eq:orbit-time-cover}. Taking the supremum over 
$\tau$ gives $\max\{R_{j}, L_{j}\} \le \epsilon$ for every $j \ge J$.
As $\epsilon>0$ was arbitrary, this proves that both $R_{j}$ and 
$L_{j}$ tend to zero as $j \to \infty$.

For $0 \le \tau \le \delta^{(\ast)}$, \eqref{eq:orbit-decay-shift} and Duhamel's formula 
applied to $\partial_{\tau}Z^{(\ast)} = -Z^{(\ast)} + N(Z^{(\ast)})$ give
\begin{equation*}
\begin{aligned}
& (Z^{(\ast)}(\tau))_{n} \\
&= (Z^{(\ast)}(\tau+\delta^{(\ast)}))_{n+1} \\
&= e^{-\delta^{(\ast)}} (Z^{(\ast)}(\tau))_{n+1} 
+ \int_{\tau}^{\tau+\delta^{(\ast)}} e^{-(\tau+\delta^{(\ast)}-s)}
(N(Z^{(\ast)}(s)))_{n+1}\dd s \\
&= e^{-\delta^{(\ast)}} (Z^{(\ast)}(\tau))_{n+1} 
+ \int_{\tau}^{\delta^{(\ast)}} e^{-(\tau+\delta^{(\ast)}-s)}
(N(Z^{(\ast)}(s)))_{n+1}\dd s
+ \int_{0}^{\tau} e^{-(\tau-s)}(N(Z^{(\ast)}(s)))_{n}\dd s.
\end{aligned}
\end{equation*}

Therefore, the definition of $N$ gives the bound
\begin{align}
R_{j} &\le e^{-\delta^{(\ast)}}R_{j+1} + \sup_{0\le\tau\le\delta^{(\ast)}}
\left [ \int_{\tau}^{\delta^{(\ast)}} e^{-(\tau+\delta^{(\ast)}-s)} 
C_{Q}R_{j-1}^{2} \dd s    
+ \int_{0}^{\tau} e^{-(\tau-s)}C_{Q}R_{j-2}^{2}\dd s \right ]
\nonumber \\
& \le e^{-\delta^{(\ast)}}R_{j+1} + C_{Q}
(1-e^{-\delta^{(\ast)}})R_{j-2}^{2} \quad (j \ge 3),
\label{eq:orbit-right-duhamel-bound}
\end{align}
where $C_{Q}=\frac{11}{4}$, and similarly
\begin{equation} \label{eq:orbit-left-tail-recurrence}
L_{j} \le e^{-\delta^{(\ast)}}L_{j-1} + 
C_{Q}(1-e^{-\delta^{(\ast)}})L_{j-3}^{2} \quad (j \ge 4).
\end{equation}

To prove the right decay, note that using $R_{j+1} \le R_{j}$ 
in \eqref{eq:orbit-right-duhamel-bound} gives
\begin{equation} \label{eq:orbit-right-tail-recurrence}
R_{j} \le C_{Q}R_{j-2}^{2} \quad (j \ge 3).
\end{equation}
Iterating this gives \[
C_{Q} R_{j+2k} \le \left ( C_{Q} R_{j}\right )^{2^{k}}
\quad (j \ge 3, k \ge 0).
\]
Therefore, if we choose $J$ large enough so that $C_{Q}R_{J} < \frac{1}{2}$,
then we have \[
R_{J+2k} \le 2^{-2^{k}} \quad (k \ge 0).
\]
This gives the desired right decay.

To prove the left decay, let $1 < \nu_{-}
< \nu_{-}+\eta < e^{\delta^{(\ast)}}$.
Choose $J$ large enough so that \[
e^{-\delta^{(\ast)}}(\nu_{-}+\eta) 
+ C_{Q}L_{J}(\nu_{-}+\eta)^{3} < 1.
\]
Then by induction on $j$ with \eqref{eq:orbit-left-tail-recurrence}, we obtain \[
L_{j} \le (\nu_{-}+\eta)^{J+2-j}L_{J} \quad (j \ge J).
\]
This gives the desired left decay.

Thus, for any fixed pair $(\nu_{-}, \nu_{+})$ in Theorem 
\ref{thm:orbit-existence}, the tails of the series defining 
$\|Z^{(\ast)}(\tau)\|_{\ell_{\nu_{-}, \nu_{+}}^{1}}$ tend to zero 
uniformly on $\tau \in [0, \delta^{(\ast)}]$. Since each component 
is continuous, we obtain
$Z^{(\ast)} \in C([0, \delta^{(\ast)}]; \ell_{\nu_{-}, \nu_{+}}^{1})$,
and by \eqref{eq:orbit-decay-shift} and the equation $\partial_{\tau}Z^{(\ast)} 
= -Z^{(\ast)} + N(Z^{(\ast)})$, we conclude that 
$Z^{(\ast)} \in C_{\tau}^{1}([0, \infty); \ell_{\nu_{-}, \nu_{+}}^{1})$.

Hence to prove Theorem \ref{thm:orbit-existence}, it remains to prove 
the existence of a bilateral travelling wave $Z^{(\ast)}$.

\subsection{Formulation as a nonlinear equation}
\label{subsec:nonlinear-formulation}
The goal is to find 
$Z^{(\ast)} \in C_{\tau}^{1}([0, \infty); \ell^{1}(\ZZ))$ solving
\eqref{eq:main-ss-Z} and a relative period $\delta^{(\ast)}>0$ such that 
$Z^{(\ast)}(\delta^{(\ast)})=RZ^{(\ast)}(0)$. We seek a solution over 
one relative period using the Taylor series ansatz
\begin{equation}\label{eq:bilateral-orbit-ansatz}
Z(\delta\theta) = \sum_{k \ge 0}a^{(k)}\theta^{k} 
\quad (0 \le \theta \le 1).
\end{equation}
Here, each $a^{(k)}=(a^{(k)}_{n})_{n \in \ZZ}$ is a sequence, and both
$\delta>0$ and the coefficient sequence $a=(a^{(k)})_{k \ge 0}$ are to be
determined. We will search for $a$ in the following space: \[
\ell_{\rho}^{1}:=
\ell_{\rho}^{1}(\ZZ_{\ge 0}; \ell^{1}(\ZZ; \RR))
:= \left \{ a=(a^{(k)})_{k \ge 0}: \|a\|_{\ell_{\rho}^{1}} := \sum_{k \ge 0}\rho^{k}\|a^{(k)}\|_{\ell^{1}} < \infty\right \},
\] where $\rho > 1$ will be fixed in the final subsection.

Under the ansatz \eqref{eq:bilateral-orbit-ansatz}, the relative 
periodicity condition and the equation \eqref{eq:main-ss-Z} become 
\begin{align}
\sum_{k \ge 0}a^{(k)} &= Ra^{(0)}, 
\label{eq:orbit-nonlinear-matching-condition} \\
(k+1)a^{(k+1)} &= \delta \left [ - a^{(k)} + \sum_{j=0}^{k} 
Q(a^{(j)}, a^{(k-j)}) \right ] \quad (k \ge 0). 
\label{eq:orbit-nonlinear-Taylor-condition}
\end{align}

Let $(\overline{\delta}, \overline{a}) \in \RR \times \ell_{\rho}^{1}$ 
be an approximate solution to
the equations above, to be specified in the final subsection. We seek 
a true solution near $(\overline{\delta}, \overline{a})$. Because of 
time-translation invariance, there is a one-dimensional phase freedom,
which we fix by imposing the temporal phase condition \[
\sum_{n \in \ZZ}\overline{q}_{n} \left ( a^{(0)}_{n} 
-\overline{a}^{(0)}_{n} \right )=0.
\]
Here, $\overline{q}$ is a finitely supported approximate tangent vector
to the numerical orbit, also to be specified later.

We now rewrite the equations above and the temporal phase condition as a 
single nonlinear equation. For 
$(\delta, a) \in \RR \times \ell_{\rho}^{1}$, define 
\begin{equation*}
\begin{aligned}
\phi(a^{(0)}) &:= \sum_{n \in \ZZ}\overline{q}_{n}(a_{n}^{(0)}-\overline{a}_{n}^{(0)}), \\
F^{(0)}(\delta, a) &:= R^{-1}\sum_{k \ge 0}a^{(k)} - a^{(0)}, \\
F^{(k)}(\delta, a) &:= a^{(k)} - \frac{\delta}{k}\left ( -a^{(k-1)} 
+ \sum_{j=0}^{k-1} Q(a^{(j)}, a^{(k-1-j)}) \right ) \quad (k \ge 1),
\end{aligned}
\end{equation*}
and define 
$\cF: \RR \times \ell_{\rho}^{1} \to \RR \times \ell_{\rho}^{1}$ 
by \[
\cF(\delta, a) = \left (\phi(a^{(0)}), F(\delta, a)\right ) = \left (\phi(a^{(0)}), (F^{(k)}(\delta, a))_{k \ge 0}\right ).
\]

Thus, Theorem \ref{thm:orbit-existence} reduces to finding a solution of 
$\cF(\delta, a)=0$ with $\delta>0$ and $a \not\equiv 0$.

\subsection{Linear preconditioning}

We precondition $F$ by replacing it with $L_{0}^{-1}F$ for a suitable
invertible linear operator $L_{0}$. The purpose of this is to make
$L_{0}^{-1}D_{a}F$ close to the identity on the infinite-dimensional
tails, which simplifies the estimates in the later fixed-point argument.

To choose such $L_{0}$, observe that for $h \in \ell_{\rho}^{1}$, the
contribution of the quadratic nonlinearity to
$(D_{a}F(\delta, a)h)_{n}^{(k)}$ is \[
-\frac{2\delta}{k} \sum_{j=0}^{k-1}\left (Q(a^{(j)}, h^{(k-1-j)})
\right )_{n} \quad (k\ge 1).
\]
The summability of $a$ makes these contributions small for large $|n|$, 
while the factor $\frac{1}{k}$ makes them small for large $k$.
Thus, on the infinite-dimensional tails, the contribution of
the nonlinear terms to $D_{a}F$ is small, and the linear part of $F$
captures the dominant behavior. We therefore define $L_{0}$ to be the 
derivative of the linear part of $F(\overline{\delta}, \cdot)$, i.e., \[
L_{0}:= D_{a}F(\overline{\delta}, 0).
\]
In particular, for $h \in \ell_{\rho}^{1}$,
\begin{align}
(L_{0}h)^{(0)} &= R^{-1}\sum_{k\ge 0} h^{(k)} - h^{(0)}, 
\label{eq:L0-formula-zeroth} \\
(L_{0}h)^{(k)} &= h^{(k)}+\frac{\overline{\delta}}{k}h^{(k-1)}
\quad (k\ge 1). \label{eq:L0-formula-high}
\end{align}

\begin{lemma}[Invertibility of $L_{0}$] \label{lem:L0-inverse}
Assume $\rho>1$ and $\overline{\delta}>0$.
Then $L_{0}$ is a bounded isomorphism of $\ell_{\rho}^{1}$, and
\begin{equation} \label{eq:L0-inverse-bound}
\|L_{0}^{-1}\|_{\ell_{\rho}^{1} \to \ell_{\rho}^{1}} \le
\frac{e^{(1+\rho)\overline{\delta}}}{1-e^{-\overline{\delta}}}
+e^{\rho\overline{\delta}}.
\end{equation}
Moreover, for every $y \in \ell_{\rho}^{1}$, $h=L_{0}^{-1}y$ is given by 
\eqref{eq:L0-explicit-inverse}.
\end{lemma}

\begin{proof}
The boundedness of $L_{0}$ follows from its definition and
$\|R^{-1}\|_{\ell^{1} \to \ell^{1}} = 1$.
To construct its inverse, fix $y\in \ell_{\rho}^{1}$ and consider
the equation $L_{0}h=y$.

The equations \eqref{eq:L0-formula-high} give \[
h^{(k)}=y^{(k)}-\frac{\overline{\delta}}{k}h^{(k-1)} \quad (k\ge 1).
\]
Iterating this recurrence gives \[
h^{(k)} = \sum_{\ell=1}^{k}\frac{\ell!}{k!}(-\overline{\delta})^{k-\ell}
y^{(\ell)} + \frac{1}{k!}(-\overline{\delta})^{k}h^{(0)} \quad (k \ge 1).
\] Substituting these into the zeroth equation 
\eqref{eq:L0-formula-zeroth} gives \[
y^{(0)} = R^{-1}\sum_{k \ge 1}\sum_{\ell=1}^{k}
\frac{\ell!}{k!}(-\overline{\delta})^{k-\ell}y^{(\ell)}
+ R^{-1} e^{-\overline{\delta}}h^{(0)} - h^{(0)}.
\]
To solve for $h^{(0)}$, we need to invert 
$I-e^{-\overline{\delta}}R^{-1}$, which is possible since \[
\|e^{-\overline{\delta}}R^{-1}\|_{\ell^{1} 
\to \ell^{1}} = e^{-\overline{\delta}} < 1.
\] In particular, by the Neumann series, \[
\left \|(I-e^{-\overline{\delta}}R^{-1})^{-1}\right \|_{
\ell^{1} \to \ell^{1}}
\le \frac{1}{1-e^{-\overline{\delta}}}.
\] This leads to the following candidate inverse $y\mapsto h$:
\begin{equation} \label{eq:L0-explicit-inverse}
\begin{aligned}
h^{(0)} &= (I-e^{-\overline{\delta}}R^{-1})^{-1} \left (
R^{-1}\sum_{k\ge 1}\sum_{\ell=1}^{k}\frac{\ell!}{k!}
(-\overline{\delta})^{k-\ell}y^{(\ell)} -y^{(0)} \right ), \\
h^{(k)} &= \frac{(-\overline{\delta})^{k}}{k!}h^{(0)} + 
\sum_{\ell=1}^{k}\frac{\ell!}{k!}(-\overline{\delta})^{k-\ell}y^{(\ell)}
\quad (k\ge 1).
\end{aligned}
\end{equation}
We now show that these formulas define a bounded linear map on
$\ell_{\rho}^{1}$ and verify that $L_{0}h=y$.

First, using $\frac{\ell!}{k!} \le \frac{1}{(k-\ell)!}$, we obtain
\begin{equation*}
\begin{aligned}
\|h^{(0)}\|_{\ell^{1}} &\le 
\left \|(I-e^{-\overline{\delta}}R^{-1})^{-1}\right \|_{\ell^{1} \to \ell^{1}}
\left ( \|R^{-1}\|_{\ell^{1} \to \ell^{1}} \sum_{k \ge 1}\sum_{\ell=1}^{k}\frac{\ell!}{k!}
\overline{\delta}^{k-\ell} \|y^{(\ell)}\|_{\ell^{1}} 
+ \|y^{(0)}\|_{\ell^{1}}\right ) \\
&\le \frac{1}{1-e^{-\overline{\delta}}} 
\left ( \sum_{\ell \ge 1}\sum_{k \ge \ell}\frac{1}{(k-\ell)!}
\overline{\delta}^{k-\ell} \|y^{(\ell)}\|_{\ell^{1}} 
+ \|y^{(0)}\|_{\ell^{1}}\right ) \\
&= \frac{1}{1-e^{-\overline{\delta}}} \left (   
e^{\overline{\delta}} \sum_{\ell \ge 1} 
\|y^{(\ell)}\|_{\ell^{1}} + \|y^{(0)}\|_{\ell^{1}}\right ) \\
&\le \frac{e^{\overline{\delta}}}{1-e^{-\overline{\delta}}} 
\|y\|_{\ell_{\rho}^{1}}, \\
\|h\|_{\ell_{\rho}^{1}} &\le \|h^{(0)}\|_{\ell^{1}} 
+ \sum_{k \ge 1} \rho^{k} \frac{\overline{\delta}^{k}}{k!} 
\|h^{(0)}\|_{\ell^{1}} + \sum_{k \ge 1}\rho^{k} \sum_{\ell=1}^{k} 
\frac{\ell!}{k!} \overline{\delta}^{k-\ell}
\|y^{(\ell)}\|_{\ell^{1}} \\
&\le e^{\rho\overline{\delta}}\|h^{(0)}\|_{\ell^{1}} 
+ \sum_{\ell \ge 1}\sum_{k \ge \ell} \frac{(\rho \overline{\delta})^{k-\ell}}
{(k-\ell)!} \rho^{\ell} \|y^{(\ell)}\|_{\ell^{1}} \\
&\le e^{\rho\overline{\delta}} \left ( \|h^{(0)}\|_{\ell^{1}} 
+ \|y\|_{\ell_{\rho}^{1}}\right ) \\
&\le \left ( \frac{e^{(1+\rho)\overline{\delta}}}
{1-e^{-\overline{\delta}}} + e^{\rho\overline{\delta}} \right ) 
\|y\|_{\ell_{\rho}^{1}}.
\end{aligned}
\end{equation*}

We now verify that the constructed $h$ solves $L_{0}h=y$. For $k\ge 1$, 
direct substitution of \eqref{eq:L0-explicit-inverse} gives \[
h^{(k)} + \frac{\overline{\delta}}{k}h^{(k-1)}=y^{(k)}.
\] To check the zeroth equation, summing the formula for $h^{(k)}$ 
and using the definition of $h^{(0)}$ gives
\begin{equation*}
\begin{aligned}
R^{-1}\sum_{k\ge 0}h^{(k)} - h^{(0)} &= 
-(I-e^{-\overline{\delta}}R^{-1})h^{(0)} + R^{-1} 
\sum_{k\ge 1}\sum_{\ell=1}^{k} \frac{\ell!}{k!}
(-\overline{\delta})^{k-\ell}y^{(\ell)}=y^{(0)}.
\end{aligned}
\end{equation*}
Thus $L_{0}h=y$. By the derivation above, every solution $h$ to 
$L_{0}h=y$ must have the form \eqref{eq:L0-explicit-inverse}.
Therefore $L_{0}$ is a bounded isomorphism, and the estimate for
$\|h\|_{\ell_{\rho}^{1}}$ proves \eqref{eq:L0-inverse-bound}.
\end{proof}

Whenever $\overline{\delta}>0$, Lemma \ref{lem:L0-inverse} 
allows us to define $\cG: \RR \times \ell_{\rho}^{1} \to    
\RR \times \ell_{\rho}^{1}$ by
\begin{equation} \label{eq:preconditioned-map}
\cG(\delta, a) := \left ( \phi(a^{(0)}), G(\delta, a) \right ) :=
\left ( \phi(a^{(0)}), L_{0}^{-1}F(\delta, a) \right ).
\end{equation}
The maps $\cF$ and $\cG$ have the same zeros.

\subsection{The Newton-like map and validation criterion}
\label{subsec:orbit-validation}

It remains to find a zero of $\cG(\delta, a)$ with
$\delta > 0$ and $a\not\equiv 0$. We will do so by constructing 
a Newton-like map and proving that it is a contraction near the 
numerical approximation $(\overline{\delta}, \overline{a})$. In this 
subsection, we set up an a posteriori framework for this argument.

At this stage, the parameter $\rho$ and the numerical data
$\overline{\delta}, \overline{a}, \overline{q}$ remain unspecified.
Their values will be specified in the final subsection, where the
conditions imposed on them will also be verified.

Assume $\overline{\delta}>0$, so that $L_{0}$ is invertible
by Lemma \ref{lem:L0-inverse}. Define \[
X:= \RR \times \ell_{\rho}^{1}, \quad 
\|(\delta, a)\|_{X} := |\delta| + \|a\|_{\ell_{\rho}^{1}},
\]
and write $\overline{x}:=(\overline{\delta}, \overline{a})$.
Throughout this subsection, all operator and multilinear norms are
induced by $\|\cdot\|_{X}$ and its restrictions.

To find a zero of $\cG$ near $\overline{x}$, let
$A:X \to X$ be a bounded linear operator to be chosen, and define
the Newton-like map \[
T(x):= x-A\cG(x).
\]
We seek $A$ such that $T$ maps a closed ball centered at $\overline{x}$ 
contractively into itself. It suffices to find $r, \epsilon>0$ such that
\begin{equation}\label{eq:orbit-contraction-conditions}
\|A\cG(\overline{x})\|\le \epsilon r, \quad 
\sup_{\|x-\overline{x}\|\le r}\|DT(x)\| \le 1-\epsilon.
\end{equation}
Indeed, the second bound shows that $T$ is a contraction, and together
these bounds give \[
\|T(x)-\overline{x}\| \le \|T(x)-T(\overline{x})\| + 
\|T(\overline{x})-\overline{x}\| \le (1-\epsilon)\|x-\overline{x}\| +
\|A\cG(\overline{x})\| \le r
\]
whenever $\|x-\overline{x}\|\le r$.

Since $\cG$ is a polynomial of degree at most three in $x$, Taylor's 
formula gives
\begin{equation}\label{eq:orbit-derivative-bound}
\sup_{\|x-\overline{x}\|\le r}\|DT(x)\| \le 
\|I-AD\cG(\overline{x})\| + r\|AD^{2}\cG(\overline{x})\|      
+\frac{r^{2}}{2}\|AD^{3}\cG(\overline{x})\|.
\end{equation}
We choose $A$ to make the first term small. The linear 
preconditioning in the previous subsection is designed to make 
$D\cG(\overline{x})$ close to the identity on the 
infinite-dimensional tail. We therefore approximate its inverse on a 
finite block and use the identity on the complement. To separate 
these two parts, let $P:X \to X$ be a finite-rank coordinate 
projection onto the finite block. Let 
$\overline{J}:PX \to PX$ be a numerical approximation to
$PD\cG(\overline{x})|_{PX}$, and let $A_{\mathrm{fin}}:PX \to PX$ be a 
numerically computed approximate inverse of $\overline{J}$. Define
\begin{equation}\label{eq:orbit-approximate-inverse}
A:=A_{\mathrm{fin}}P + (I-P).
\end{equation}

Extending $\overline{J}$ by the identity on the complement gives an 
approximation to $D\cG(\overline{x})$. Define the approximation
error by \[
E_{\mathrm{app}} := D\cG(\overline{x}) - \overline{J}P - (I-P).
\]
Then \[
I-AD\cG(\overline{x}) = (I_{PX}-A_{\mathrm{fin}}\overline{J})P 
-AE_{\mathrm{app}},
\]
and hence
\begin{equation}\label{eq:orbit-linear-error-bound}
\|I-AD\cG(\overline{x})\| \le \|I_{PX}-A_{\mathrm{fin}}\overline{J}\|
+ \|AE_{\mathrm{app}}\|.
\end{equation}

The following criterion combines the estimates 
\eqref{eq:orbit-contraction-conditions},
\eqref{eq:orbit-derivative-bound}, and 
\eqref{eq:orbit-linear-error-bound}.

\begin{proposition}\label{prop:orbit-validation}
Let $Z_{0}, Z_{1, \mathrm{inv}}, Z_{1, \mathrm{app}}, Z_{2}, Z_{3}>0$
be constants satisfying
\begin{equation}\label{eq:orbit-validation-estimates}
\begin{aligned}
\|A\cG(\overline{x})\| &\le Z_{0}, \\
\|I_{PX}-A_{\mathrm{fin}}\overline{J}\| &\le Z_{1, \mathrm{inv}}, \\
\|AE_{\mathrm{app}}\| &\le Z_{1, \mathrm{app}}, \\
\|AD^{2}\cG(\overline{x})\| &\le Z_{2}, \\
\frac{1}{2}\|AD^{3}\cG(\overline{x})\| &\le Z_{3}.
\end{aligned}
\end{equation}
Assume that there exists $r>0$ such that
\begin{equation}\label{eq:orbit-radii-inequality}
Z_{0}+(Z_{1, \mathrm{inv}}+Z_{1, \mathrm{app}}-1)r+Z_{2}r^{2} 
+Z_{3}r^{3} < 0.
\end{equation}
Then $\cG$ has a unique zero $(\delta^{(\ast)}, a^{(\ast)})$ satisfying
\[
|\delta^{(\ast)}-\overline{\delta}|
+\|a^{(\ast)}-\overline{a}\|_{\ell_{\rho}^{1}} \le r.
\]
If, in addition,
\begin{equation}\label{eq:orbit-radius-restrictions}
r< \min \left \{ \overline{\delta}, 
\|\overline{a}^{(0)}\|_{\ell^{1}} \right \},
\end{equation}
then $\delta^{(\ast)} > 0$ and $a^{(\ast)}\not\equiv 0$.
\end{proposition}

\begin{proof}
Set \[
\epsilon:= 1-Z_{1, \mathrm{inv}}-Z_{1, \mathrm{app}}-Z_{2}r-Z_{3}r^{2}.
\]
Inequality \eqref{eq:orbit-radii-inequality} gives $\epsilon>0$ 
and $Z_{0}<\epsilon r$. The estimates in 
\eqref{eq:orbit-validation-estimates}, together with 
\eqref{eq:orbit-derivative-bound} and
\eqref{eq:orbit-linear-error-bound}, therefore imply
\eqref{eq:orbit-contraction-conditions}.
Thus $T$ is a contraction of the ball
$\{x: \|x-\overline{x}\|\le r\}$ into itself, and the Banach 
fixed-point theorem gives a unique fixed point of $T$ in the ball.
To prove that this is also the unique zero of $\cG$, we are left to 
check that $A$ is injective.

To this end, observe that $Z_{1, \mathrm{inv}} < 1$ implies that
$A_{\mathrm{fin}}\overline{J}:PX \to PX$ is invertible. Since $PX$ is 
finite-dimensional, $A_{\mathrm{fin}}$ is invertible, and hence 
$A=A_{\mathrm{fin}}P+(I-P)$ is injective. This completes the proof.
\end{proof}

It remains to choose
$\rho, \overline{\delta}, \overline{a}, \overline{q}, P, 
\overline{J}$, and $A_{\mathrm{fin}}$ for which all the conditions 
required for validation are satisfied. These include $\rho>1$ and 
$\overline{\delta}>0$, the estimates 
\eqref{eq:orbit-validation-estimates}, and the existence of $r>0$
satisfying both \eqref{eq:orbit-radii-inequality} and
\eqref{eq:orbit-radius-restrictions}.

\subsection{Reduction to finite computations}
\label{subsec:orbit-analytic-estimates}

We now reduce the verification of the estimates
\eqref{eq:orbit-validation-estimates} to finite computations.
The estimate involving $Z_{1, \mathrm{inv}}$ is already a finite
matrix calculation. For the remaining estimates, we derive upper 
bounds that can be evaluated by finite computations. The bound 
corresponding to $Z_{0}$ is given in \eqref{eq:orbit-residual-bound}, 
and those corresponding to $Z_{2}$ and $Z_{3}$ are given in Lemma 
\ref{lem:orbit-higher-derivatives}. The bound corresponding to 
$Z_{1, \mathrm{app}}$ is obtained in Proposition 
\ref{prop:orbit-derivative-error}.

We continue to assume $\rho>1$ and $\overline{\delta}>0$.
We now specify the finite-rank coordinate projection $P:X\to X$ used
to define $A$ in \eqref{eq:orbit-approximate-inverse}. Choose integers 
$L\le U$ and $K\ge 0$, and define
\begin{equation}\label{eq:orbit-finite-block}
P(\delta, a) := (\delta, \widetilde{a}),
\quad \widetilde{a}_{n}^{(k)} := 
\begin{cases}
a_{n}^{(k)} & \text{if }L \le n \le U, \, 0 \le k \le K, \\
0 & \text{otherwise.}
\end{cases}
\end{equation}
We also assume that $\overline{a}$ is finitely supported, with
$\overline{a}^{(k)}=0$ for $k>K$. We do not require the spatial 
support of $\overline{a}$ to lie in $[L, U]$.

As in the previous subsection, all operator and multilinear norms
below are induced by $\|\cdot\|_{X}$ and its restrictions.

We begin with the estimate involving $Z_{0}$. Since
$\phi(\overline{a}^{(0)})=0$, we have
\begin{equation}\label{eq:orbit-residual-bound}
\|A\cG(\overline{x})\| \le \|A\| \|L_{0}^{-1}\|
\|F(\overline{x})\|_{\ell_{\rho}^{1}}.
\end{equation}
Every factor on the right can be bounded by finite computations: \[
\|A\| = \max \{ 1, \|A_{\mathrm{fin}}\| \},
\]
\eqref{eq:L0-inverse-bound} bounds $\|L_{0}^{-1}\|$, and
$F(\overline{x})$ is finitely supported because $\overline{a}$ is
finitely supported and $Q$ has finite interaction range.

We next obtain finitely computable bounds for the terms corresponding 
to $Z_{2}$ and $Z_{3}$.

\begin{lemma}[Higher-derivative bounds]
\label{lem:orbit-higher-derivatives}
We have
\begin{equation}\label{eq:orbit-higher-derivative-bounds}
\begin{aligned}
\|AD^{2}\cG(\overline{x})\| &\le \|A\| \|L_{0}^{-1}\|
\max \left \{ \rho(1+2C_{Q}\|\overline{a}\|_{\ell_{\rho}^{1}}),
2\overline{\delta}\rho C_{Q} \right \}, \\
\frac{1}{2}\|AD^{3}\cG(\overline{x})\| &\le \|A\| \|L_{0}^{-1}\|
\rho C_{Q}.
\end{aligned}
\end{equation}
In particular, both right-hand sides can be bounded using finite 
computations.
\end{lemma}

\begin{proof}
It suffices to estimate the second and third derivatives of
$F^{(k)}$ for $k\ge 1$, since $\phi$ and $F^{(0)}$ are affine.
Differentiation and \eqref{eq:orbit-bilinear-bound} give
\begin{equation*}
\begin{aligned}
&\|D^{2}F(\overline{x})[(\delta_{1}, h_{1}), (\delta_{2}, h_{2})]
\|_{\ell_{\rho}^{1}} \\
&= \sum_{k \ge 1}\rho^{k} \bigg \| \frac{\delta_{1}}{k}h_{2}^{(k-1)}
+\frac{\delta_{2}}{k}h_{1}^{(k-1)}
-\frac{2\delta_{1}}{k}\sum_{j=0}^{k-1}Q(\overline{a}^{(j)}, 
h_{2}^{(k-1-j)})
 \\
&\qquad \qquad \quad 
-\frac{2\delta_{2}}{k}\sum_{j=0}^{k-1}Q(\overline{a}^{(j)}, 
h_{1}^{(k-1-j)})
  -\frac{2\overline{\delta}}{k}\sum_{j=0}^{k-1}Q(h_{1}^{(j)},
h_{2}^{(k-1-j)}) \bigg \|_{\ell^{1}} \\
&\le \rho (|\delta_{1}|\|h_{2}\|_{\ell_{\rho}^{1}} 
+|\delta_{2}|\|h_{1}\|_{\ell_{\rho}^{1}}) \\
&+ 2C_{Q} \sum_{k \ge 1}\rho^{k} \sum_{j=0}^{k-1} \left [ 
|\delta_{1}|\|\overline{a}^{(j)}\|_{\ell^{1}}
\|h_{2}^{(k-1-j)}\|_{\ell^{1}} 
+|\delta_{2}|\|\overline{a}^{(j)}\|_{\ell^{1}}
\|h_{1}^{(k-1-j)}\|_{\ell^{1}} 
+\overline{\delta}\|h_{1}^{(j)}\|_{\ell^{1}} 
\|h_{2}^{(k-1-j)}\|_{\ell^{1}} \right ] \\
&\le \rho(1+2C_{Q}\|\overline{a}\|_{\ell_{\rho}^{1}})
(|\delta_{1}|\|h_{2}\|_{\ell_{\rho}^{1}} 
+|\delta_{2}|\|h_{1}\|_{\ell_{\rho}^{1}} )
+ 2\overline{\delta}\rho C_{Q} \|h_{1}\|_{\ell_{\rho}^{1}}
\|h_{2}\|_{\ell_{\rho}^{1}} \\
&\le \max \left \{ \rho (1+2C_{Q}\|\overline{a}\|_{\ell_{\rho}^{1}}), 
2\overline{\delta}\rho C_{Q} \right \}
\prod_{i=1}^{2}\|(\delta_{i}, h_{i})\|_{X}.
\end{aligned}
\end{equation*}
Similarly,
\begin{equation*}
\begin{aligned}
& \|D^{3}F(\overline{x})[(\delta_{1}, h_{1}), 
(\delta_{2}, h_{2}), (\delta_{3}, h_{3})]\|_{\ell_{\rho}^{1}} \\ 
&= \sum_{k \ge 1}\rho^{k}\frac{1}{k}\left \| \sum_{j=0}^{k-1} 
\left [ 2\delta_{1}Q(h_{2}^{(j)}, h_{3}^{(k-1-j)}) 
+ 2\delta_{2}Q(h_{3}^{(j)}, h_{1}^{(k-1-j)})
+ 2\delta_{3}Q(h_{1}^{(j)}, h_{2}^{(k-1-j)}) \right ] 
\right \|_{\ell^{1}} \\
&\le 2\rho C_{Q} \left ( |\delta_{1}|\|h_{2}\|_{\ell_{\rho}^{1}} 
\|h_{3}\|_{\ell_{\rho}^{1}} 
+ |\delta_{2}|\|h_{3}\|_{\ell_{\rho}^{1}} 
\|h_{1}\|_{\ell_{\rho}^{1}} 
+ |\delta_{3}|\|h_{1}\|_{\ell_{\rho}^{1}} 
\|h_{2}\|_{\ell_{\rho}^{1}} \right ) \\
&\le 2\rho C_{Q} \prod_{i=1}^{3}\|(\delta_{i}, h_{i})\|_{X}.
\end{aligned}
\end{equation*}
Applying $A(0, L_{0}^{-1}\,\cdot\,)$ gives the stated bounds.
\end{proof}

It remains to estimate $\|AE_{\mathrm{app}}\|$. For $m \in \ZZ$ and 
$\ell \ge 0$, let $e_{m, \ell} \in \ell_{\rho}^{1}$ be the 
normalized coordinate vector corresponding to $a_{m}^{(\ell)}$: \[
(e_{m, \ell})_{n}^{(k)} := \begin{cases}
\rho^{-\ell} & \text{if }(n, k)=(m, \ell), \\
0 & \text{otherwise.}
\end{cases}
\]
Since $X$ has a weighted $\ell^{1}$ norm, $\|AE_{\mathrm{app}}\|$
is the supremum of the column norms:
\begin{equation}\label{eq:orbit-column-norm}
\|AE_{\mathrm{app}}\|=\max\left \{ \|AE_{\mathrm{app}}(1, 0)\|,
\sup_{m\in \ZZ, \ell\ge 0} \|AE_{\mathrm{app}}(0, e_{m, \ell})\|
\right \}.
\end{equation}

In Subsubsection \ref{subsubsec:orbit-fixed-columns}, we derive a 
finitely computable bound for $\|AE_{\mathrm{app}}\mathbf{e}\|$ for 
any fixed input $\mathbf{e}=(1, 0)$ or $(0, e_{m, \ell})$. This bound 
will be applied individually to finitely many inputs $\mathbf{e}$. 
For all other inputs, we will either show that the output 
$AE_{\mathrm{app}}\mathbf{e}$ vanishes or apply the uniform bound in Lemma 
\ref{lem:orbit-Taylor-comparison}. Together, these estimates bound 
all the output norms in \eqref{eq:orbit-column-norm} and give the 
desired estimate for $\|AE_{\mathrm{app}}\|$ in Proposition 
\ref{prop:orbit-derivative-error}.

\subsubsection{Bounds for a fixed input}
\label{subsubsec:orbit-fixed-columns}

Fix an input $\mathbf{e}=(1, 0)$ or $(0, e_{m, \ell})$. 
We split the output 
$AE_{\mathrm{app}}\mathbf{e}$ into its $PX$ and $(I-P)X$ parts:
\begin{equation}\label{eq:orbit-column-splitting}
\begin{aligned}
AE_{\mathrm{app}}\mathbf{e} &= (A_{\mathrm{fin}}P+(I-P))
(D\cG(\overline{x})-\overline{J}P-(I-P))\mathbf{e} \\
&= \underbrace{
A_{\mathrm{fin}}(PD\cG(\overline{x})\mathbf{e} 
-\overline{J}P\mathbf{e})}_{\in PX} + 
\underbrace{(I-P)(D\cG(\overline{x})\mathbf{e} 
-\mathbf{e})}_{\in (I-P)X}.
\end{aligned}
\end{equation}
The only terms that still need to be reduced to finite computations 
are those involving $D\cG(\overline{x})\mathbf{e}$. For the 
$PX$ component, we need bounds for the finitely many coordinates of
$PD\cG(\overline{x})\mathbf{e}$, and for the $(I-P)X$
component, we need a bound for 
$\|(I-P)(D\cG(\overline{x})\mathbf{e}-\mathbf{e})\|$.

Write $\mathbf{e} = (\eta, b)$. Then \eqref{eq:preconditioned-map} 
gives 
\begin{equation}\label{eq:orbit-derivative-output}
D\cG(\overline{x})\mathbf{e} = \left ( \sum_{n \in \ZZ} 
\overline{q}_{n}b_{n}^{(0)}, L_{0}^{-1}DF(\overline{x})\mathbf{e}
\right ).
\end{equation}
As $\overline{q}$ is finitely supported, the first entry is a finite 
sum. So we only need to study $L_{0}^{-1}DF(\overline{x})\mathbf{e}$.

For the period input $\mathbf{e}=(1, 0)$, we have \[
L_{0}^{-1}DF(\overline{x})\mathbf{e}=L_{0}^{-1}y_{\delta},
\]
where $y_{\delta}:=D_{\delta}F(\overline{x})$ is finitely supported,
with
\begin{equation}\label{eq:orbit-period-source}
\begin{aligned}
y_{\delta}^{(0)} &= 0, \\
y_{\delta}^{(k)} &= \frac{1}{k}\left ( \overline{a}^{(k-1)} - 
\sum_{j=0}^{k-1}Q(\overline{a}^{(j)}, \overline{a}^{(k-1-j)}) \right )
\quad (k \ge 1).
\end{aligned}
\end{equation}
In particular, $y_{\delta}^{(k)}=0$ for $k > 2K+1$, since 
$\overline{a}^{(k)}=0$ for $k>K$.

For a coefficient input $\mathbf{e}=(0, e_{m, \ell})$, one can 
compute that \begin{equation}\label{eq:orbit-output-coordinate-input}
L_{0}^{-1}DF(\overline{x})\mathbf{e} = L_{0}^{-1}D_{a}F(\overline{x})
e_{m, \ell} = e_{m, \ell} + L_{0}^{-1}y_{m, \ell},
\end{equation}
where
\begin{equation}\label{eq:orbit-coefficient-source}
(y_{m, \ell})^{(k)} := \begin{cases}
-\frac{2\overline{\delta}}{k}\rho^{-\ell} 
Q(\overline{a}^{(k-1-\ell)}, \mathbf{1}_{\{m\}})
& \text{if }\ell + 1 \le k \le \ell+K+1, \\
0 & \text{otherwise},
\end{cases}
\end{equation}
where $\mathbf{1}_{\{m\}} \in \ell^{1}(\ZZ)$ is the sequence
equal to one at $m$ and zero elsewhere.
Note that the spatial support of $y_{m, \ell}$ lies in 
$[m-2, m+2]$, so it is also finitely supported.

For either type of input $\mathbf{e}$, to bound 
$AE_{\mathrm{app}}\mathbf{e}$ using finite computations, it is 
therefore enough to bound the coordinates of 
$P(0, L_{0}^{-1}y)$ and the norm $\|(I-P)(0, L_{0}^{-1}y)\|$ for 
finitely supported $y$; all other terms in 
\eqref{eq:orbit-column-splitting} require only finite computations.

However, the formula \eqref{eq:L0-explicit-inverse} for
$L_{0}^{-1}$ involves infinite sums, so even the evaluation of 
$P(0, L_{0}^{-1}y)$ is not yet reduced to finite computations. The next 
lemma gives finite expressions for the coordinates of $L_{0}^{-1}y$.

\begin{lemma}[Finite-coordinate formulas]
\label{lem:orbit-finite-coordinates}
Let $y\in\ell_{\rho}^{1}$ be finitely supported.
Choose integers $M\ge 0$ and $n_{+} \in \ZZ$ such that \[
y_{n}^{(k)}=0 \quad \text{if }k>M \text{ or }n>n_{+}.
\]
Then $(L_{0}^{-1}y)_{n}^{(k)}=0$ for any $n>n_{+}$ and $k \ge 0$. 
For $n \le n_{+}$, we have
\begin{align}
(L_{0}^{-1}y)_{n}^{(0)} &= \sum_{j=0}^{n_{+}-n} e^{-j\overline{\delta}}
\left [ -y_{n+j}^{(0)} + \sum_{\ell=1}^{M} 
\ell! (-\overline{\delta})^{-\ell} \left ( e^{-\overline{\delta}} 
-\sum_{k=0}^{\ell-1} \frac{(-\overline{\delta})^{k}}{k!} \right )   
y_{n+j+1}^{(\ell)} \right ], \label{eq:orbit-finite-zeroth-coefficient}
\\
(L_{0}^{-1}y)_{n}^{(k)} &= \frac{(-\overline{\delta})^{k}}{k!}
(L_{0}^{-1}y)_{n}^{(0)} + \sum_{p=1}^{\min\{k, M\}} \frac{p!}{k!}
(-\overline{\delta})^{k-p}y_{n}^{(p)} \quad (k \ge 1).
\label{eq:orbit-finite-higher-coefficients}
\end{align}
In particular, any fixed finite collection of these coordinates, for 
instance $P(0, L_{0}^{-1}y)$, can be bounded using finite computations.
\end{lemma}

\begin{remark}\label{rmk:orbit-coefficient-evaluation}
For large $\ell$, the expression \[
e^{-\overline{\delta}}-\sum_{k=0}^{\ell-1}\frac{(-\overline{\delta})^{k}}
{k!}
\]
in \eqref{eq:orbit-finite-zeroth-coefficient} involves subtracting 
two nearly equal numbers, which may lead to a substantial 
loss of relative accuracy. Although direct evaluation is possible 
at sufficiently high precision, we use the equivalent representation \[
\ell!(-\overline{\delta})^{-\ell} \left ( e^{-\overline{\delta}}
-\sum_{k=0}^{\ell-1}\frac{(-\overline{\delta})^{k}}{k!}\right )
= \sum_{j=0}^{\infty} \frac{\ell!(-\overline{\delta})^{j}}{(\ell+j)!}
\quad \ell\ge 1.
\]
When $0<\overline{\delta}<1$, this is an alternating series whose 
terms decrease in absolute value. Hence, for any integer $N \ge 1$, \[
\sum_{j=0}^{\infty}\frac{\ell!(-\overline{\delta})^{j}}{(\ell+j)!}
\in \left [ 
\sum_{j=0}^{2N-1}\frac{\ell!(-\overline{\delta})^{j}}{(\ell+j)!},
\sum_{j=0}^{2N}\frac{\ell!(-\overline{\delta})^{j}}{(\ell+j)!}
\right ].
\]
This will be used in the computer-assisted verification.
\end{remark}

\begin{proof}[Proof of Lemma \ref{lem:orbit-finite-coordinates}]
Expanding $(I-e^{-\overline{\delta}}R^{-1})^{-1}
=\sum_{j\ge 0}e^{-j\overline{\delta}}R^{-j}$ in
\eqref{eq:L0-explicit-inverse} gives
\begin{equation*}
\begin{aligned}
(L_{0}^{-1}y)_{n}^{(0)} &= \sum_{j \ge 0}e^{-j\overline{\delta}} 
\left ( \sum_{k\ge 1}\sum_{\ell=1}^{k}\frac{\ell!}{k!}
(-\overline{\delta})^{k-\ell}y_{n+j+1}^{(\ell)} -y_{n+j}^{(0)} \right ) \\
&= \sum_{0 \le j \le n_{+}-n} e^{-j\overline{\delta}} \left [ -y_{n+j}^{(0)} 
+ \sum_{\ell \ge 1}\sum_{k \ge \ell}\frac{\ell!}{k!}
(-\overline{\delta})^{k-\ell}y_{n+j+1}^{(\ell)}\right ] \\
&= \sum_{0 \le j \le n_{+}-n} 
e^{-j\overline{\delta}} \left [ -y_{n+j}^{(0)} 
+ \sum_{\ell=1}^{M}\ell! (-\overline{\delta})^{-\ell} 
\sum_{k \ge \ell} \frac{(-\overline{\delta})^{k}}{k!} 
y_{n+j+1}^{(\ell)}\right ].
\end{aligned}
\end{equation*}
This proves \eqref{eq:orbit-finite-zeroth-coefficient}. The formula 
\eqref{eq:orbit-finite-higher-coefficients} follows directly from 
\eqref{eq:L0-explicit-inverse}. These formulas also show that 
$(L_{0}^{-1}y)_{n}^{(k)}=0$ for $n>n_{+}$.
\end{proof}

We next need to bound $\|(I-P)(0, L_{0}^{-1}y)\|$. The following 
lemma bounds this
norm using only finitely many output coordinates, each of which can be 
computed using Lemma \ref{lem:orbit-finite-coordinates}.

\begin{lemma}[Bound outside the finite block]
\label{lem:orbit-output-complement}
Let $y \in \ell_{\rho}^{1}$ be finitely supported. Choose 
integers $M \ge K$, $n_{-} \le L$, and $n_{+} \ge U$, with 
$M+1>\rho\overline{\delta}$, such that \[
y_{n}^{(k)}=0 \quad \text{if }k>M, \, n <n_{-}, \, \text{or }n>n_{+}.
\]
Then
\begin{equation}\label{eq:orbit-output-complement-bound}
\begin{aligned}
\|(I-P)(0, L_{0}^{-1}y)\| &\le 
\sum_{\substack{n_{-}\le n \le n_{+}, \, 0 \le k \le M \\
n \not\in [L, U] \text{ or }k>K}} \rho^{k}
|(L_{0}^{-1}y)_{n}^{(k)}| \\
&+ \frac{\rho\overline{\delta}}{M+1-\rho\overline{\delta}}
\rho^{M}\sum_{n_{-}\le n\le n_{+}}|(L_{0}^{-1}y)_{n}^{(M)}| \\
&+ \frac{ e^{\rho\overline{\delta}}}
{1-e^{-\overline{\delta}}}|(L_{0}^{-1}y)_{n_{-}-1}^{(0)}|.
\end{aligned}
\end{equation}
The right-hand side can be bounded by finite computations using
Lemma \ref{lem:orbit-finite-coordinates}.
\end{lemma}

\begin{proof}
It is enough to prove that \begin{equation}
\label{eq:orbit-taylor-tail-bound}
\sum_{n_{-}\le n \le n_{+}, \, k>M} \rho^{k}
|(L_{0}^{-1}y)_{n}^{(k)}| \le
\frac{\rho\overline{\delta}}{M+1-\rho\overline{\delta}}
\rho^{M}\sum_{n_{-}\le n\le n_{+}}|(L_{0}^{-1}y)_{n}^{(M)}|,
\end{equation} and
\begin{equation}
\label{eq:orbit-spatial-tail-bound}
\sum_{n \not\in [n_{-}, n_{+}], \, k \ge 0} \rho^{k}
|(L_{0}^{-1}y)_{n}^{(k)}| \le       
\frac{ e^{\rho\overline{\delta}}}
{1-e^{-\overline{\delta}}}|(L_{0}^{-1}y)_{n_{-}-1}^{(0)}|.
\end{equation}
To prove \eqref{eq:orbit-taylor-tail-bound}, observe that
the recurrence \eqref{eq:L0-formula-high} gives \[
(L_{0}^{-1}y)_{n}^{(k)} = \frac{M!}{k!}(-\overline{\delta})^{k-M} 
(L_{0}^{-1}y)_{n}^{(M)} \quad (k > M).
\]
Therefore,
\begin{equation*}
\begin{aligned}
\sum_{n_{-}\le n \le n_{+}, \, k>M} \rho^{k}
|(L_{0}^{-1}y)_{n}^{(k)}| &= 
\sum_{n_{-}\le n\le n_{+}} \rho^{M}|(L_{0}^{-1}y)_{n}^{(M)}|
\sum_{k \ge M+1}\rho^{k-M} \frac{M!}{k!}\overline{\delta}^{k-M} \\
& \le \sum_{n_{-}\le n\le n_{+}}\rho^{M}|(L_{0}^{-1}y)_{n}^{(M)}| 
\sum_{j \ge 1} (\rho\overline{\delta})^{j}\frac{1}{(M+1)^{j}} \\
&= \frac{\rho\overline{\delta}}{M+1-\rho\overline{\delta}}
\rho^{M}\sum_{n_{-}\le n\le n_{+}}|(L_{0}^{-1}y)_{n}^{(M)}|.
\end{aligned}
\end{equation*}
We are left to prove \eqref{eq:orbit-spatial-tail-bound}. Note 
that for $n>n_{+}$, $(L_{0}^{-1}y)_{n}^{(k)}=0$ by Lemma 
\ref{lem:orbit-finite-coordinates}. Also for $n<n_{-}-1$, 
\eqref{eq:orbit-finite-zeroth-coefficient} gives
\begin{equation*}
\begin{aligned}
(L_{0}^{-1}y)_{n+1}^{(0)} &= \sum_{j=0}^{n_{+}-n-1} e^{-j\overline{\delta}}
\left [ -y_{n+j+1}^{(0)} + \sum_{\ell=1}^{M} 
\ell! (-\overline{\delta})^{-\ell} \left ( e^{-\overline{\delta}} 
-\sum_{k=0}^{\ell-1} \frac{(-\overline{\delta})^{k}}{k!} \right )   
y_{n+j+2}^{(\ell)} \right ] \\
&=e^{\overline{\delta}} \sum_{j=1}^{n_{+}-n} e^{-j\overline{\delta}}
\left [ -y_{n+j}^{(0)} + \sum_{\ell=1}^{M} 
\ell! (-\overline{\delta})^{-\ell} \left ( e^{-\overline{\delta}} 
-\sum_{k=0}^{\ell-1} \frac{(-\overline{\delta})^{k}}{k!} \right )   
y_{n+j+1}^{(\ell)} \right ] \\
&= e^{\overline{\delta}}(L_{0}^{-1}y)_{n}^{(0)},
\end{aligned}
\end{equation*}
so that \[
(L_{0}^{-1}y)_{n}^{(0)} = e^{-\overline{\delta}(n_{-}-n-1)} 
(L_{0}^{-1}y)_{n_{-}-1}^{(0)} \quad (n < n_{-}).
\]
Then \eqref{eq:L0-formula-high} gives \[
(L_{0}^{-1}y)_{n}^{(k)} = \frac{(-\overline{\delta})^{k}}{k!}
e^{-\overline{\delta}(n_{-}-n-1)}
(L_{0}^{-1}y)_{n_{-}-1}^{(0)} \quad (n < n_{-}, \, k \ge 1).
\]
Thus,
\begin{equation*}
\begin{aligned}
\sum_{n \not\in [n_{-}, n_{+}], \, k \ge 0} \rho^{k}
|(L_{0}^{-1}y)_{n}^{(k)}|
&= \sum_{n < n_{-}}\sum_{k=0}^{\infty} \rho^{k} 
\frac{\overline{\delta}^{k}}{k!} e^{-\overline{\delta}(n_{-}-n-1)}
|(L_{0}^{-1}y)_{n_{-}-1}^{(0)}| \\
&= |(L_{0}^{-1}y)_{n_{-}-1}^{(0)}|\sum_{n < n_{-}} 
e^{\rho\overline{\delta}}  e^{-\overline{\delta}(n_{-}-n-1)} \\
&= \frac{e^{\rho\overline{\delta}}}
{1-e^{-\overline{\delta}}}|(L_{0}^{-1}y)_{n_{-}-1}^{(0)}|.
\end{aligned}
\end{equation*}
This completes the proof.
\end{proof}

\subsubsection{Uniform bounds for infinitely many inputs}
\label{subsubsec:orbit-input-columns}

The preceding two lemmas allow us to bound
$\|AE_{\mathrm{app}}\mathbf{e}\|$ by finite computations for each 
fixed input $\mathbf{e}$.
Finitely many inputs will be treated in this way. We now obtain
a uniform bound that deals with the remaining inputs.

For inputs $\mathbf{e}=(0, e_{m, \ell})$, 
\eqref{eq:orbit-derivative-output}, \eqref{eq:orbit-output-coordinate-input}, 
and the definition of $E_{\mathrm{app}}$ give
\begin{equation}\label{eq:orbit-coefficient-column}
\begin{aligned}
E_{\mathrm{app}}(0, e_{m, \ell})
&= D\cG(\overline{x})(0, e_{m, \ell}) 
-\overline{J}P(0, e_{m, \ell}) - (I-P)(0, e_{m, \ell}) \\
&= \left ( \overline{q}_{m} \mathbf{1}_{\{\ell=0\}},
L_{0}^{-1}y_{m, \ell} \right ) + (I_{PX}-\overline{J})P(0, e_{m, \ell}).
\end{aligned}
\end{equation}
Since $Q$ has finite interaction range,
\eqref{eq:orbit-coefficient-source} shows that $y_{m, \ell}=0$
for every $\ell$ when $m$ is sufficiently far from the spatial
support of $\overline{a}$.
The other terms in \eqref{eq:orbit-coefficient-column} also vanish
for sufficiently large $|m|$.
Thus, if we choose a finite set $\cS\subset\ZZ$ such that
\begin{equation}\label{eq:orbit-spatial-input-set}
\begin{aligned}
\cS\supset
[L, U] \cup \mathrm{supp}\overline{q} \cup \left \{ m \in \ZZ: 
Q(\overline{a}^{(j)}, \mathbf{1}_{\{m\}}) \neq 0 \text{ for some }
0\le j\le K \right \},
\end{aligned}
\end{equation}
then we have
\begin{equation}\label{eq:orbit-remote-columns}
E_{\mathrm{app}}(0, e_{m, \ell})=0 \quad (m \not\in\cS, \, 
\ell \ge 0).
\end{equation}

It therefore remains to consider inputs $(0, e_{m, \ell})$ with
$m \in \cS$. Fix an integer $\ell_{0}>K$. For the finitely many 
inputs with $0 \le \ell < \ell_{0}$, we bound the output norms as 
in Subsubsection \ref{subsubsec:orbit-fixed-columns}. To handle the 
infinitely many indices $\ell \ge \ell_{0}$, we need a bound uniform 
in $\ell$. The next lemma gives such a bound.

\begin{lemma}[Uniform bound for high Taylor degrees]
\label{lem:orbit-Taylor-comparison}
For each $m\in \cS$,
\begin{equation}\label{eq:orbit-high-Taylor-column-bound}
\begin{aligned}
\sup_{\ell\ge \ell_{0}} \|AE_{\mathrm{app}}(0, e_{m, \ell})\| \le     
\sum_{k=\ell_{0}+1}^{\ell_{0}+K+1}\sum_{n=m-2}^{m+2}\rho^{k}
|(y_{m, \ell_{0}})_{n}^{(k)}| \|A(0, L_{0}^{-1}e_{n, k})\|.
\end{aligned}
\end{equation}
Note that, since \[
\|A(0, L_{0}^{-1}e_{n, k})\| = 
\|A_{\mathrm{fin}}P(0, L_{0}^{-1}e_{n, k})\| 
+ \|(I-P)(0, L_{0}^{-1}e_{n, k})\|,
\]
the right-hand side of \eqref{eq:orbit-high-Taylor-column-bound} can 
be bounded by finite computations using Lemmas 
\ref{lem:orbit-finite-coordinates} and \ref{lem:orbit-output-complement}.
\end{lemma}

\begin{proof}
For $\ell=\ell_{0}+s$ with $s \ge 0$,
\eqref{eq:orbit-coefficient-column} gives \[
AE_{\mathrm{app}}(0, e_{m, \ell_{0}+s})
= A(0, L_{0}^{-1}y_{m, \ell_{0}+s}).
\]
Also, \eqref{eq:orbit-coefficient-source} gives
\begin{equation*}
\begin{aligned}
(y_{m, \ell_{0}+s})^{(k+s)} &= \begin{cases}
-\frac{2\overline{\delta}}{k+s}\rho^{-\ell_{0}-s} 
Q(\overline{a}^{(k-1-\ell_{0})}, \mathbf{1}_{\{m\}})
& \text{if }\ell_{0} + 1 \le k \le \ell_{0}+K+1, \\
0 & \text{otherwise}
\end{cases} \\
&= \begin{cases}
(y_{m, \ell_{0}})^{(k)} \frac{k}{k+s} \rho^{-s} & 
\text{if } \ell_{0}+1 \le k \le \ell_{0}+K+1, \\
0 & \text{otherwise},
\end{cases}
\end{aligned}
\end{equation*}
so that
\begin{equation*}
\begin{aligned}
y_{m, \ell_{0}+s} &= \sum_{k=\ell_{0}+1}^{\ell_{0}+K+1}
\sum_{n=m-2}^{m+2} (y_{m, \ell_{0}+s})_{n}^{(k+s)}
 \rho^{k+s} e_{n, k+s} \\
&= \sum_{k=\ell_{0}+1}^{\ell_{0}+K+1} \sum_{n=m-2}^{m+2}
\rho^{k} \frac{k}{k+s} (y_{m, \ell_{0}})_{n}^{(k)}e_{n, k+s}.
\end{aligned}
\end{equation*}
Therefore, to prove \eqref{eq:orbit-high-Taylor-column-bound},
it suffices to prove that 
\begin{equation} \label{eq:orbit-Taylor-comparison}
\frac{k}{k+s}\|A(0, L_{0}^{-1}e_{n, k+s})\| 
\le \|A(0, L_{0}^{-1}e_{n, k})\| \quad (k>\ell_{0}, \, n \in \ZZ, \, s \ge 0).
\end{equation}
Let $k>\ell_{0}>K$, $n \in \ZZ$, and $s \ge 0$. Observe that by 
\eqref{eq:L0-explicit-inverse}, we have 
\begin{equation*}
\begin{aligned}
(L_{0}^{-1}e_{n, k})^{(0)} &= (I-e^{-\overline{\delta}}R^{-1})^{-1} 
\left ( R^{-1}\sum_{p \ge k}\frac{k!}{p!}
(-\overline{\delta})^{p-k}\rho^{-k}\mathbf{1}_{\{n\}} \right ), \\
(L_{0}^{-1}e_{n-1, 0})^{(0)} &= (I-e^{-\overline{\delta}}R^{-1})^{-1}
\left ( - \mathbf{1}_{\{n-1\}}\right ),
\end{aligned}
\end{equation*}
so that \[
(L_{0}^{-1}e_{n, k})^{(0)} = - \rho^{-k} 
\sum_{p \ge k}\frac{k!}{p!}(-\overline{\delta})^{p-k}
(L_{0}^{-1}e_{n-1, 0})^{(0)}.
\] The zeroth coefficient $(L_{0}^{-1}e_{n, k})^{(0)}$ is therefore supported 
in $\{n': n' \le n-1\}$.
Again by \eqref{eq:L0-explicit-inverse}, we have 
\begin{equation*}
\begin{aligned}
(L_{0}^{-1}e_{n, k})^{(q)} &= \frac{(-\overline{\delta})^{q}}{q!}
(L_{0}^{-1}e_{n, k})^{(0)} + \mathbf{1}_{\{q \ge k\}} \frac{k!}{q!}
(-\overline{\delta})^{q-k}\rho^{-k}\mathbf{1}_{\{n\}}, \\
(L_{0}^{-1}e_{n-1, 0})^{(q)} &= \frac{(-\overline{\delta})^{q}}{q!}
(L_{0}^{-1}e_{n-1, 0})^{(0)}.
\end{aligned}
\end{equation*}
Therefore, \[
L_{0}^{-1}e_{n, k} = -\rho^{-k} 
\sum_{p \ge k}\frac{k!}{p!}(-\overline{\delta})^{p-k} 
L_{0}^{-1}e_{n-1, 0} + E_{n, k},
\] where \[
(E_{n, k})^{(q)} = \mathbf{1}_{\{q \ge k\}}\frac{k!}{q!}
(-\overline{\delta})^{q-k}\rho^{-k} \mathbf{1}_{\{n\}}.
\] Note that the support of $E_{n, k}$ is contained in 
$\{n\} \times \{q: q \ge k\}$, which lies outside the finite block 
in \eqref{eq:orbit-finite-block}. Thus, 
$A(0, E_{n, k})=(0, E_{n, k})$. On the other hand, the support of 
$L_{0}^{-1}e_{n-1, 0}$ is contained in 
$\{n': n' \le n-1\} \times \ZZ_{\ge 0}$,
which in particular is disjoint from the support of $E_{n, k}$. Since 
$A$ maps $PX$ into itself and acts as the identity on $(I-P)X$, 
the support of $A(0, L_{0}^{-1}e_{n-1, 0})$ is still disjoint from 
the support of $E_{n, k}$. Therefore, 
\begin{equation*}
\begin{aligned}
\|A(0, L_{0}^{-1}e_{n, k})\| &= \rho^{-k} 
\left | \sum_{p \ge k} \frac{k!}{p!}(-\overline{\delta})^{p-k} \right | 
\|A(0, L_{0}^{-1}e_{n-1, 0})\| + \|(0, E_{n, k})\| \\ 
&= \rho^{-k} 
\left | \sum_{p \ge k} \frac{k!}{p!}(-\overline{\delta})^{p-k} \right |
\|A(0, L_{0}^{-1}e_{n-1, 0})\|
+ \sum_{q \ge k}\rho^{q-k}\frac{k!}{q!}\overline{\delta}^{q-k}.
\end{aligned}
\end{equation*}
Replacing $k$ by $k+s$, we have 
\begin{equation*}
\begin{aligned}
&\|A(0, L_{0}^{-1}e_{n, k+s})\| \\
&= \rho^{-k-s} 
\left | \sum_{p \ge k+s} \frac{(k+s)!}{p!}(-\overline{\delta})^{p-k-s}
\right | \|A(0, L_{0}^{-1}e_{n-1, 0})\|
+ \sum_{q \ge k+s}\rho^{q-k-s}\frac{(k+s)!}{q!}
\overline{\delta}^{q-k-s}.
\end{aligned}
\end{equation*}
Therefore, to prove \eqref{eq:orbit-Taylor-comparison}, it is enough 
to prove the following two inequalities:
\begin{equation*}
\begin{aligned}
\frac{k}{k+s}\rho^{-s} 
\left | \sum_{p \ge k+s} \frac{(k+s)!}{p!}(-\overline{\delta})^{p-k-s}
\right | &\le \left | \sum_{p \ge k} \frac{k!}{p!}
(-\overline{\delta})^{p-k} \right |, \\
\frac{k}{k+s}\sum_{q \ge k+s} (\rho\overline{\delta})^{q-k-s}
\frac{(k+s)!}{q!} &\le \sum_{q \ge k}(\rho\overline{\delta})^{q-k}
\frac{k!}{q!}.
\end{aligned}
\end{equation*}
The first inequality follows from 
\begin{equation*}
\begin{aligned}
&\frac{k}{k+s}\rho^{-s} 
\left | \sum_{p \ge k+s} \frac{(k+s)!}{p!}(-\overline{\delta})^{p-k-s}
\right | \\
&\le \frac{k}{k+s} \cdot 1 \cdot \left | \sum_{j \ge 0} 
\frac{(k+s)!}{(k+s+j)!}(-\overline{\delta})^{j} \right | \\
&= \frac{k}{k+s} \left | \sum_{j \ge 0} (k+s)
\int_{0}^{1}(1-t)^{k+s-1}t^{j}\dd t
\frac{(-\overline{\delta})^{j}}{j!} \right | \\
&= k\int_{0}^{1}(1-t)^{k+s-1}e^{-\overline{\delta}t} \dd t \\
&\le k\int_{0}^{1}(1-t)^{k-1}e^{-\overline{\delta}t}\dd t \\
&= k\sum_{j \ge 0} \int_{0}^{1}(1-t)^{k-1}t^{j} \dd t 
\frac{(-\overline{\delta})^{j}}{j!} \\
&= \sum_{j \ge 0} \frac{k!}{(k+j)!}(-\overline{\delta})^{j}.
\end{aligned}
\end{equation*}
The second inequality follows from  
\begin{equation*}
\begin{aligned}
&\frac{k}{k+s}\sum_{q \ge k+s} (\rho\overline{\delta})^{q-k-s}
\frac{(k+s)!}{q!} \\
&\le
1 \cdot \sum_{j \ge 0} (\rho\overline{\delta})^{j}
\frac{(k+s)!}{(k+s+j)!} \\
& \le \sum_{j \ge 0}(\rho\overline{\delta})^{j}
\frac{k!}{(k+j)!} = \sum_{q \ge k}(\rho\overline{\delta})^{q-k}
\frac{k!}{q!}.
\end{aligned}
\end{equation*}
This completes the proof.
\end{proof}

\subsubsection{Collecting the bounds}
\label{subsubsec:orbit-column-bounds}

We now combine the estimates we obtained into a finite criterion for
$\|AE_{\mathrm{app}}\| \le Z_{1, \mathrm{app}}$.

\begin{proposition}[Bound for $\|AE_{\mathrm{app}}\|$]
\label{prop:orbit-derivative-error}
Let $P, \cS, \ell_{0}$ be as above, and let $y_{\delta}$ 
and $y_{m, \ell}$ be given by \eqref{eq:orbit-period-source} and
\eqref{eq:orbit-coefficient-source}, respectively.
For each finitely supported sequence $y$ appearing below, denote
by $\mathcal{T}(y)$ the right-hand side of
\eqref{eq:orbit-output-complement-bound}, with $M, n_{-}, n_{+}$
chosen for $y$ as in Lemma \ref{lem:orbit-output-complement}.
Suppose that $Z_{1, \mathrm{app}}>0$ satisfies
\begin{align}
\|A_{\mathrm{fin}}(P(0, L_{0}^{-1}y_{\delta}) - \overline{J}(1, 0))\|     
+\mathcal{T}(y_{\delta}) &\le Z_{1, \mathrm{app}}, 
\label{eq:orbit-period-check}\\
\bigg \| A_{\mathrm{fin}}\bigg [ 
P \left ( \overline{q}_{m}\mathbf{1}_{\{\ell=0\}}, 
L_{0}^{-1}y_{m, \ell} \right ) \qquad \qquad \qquad & 
\label{eq:orbit-coefficient-check} \\
+ (I_{PX}-\overline{J})P(0, e_{m, \ell})  
\bigg ] \bigg \| + 
\mathcal{T}(y_{m, \ell}) &\le Z_{1, \mathrm{app}} \quad (m \in \cS, 
\, 0 \le \ell < \ell_{0}), \nonumber \\
\sum_{k=\ell_{0}+1}^{\ell_{0}+K+1}\sum_{n=m-2}^{m+2}
\rho^{k}|(y_{m, \ell_{0}})_{n}^{(k)}| 
\left [ \|A_{\mathrm{fin}}P(0, L_{0}^{-1}e_{n, k})\| + 
\mathcal{T}(e_{n, k}) \right ] &\le Z_{1, \mathrm{app}} \quad 
(m \in \cS). \label{eq:orbit-high-Taylor-check}
\end{align}
Then
\[
\|AE_{\mathrm{app}}\|\le Z_{1, \mathrm{app}}.
\]
All quantities on the left-hand sides of
\eqref{eq:orbit-period-check}, \eqref{eq:orbit-coefficient-check}, 
and \eqref{eq:orbit-high-Taylor-check} can be bounded by finite
computations using Lemmas \ref{lem:orbit-finite-coordinates} and
\ref{lem:orbit-output-complement}.
\end{proposition}

\begin{proof}
By \eqref{eq:orbit-column-norm}, it suffices to show that
$\|AE_{\mathrm{app}}\mathbf{e}\|\le Z_{1, \mathrm{app}}$ for every 
input $\mathbf{e}=(1, 0)$ or $(0, e_{m, \ell})$. By 
\eqref{eq:orbit-column-splitting}, we have
\begin{equation} \label{eq:orbit-output-estimate}
\|AE_{\mathrm{app}}\mathbf{e}\| \le \|A_{\mathrm{fin}}
(PD\cG(\overline{x})\mathbf{e}-\overline{J}P\mathbf{e})\|
+ \|(I-P)(D\cG(\overline{x})\mathbf{e}-\mathbf{e})\|.
\end{equation}
We consider four cases.

\medskip
\noindent\textbf{Case 1}: $\mathbf{e}=(1, 0)$.
Recall that \[
D\cG(\overline{x})(1, 0) = (0, L_{0}^{-1}y_{\delta}).
\]
Since $P(1, 0)=(1, 0)$, \eqref{eq:orbit-output-estimate} gives \[
\|AE_{\mathrm{app}}(1, 0)\| \le \|A_{\mathrm{fin}} ( 
P(0, L_{0}^{-1}y_{\delta})-\overline{J}(1, 0) ) \| + \|(I-P)
(0, L_{0}^{-1}y_{\delta})\| \le Z_{1, \mathrm{app}},
\]
where the last inequality follows from Lemma 
\ref{lem:orbit-output-complement} and \eqref{eq:orbit-period-check}.

\medskip
\noindent\textbf{Case 2}: $\mathbf{e}=(0, e_{m, \ell})$ with
$m\in \cS$ and $0 \le \ell < \ell_{0}$.
Recall that \[
D\cG(\overline{x})\mathbf{e} = \left ( 
\sum_{n \in \ZZ}\overline{q}_{n}(e_{m, \ell})_{n}^{(0)},
e_{m, \ell} + L_{0}^{-1}y_{m, \ell}\right ).
\]
Plugging this into \eqref{eq:orbit-output-estimate} gives
\begin{equation*}
\begin{aligned}
\|AE_{\mathrm{app}}(0, e_{m, \ell})\| &\le \left \|
A_{\mathrm{fin}}\left [ P\left ( \overline{q}_{m}
\mathbf{1}_{\{\ell=0\}}, L_{0}^{-1}y_{m, \ell}\right ) 
+ (I_{PX}-\overline{J})P(0, e_{m, \ell})\right ]\right \| \\
&+ \|(I-P)(0, L_{0}^{-1}y_{m, \ell})\| \le Z_{1, \mathrm{app}},
\end{aligned}
\end{equation*}
where the last inequality follows from 
Lemma \ref{lem:orbit-output-complement} and 
\eqref{eq:orbit-coefficient-check}.

\medskip
\noindent\textbf{Case 3}: $\mathbf{e}=(0, e_{m, \ell})$ with
$m\in\cS$ and $\ell \ge \ell_{0}$.
Lemma \ref{lem:orbit-Taylor-comparison} gives
\begin{equation*}
\begin{aligned}
\|AE_{\mathrm{app}}(0, e_{m, \ell})\| &\le
\sum_{k=\ell_{0}+1}^{\ell_{0}+K+1}\sum_{n=m-2}^{m+2}\rho^{k}
|(y_{m, \ell_{0}})_{n}^{(k)}| 
\|A(0, L_{0}^{-1}e_{n, k})\| \le Z_{1, \mathrm{app}},
\end{aligned}
\end{equation*}
where the last inequality follows from Lemma 
\ref{lem:orbit-output-complement} and 
\eqref{eq:orbit-high-Taylor-check}, together with \[
\|A(0, L_{0}^{-1}e_{n, k})\|
\le \|A_{\mathrm{fin}}P(0, L_{0}^{-1}e_{n, k})\|      
+\|(I-P)(0, L_{0}^{-1}e_{n, k})\| 
\le \|A_{\mathrm{fin}}P(0, L_{0}^{-1}e_{n, k})\|    
+\mathcal{T}(e_{n, k}).
\]

\medskip
\noindent\textbf{Case 4}: $\mathbf{e}=(0, e_{m, \ell})$ with
$m\not\in \cS$.
By \eqref{eq:orbit-remote-columns}, \[
AE_{\mathrm{app}}(0, e_{m, \ell})=0 \quad (\ell \ge 0).
\]

These four cases complete the proof.
\end{proof}

\begin{remark}\label{rmk:orbit-verification}
All five estimates in Proposition \ref{prop:orbit-validation}
have now been reduced to finite computations.
It remains to specify
$\rho, \overline{\delta}, \overline{a}, \overline{q}, L, U, K,
\overline{J}, A_{\mathrm{fin}}, \cS, \ell_{0}$
as above, and certify
the following:
\begin{itemize}
\item The parameter conditions
$\rho>1$ and $\overline{\delta}>0$.
\item The finite matrix bound
$\|I_{PX}-A_{\mathrm{fin}}\overline{J}\|\le Z_{1, \mathrm{inv}}$.
\item The right-hand side of \eqref{eq:orbit-residual-bound}
is at most $Z_{0}$.
\item The right-hand sides of
\eqref{eq:orbit-higher-derivative-bounds}
are at most $Z_{2}$ and $Z_{3}$, respectively.
\item The inequalities \eqref{eq:orbit-period-check}, 
\eqref{eq:orbit-coefficient-check}, 
\eqref{eq:orbit-high-Taylor-check} in
Proposition \ref{prop:orbit-derivative-error}.
\item A radius $r>0$ satisfies
\eqref{eq:orbit-radii-inequality} and
\eqref{eq:orbit-radius-restrictions}.
\end{itemize}
These checks will be carried out in the next subsection.
\end{remark}

\subsection{Computer-assisted verification and completion of the proof}
\label{subsec:orbit-finite-verification}

In this subsection, we describe the computer-assisted verification 
needed to complete the proof of Theorem \ref{thm:orbit-existence}.
All code and input data used here are available in the supplementary
material and in the GitHub repository \cite{CAParchive}.

We first specify the parameters and numerical data listed in Remark 
\ref{rmk:orbit-verification}.

For $\overline{\delta}$ and $\overline{a}$, we use the numerical 
approximation constructed in Appendix \ref{app:numerical-approximation}.
These data are stored in \texttt{approximation.npz} in the directory 
\texttt{codes/appendix\_a}, and satisfy \[
\overline{\delta} \approx 0.49840771506392201, \quad I=[-51, 21], 
\quad K=19,
\]
with $\overline{a}_{n}^{(k)}=0$ for all $(n, k)$ outside
$I \times [0, K]$.

We use \[
\rho = \frac{101}{100}, \quad L=-20, \quad U=8,
\]
so that $\dim PX=29 \times 20 + 1 = 581$.

For $\overline{q}$ in the definition of $\cF$, recall that its
purpose is to fix the temporal phase freedom arising from time translation.
It is therefore natural to choose $\overline{q}$ in the direction
of a tangent to the orbit. We take 
$-\overline{a}^{(0)}+N(\overline{a}^{(0)})$, restrict it to $[L, U]$,
and numerically normalize it in $\ell^{2}$.

Next, we choose $\overline{J}$ and $A_{\mathrm{fin}}$, which are 
intended to be numerical approximations to 
$PD\cG(\overline{x})|_{PX}$ and its inverse. To compute 
$\overline{J}$, we use \eqref{eq:orbit-derivative-output} to obtain the 
columns of $PD\cG(\overline{x})|_{PX}$, and the 
finite-coordinate formulas in Lemma \ref{lem:orbit-finite-coordinates}
are used to evaluate these columns numerically. Then $A_{\mathrm{fin}}$ is 
computed as a numerical inverse of $\overline{J}$. These computations 
use ordinary floating-point arithmetic.

Finally, we choose \[
\cS = [-53, 23], \quad \ell_{0}=39, 
\]
so that \eqref{eq:orbit-spatial-input-set} and $\ell_{0}>K$ are satisfied.

The quantities $\overline{q}, \overline{J}$, and $A_{\mathrm{fin}}$
are computed in \texttt{prepare\_data.py}. They are stored in
\texttt{data.npz} in the directory \texttt{codes/section\_2},
together with $\rho, L, U, \cS$, and $\ell_{0}$.
This file also stores a relative path to the approximation file
from Appendix \ref{app:numerical-approximation}, from which
$I, K, \overline{\delta}$, and $\overline{a}$ are read directly.

We now verify the conditions listed in Remark
\ref{rmk:orbit-verification}. The resulting bounds are recorded 
in the following lemma.

\begin{lemma}[Verified bounds]\label{lem:orbit-finite-verification}
For the data specified above, the estimates
\eqref{eq:orbit-validation-estimates} hold with
\begin{equation}\label{eq:orbit-validation-constants}
Z_{0}=10^{-6}, \quad Z_{1, \mathrm{inv}}=10^{-13}, \quad
Z_{1, \mathrm{app}}=0.52, \quad Z_{2}=3200, \quad 
Z_{3}=260.
\end{equation}
Moreover, $\overline{\delta}>0$, and $r=10^{-5}$ satisfies 
\eqref{eq:orbit-radii-inequality} and
\eqref{eq:orbit-radius-restrictions}.
\end{lemma}

\begin{proof}
The script \texttt{verify\_orbit.py} verifies the conditions listed in Remark
\ref{rmk:orbit-verification} in the order given there. Interval 
computations are performed with \texttt{python-flint}, which provides 
Python access to FLINT and Arb; see 
\cite{pythonflint, Johansson2017arb}. We use a working precision of 
$64$ bits. 

The function \texttt{check\_parameters} verifies $\rho>1$ and 
$\overline{\delta}>0$. It also checks that $\overline{\delta}<1$,
which is needed in Remark \ref{rmk:orbit-coefficient-evaluation}.

The function \texttt{check\_z1\_inv} computes an interval enclosure of
$I_{PX}-A_{\mathrm{fin}}\overline{J}$. In the normalized basis of 
$PX$, the operator norm is the maximum absolute column sum. The 
function evaluates this quantity using interval arithmetic and returns 
the upper endpoint of the resulting scalar enclosure, denoted by 
$\widehat{Z}_{1, \mathrm{inv}}$. It also checks that
$\widehat{Z}_{1, \mathrm{inv}} \le Z_{1, \mathrm{inv}}$.

The function \texttt{check\_z0} evaluates the right-hand side of
\eqref{eq:orbit-residual-bound}, with $\|L_{0}^{-1}\|$ replaced
by its upper bound in \eqref{eq:L0-inverse-bound}. $F(\overline{x})$ 
is finitely supported, so its weighted norm is evaluated by a finite 
sum using interval arithmetic. The upper endpoint of the enclosure 
of the full right-hand side is returned as $\widehat{Z}_{0}$, and 
the function checks that $\widehat{Z}_{0}\le Z_{0}$.

The function \texttt{check\_z2\_z3} similarly evaluates the two
right-hand sides of \eqref{eq:orbit-higher-derivative-bounds}.
It returns the upper endpoints of the resulting enclosures, denoted 
by $\widehat{Z}_{2}$ and $\widehat{Z}_{3}$, and checks that 
$\widehat{Z}_{2}\le Z_{2}$ and $\widehat{Z}_{3}\le Z_{3}$.

The function \texttt{check\_z1\_app} computes rigorous upper bounds for the 
left-hand sides of \eqref{eq:orbit-period-check},
\eqref{eq:orbit-coefficient-check}, and
\eqref{eq:orbit-high-Taylor-check}, over all the indices specified 
there. For the first two inequalities, $y_{\delta}$ and $y_{m, \ell}$ 
are computed using \eqref{eq:orbit-period-source} and
\eqref{eq:orbit-coefficient-source}. The required finitely many 
coordinates of $L_{0}^{-1}y_{\delta}$ and $L_{0}^{-1}y_{m, \ell}$ 
are enclosed using Lemma \ref{lem:orbit-finite-coordinates}
(see also Remark \ref{rmk:orbit-coefficient-evaluation}). 
$\mathcal{T}(y_{\delta})$ and $\mathcal{T}(y_{m, \ell})$ are 
evaluated using \eqref{eq:orbit-output-complement-bound}. The remaining 
matrix computations are also performed using interval arithmetic.
For \eqref{eq:orbit-high-Taylor-check}, the same computations are 
applied to $e_{n, k}$. The function returns the upper endpoint of an 
enclosure of the maximum over all three families of checks, denoted 
by $\widehat{Z}_{1, \mathrm{app}}$, and verifies that
$\widehat{Z}_{1, \mathrm{app}}\le Z_{1, \mathrm{app}}$.

Finally, the function \texttt{check\_radius} uses the prescribed
constants from \eqref{eq:orbit-validation-constants} to verify
\eqref{eq:orbit-radii-inequality} and
\eqref{eq:orbit-radius-restrictions} for $r=10^{-5}$.

The script \texttt{verify\_orbit.py} runs these functions in the order 
described above and records in \texttt{results\_orbit.txt} the returned 
upper bounds, the prescribed constants, and whether each required 
inequality was verified.

Table \ref{tab:orbit-verification} summarizes the output and the 
time spent on each verification step. The times were measured on an 
Intel Core Ultra 5 125H processor. On the same machine, 
\texttt{prepare\_data.py} generated \texttt{data.npz} in approximately 
$3$ seconds.

\begin{table}[htbp]
\centering
\small
\begin{tabular}{lccc}
\hline
Function & Returned upper bound & Prescribed constant & Time (s) \\
\hline
\texttt{check\_parameters} & --- & --- & $<0.001$ \\
\texttt{check\_z1\_inv} & $6.801 \times 10^{-15}$ & $10^{-13}$ & $\approx 1.120$ \\
\texttt{check\_z0} & $5.857 \times 10^{-7}$ & $10^{-6}$ & $\approx 0.288$ \\
\texttt{check\_z2\_z3} & $(3138, 259.6)$ & $(3200, 260)$ & $\approx 0.141$ \\
\texttt{check\_z1\_app} & $0.5110$ & $0.52$ & $\approx 82.804$ \\
\texttt{check\_radius} & --- & --- & $<0.001$ \\
\hline
Total & & & $\approx 85.249$ \\
\hline
\end{tabular}
\caption{Returned upper bounds and computation times for the 
verification. The returned upper bounds are rounded upward for 
display. They still satisfy the comparisons with the prescribed 
constants in \eqref{eq:orbit-validation-constants}. Times are from 
one recorded run and are rounded to the nearest millisecond. The total 
time refers to \texttt{verify\_orbit.py} and includes reading the 
input data. Timings may vary between runs.}
\label{tab:orbit-verification}
\end{table}
\end{proof}

\begin{proof}[Proof of Theorem \ref{thm:orbit-existence}]
By Lemma \ref{lem:orbit-finite-verification}, all the assumptions 
of Proposition \ref{prop:orbit-validation} are satisfied. Therefore 
there exists a zero $(\delta^{(\ast)}, a^{(\ast)})$ of $\cG$ 
satisfying \[
|\delta^{(\ast)}-\overline{\delta}| + 
\|a^{(\ast)}-\overline{a}\|_{\ell_{\rho}^{1}} \le 10^{-5},
\] with $\delta^{(\ast)} > 0$ and $a^{(\ast)} \not\equiv 0$.
By Lemma \ref{lem:L0-inverse}, this is also a zero of $\cF$.

The series \eqref{eq:bilateral-orbit-ansatz} therefore
defines a solution $Z^{(\ast)}\not\equiv 0$ of \eqref{eq:main-ss-Z}
on $[0, \delta^{(\ast)}]$ satisfying \[
Z^{(\ast)}(\delta^{(\ast)})=RZ^{(\ast)}(0).
\]
Since $\rho>1$, the series and its first time derivative converge 
in $\ell^{1}$. By shift equivariance and uniqueness of the 
solution to \eqref{eq:main-ss-Z}, this extends to a solution on 
$[0, \infty)$ satisfying \[
Z^{(\ast)}(\tau+\delta^{(\ast)})=RZ^{(\ast)}(\tau) \quad (\tau \ge 0).
\]
The spatial decay asserted in the theorem follows from Subsection
\ref{subsec:orbit-spatial-decay}. This completes the proof.
\end{proof}

\section{Linearized dynamics near the travelling wave}
\label{sec:monodromy}

Throughout this section, $N$ denotes either $N_{\mathrm{GOY}}$
or $N_{\mathrm{Sabra}}$. Unless a model is specified, all
statements apply to both choices of $N$. $\ell^{\infty}$ denotes
$\ell^{\infty}(\ZZ; \CC)$, which we regard as a Banach space over
$\RR$, so that all linear maps, derivatives, and spans are
over $\RR$ rather than $\CC$.

Let $Z^{(\ast)}$ and $\delta^{(\ast)}$ be the bilateral travelling wave and 
relative period constructed in the proof of Theorem \ref{thm:orbit-existence}.
The GOY and Sabra equations share the real subsystem used in this
construction, so $Z^{(\ast)}$ is a travelling wave for
\eqref{eq:main-ss-Z} in both models. Thus, for either model,
$Z^{(\ast)}(0)$ is a fixed point of the bilateral relative return map
$R^{-1}\Phi^{\delta^{(\ast)}}$, where $\Phi^{s}$ denotes the 
corresponding local flow of \eqref{eq:main-ss-Z} on $\ell^{\infty}$.
In this section, we study the linear stability of this fixed point
$Z^{(\ast)}(0)$. That is, we study the linearized bilateral relative
return map \[
M:= R^{-1}D\Phi^{\delta^{(\ast)}}(Z^{(\ast)}(0)).
\]

This linearized map $M$ is not contractive, since the time-translation and 
phase symmetries give nonstable modes, as we verify in Subsection 
\ref{subsec:monodromy-modes}. The corresponding modes span the space
\begin{equation}\label{eq:monodromy-V}
V:= \begin{cases}
\mathrm{span}\{ v_{\mathrm{r}}^{u}, v_{\mathrm{r}}^{c}, 
v_{\mathrm{G}, 1}^{c}, v_{\mathrm{G}, 2}^{c} \} & \text{for GOY}, \\
\mathrm{span}\{ v_{\mathrm{r}}^{u}, v_{\mathrm{r}}^{c},
v_{\mathrm{S}}^{u} \} & \text{for Sabra},
\end{cases}
\end{equation}
where 
\begin{equation}\label{eq:monodromy-modes}
\begin{aligned}
v_{\mathrm{r}}^{u} &:= N(Z^{(\ast)}(0)), \\
v_{\mathrm{r}}^{c} &:= -Z^{(\ast)}(0) + N(Z^{(\ast)}(0)), \\
(v_{\mathrm{G}, 1}^{c})_{n} &:= 
\begin{cases}
-i(Z^{(\ast)}(0))_{n} & n \equiv 0 \pmod{3}, \\
i(Z^{(\ast)}(0))_{n} & n \equiv 1 \pmod{3}, \\
0 & n \equiv 2 \pmod{3},
\end{cases} \\
(v_{\mathrm{G}, 2}^{c})_{n} &:= 
\begin{cases}
-i(Z^{(\ast)}(0))_{n} & n \equiv 0 \pmod{3}, \\
0 & n \equiv 1 \pmod{3}, \\
i(Z^{(\ast)}(0))_{n} & n \equiv 2 \pmod{3},
\end{cases} \\
(v_{\mathrm{S}}^{u})_{n} &:= i\left ( \frac{1+\sqrt{5}}{2} \right )^{n} 
(Z^{(\ast)}(0))_{n}.
\end{aligned}
\end{equation}
In the modulation argument later, we will remove these nonstable modes
by applying suitable symmetries. Formally, this modulation amounts
to projecting onto a complement of $V$. More precisely, for a bounded linear 
projection $\Pi: \ell^{\infty} \to V$, set $Y:= \ker\Pi$. Then the 
modulation corresponds to applying $I-\Pi$ to the error. Therefore, 
instead of studying $M$, we study the modulated operator 
\begin{equation}\label{eq:modulated-operator}
\widetilde{M}:=(I-\Pi)M|_{Y}.
\end{equation}
See Section \ref{sec:boundary} for more details on why this operator 
appears naturally. The goal of this section is to prove
that we may choose $\Pi$ so that some power of $\widetilde{M}:Y \to Y$ is 
contractive, where the norm on $Y$ is the restriction of the standard 
norm on $\ell^{\infty}$.

\begin{theorem} \label{thm:monodromy-splitting}
Let $V$ be as in \eqref{eq:monodromy-V}. Then $V$ is $M$-invariant, 
and the generators listed in \eqref{eq:monodromy-V} are linearly 
independent in each model, GOY or Sabra.

Also, there is a topological direct sum decomposition \[
\ell^{\infty}(\ZZ; \CC) = Y \oplus V
\]
and a corresponding bounded linear projection 
$\Pi: \ell^{\infty}(\ZZ; \CC) \to V$ with $\ker \Pi = Y$, such that
\begin{equation}\label{eq:modulated-operator-contractivity}
\|\widetilde{M}^{p}\|_{Y \to Y} < 1
\end{equation}
for some integer $p \ge 1$, where $\widetilde{M}=(I-\Pi)M|_{Y}$ and 
the norm on $Y$ is the restriction of the standard norm on 
$\ell^{\infty}(\ZZ; \CC)$. Here, $\Pi v$ depends only on a fixed 
finite set of coordinates of $v \in \ell^{\infty}(\ZZ; \CC)$.
\end{theorem}

The rest of this section is devoted to the proof of Theorem 
\ref{thm:monodromy-splitting}. In Subsection 
\ref{subsec:monodromy-modes}, we use the variational equation to 
compute the action of $M$ on the modes in \eqref{eq:monodromy-modes}.
In Subsection \ref{subsec:monodromy-criterion}, we choose the projection 
$\Pi$ used in the proof of Theorem \ref{thm:monodromy-splitting} and 
give a contraction criterion. To verify the required bounds, we construct 
an approximation $\overline{\widetilde{M}}$ to $\widetilde{M}$ in
Subsection \ref{subsec:monodromy-estimates} and reduce the estimates 
to finite computations. Finally, in Subsection 
\ref{subsec:monodromy-verification}, we carry out these finite 
computations using interval arithmetic.

\subsection{The variational equation and the nonstable modes}
\label{subsec:monodromy-modes}

We first compute the action of $M$ on the modes in \eqref{eq:monodromy-modes}.
These modes arise from time-translation and phase symmetries. 
We explain their origin and then verify the claimed action using the 
variational equation.

\subsubsection{Time-translation symmetries}

The common unstable mode and the center mode, $v_{\mathrm{r}}^{u}$ and 
$v_{\mathrm{r}}^{c}$, reflect two time-translation symmetries. Translating 
the self-similar time changes the temporal phase of the travelling wave 
but preserves its relative period. Thus the tangent vector to the 
bilateral travelling wave $\tau \mapsto Z^{(\ast)}(\tau)$ 
at $\tau=0$, namely, $-Z^{(\ast)}(0)+N(Z^{(\ast)}(0))$, corresponds 
to the center mode.

The unstable mode comes from translating the original time. Returning 
to the original time variable, write \[
z^{(\ast)}(t) := \frac{1}{1-t}Z^{(\ast)}(-\log(1-t)) \quad (0 \le t < 1),
\]
which solves $\partial_{t}z^{(\ast)}=N(z^{(\ast)})$. For small 
$\epsilon>0$, the translated solution $z^{(\ast)}(t+\epsilon)$ has 
blowup time $1-\epsilon$. Returning to self-similar variables centered 
at $t=1$, we obtain another solution of \eqref{eq:main-ss-Z} for 
$0 \le \tau <-\log\epsilon$: \[
Z_{\epsilon}(\tau) := e^{-\tau}z^{(\ast)}(1-e^{-\tau}+\epsilon) 
= \frac{1}{1-\epsilon e^{\tau}} Z^{(\ast)}
(\tau-\log(1-\epsilon e^{\tau})).
\]
At $\tau=0$, the initial data satisfy
\begin{align*}
Z_{\epsilon}(0) &= \frac{1}{1-\epsilon}Z^{(\ast)}(-\log(1-\epsilon)) \\
&\approx (1+\epsilon)(Z^{(\ast)}(0) + \epsilon \partial_{\tau} 
Z^{(\ast)}(0)) \\
&\approx Z^{(\ast)}(0) + \epsilon N(Z^{(\ast)}(0)).
\end{align*}
Taking the bilateral relative return map to these initial data gives 
\begin{align*}
R^{-1}\Phi^{\delta^{(\ast)}}(Z_{\epsilon}(0))
&= R^{-1}Z_{\epsilon}(\delta^{(\ast)}) \\
&= R^{-1}\frac{1}{1-\epsilon e^{\delta^{(\ast)}}}Z^{(\ast)}
(\delta^{(\ast)} - \log(1-\epsilon e^{\delta^{(\ast)}})) \\
&= \frac{1}{1-\epsilon e^{\delta^{(\ast)}}}Z^{(\ast)}
(- \log(1-\epsilon e^{\delta^{(\ast)}})) \\
&\approx Z^{(\ast)}(0) + \epsilon e^{\delta^{(\ast)}}N(Z^{(\ast)}(0)).
\end{align*}
Thus the direction $N(Z^{(\ast)}(0))$ is amplified by the factor
$e^{\delta^{(\ast)}}$ under the linearized bilateral relative return map.

We now verify the claimed action of $M$ on these two modes using the 
variational equation. Since $-I+N$ is a $C^{1}$ vector field on 
$\ell^{\infty}$, the derivative of its flow satisfies the variational 
equation:
\begin{equation}\label{eq:variational-ODE}
\partial_{s}D\Phi^{s}(Z^{(\ast)}(0)) = \left ( 
-I + DN(Z^{(\ast)}(s)) \right ) D\Phi^{s}(Z^{(\ast)}(0)), \quad 
D\Phi^{0}(Z^{(\ast)}(0)) = I.
\end{equation}
On the other hand, differentiating \[
\partial_{s}Z^{(\ast)}(s)=-Z^{(\ast)}(s)+N(Z^{(\ast)}(s))
\]
with respect to $s$ gives
\begin{equation}\label{eq:variational-ODE-center}
\partial_{s}\left ( \partial_{s}Z^{(\ast)}(s)\right )        
=\left ( -I + DN(Z^{(\ast)}(s)) \right ) \partial_{s}Z^{(\ast)}(s),
\quad \partial_{s}Z^{(\ast)}(0) = -Z^{(\ast)}(0)+N(Z^{(\ast)}(0)).
\end{equation}
Recall that $v_{\mathrm{r}}^{c} = -Z^{(\ast)}(0)+N(Z^{(\ast)}(0))$. By 
\eqref{eq:variational-ODE} and \eqref{eq:variational-ODE-center}, 
the functions $s \mapsto D\Phi^{s}(Z^{(\ast)}(0))v_{\mathrm{r}}^{c}$ and
$s \mapsto \partial_{s}Z^{(\ast)}(s)$ satisfy the same variational
equation with initial value $v_{\mathrm{r}}^{c}$. They therefore coincide by 
uniqueness. Using relative periodicity of $Z^{(\ast)}$, we conclude that
\begin{align*}
Mv_{\mathrm{r}}^{c} &= R^{-1}D\Phi^{\delta^{(\ast)}}(Z^{(\ast)}(0))
v_{\mathrm{r}}^{c} \\
&= R^{-1}\partial_{s}Z^{(\ast)}(\delta^{(\ast)}) \\
&= R^{-1}\left ( -Z^{(\ast)}(\delta^{(\ast)}) + 
N(Z^{(\ast)}(\delta^{(\ast)})) \right ) \\
&= -Z^{(\ast)}(0)+N(Z^{(\ast)}(0)) = v_{\mathrm{r}}^{c}.
\end{align*}

For the unstable mode, observe that \[
DN(z)z = \lim_{h \to 0} \frac{N(z+hz)-N(z)}{h} 
= \lim_{h \to 0} \frac{(1+h)^{2}N(z) - N(z)}{h} = 2N(z), 
\]
and hence
\begin{align}
\frac{\dd}{\dd s} \left ( e^{s}N(Z^{(\ast)}(s))\right )
&= e^{s} \left [ N(Z^{(\ast)}(s)) + DN(Z^{(\ast)}(s)) 
\left ( -Z^{(\ast)}(s) + N(Z^{(\ast)}(s)) \right ) \right ] 
\nonumber \\
&= e^{s} \left [ - N(Z^{(\ast)}(s)) + DN(Z^{(\ast)}(s))N(Z^{(\ast)}(s)) 
\right ] \nonumber \\
&= \left ( -I+DN(Z^{(\ast)}(s))\right )e^{s}N(Z^{(\ast)}(s)).
\label{eq:variational-ODE-unstable}
\end{align}
Recall that $v_{\mathrm{r}}^{u} = N(Z^{(\ast)}(0))$. Comparing 
\eqref{eq:variational-ODE-unstable} and \eqref{eq:variational-ODE} 
and using uniqueness as above, we obtain \[
D\Phi^{s}(Z^{(\ast)}(0))v_{\mathrm{r}}^{u} = e^{s}N(Z^{(\ast)}(s)).
\]
Thus we conclude that
\begin{align*}
Mv_{\mathrm{r}}^{u} &= R^{-1}D\Phi^{\delta^{(\ast)}}(Z^{(\ast)}(0))
v_{\mathrm{r}}^{u} \\
&= e^{\delta^{(\ast)}}R^{-1}N(Z^{(\ast)}(\delta^{(\ast)})) \\
&= e^{\delta^{(\ast)}}N(Z^{(\ast)}(0)) 
= e^{\delta^{(\ast)}}v_{\mathrm{r}}^{u}.
\end{align*}

\subsubsection{Phase symmetries}

The remaining modes $v_{\mathrm{G}, 1}^{c}$, $v_{\mathrm{G}, 2}^{c}$
or $v_{\mathrm{S}}^{u}$ come from phase symmetries described in
\cite[Section III.B]{LvovEtAl1998}. We consider the two models
separately.

For the GOY model, define
\begin{equation}\label{eq:GOYphase}
(T_{\theta_{1}, \theta_{2}}Z)_{n} :=
\begin{cases}
e^{-i(\theta_{1}+\theta_{2})}Z_{n} & n \equiv 0 \pmod{3}, \\
e^{i\theta_{1}}Z_{n} & n \equiv 1 \pmod{3}, \\
e^{i\theta_{2}}Z_{n} & n \equiv 2 \pmod{3},
\end{cases}
\quad (\theta_{1}, \theta_{2}) \in \RR^{2}.
\end{equation}
Substitution into the nonlinearity gives
$N_{\mathrm{GOY}}(T_{\theta_{1}, \theta_{2}}Z)
=T_{\theta_{1}, \theta_{2}}N_{\mathrm{GOY}}(Z)$,
and therefore if $Z$ solves \eqref{eq:main-ss-Z} for
$N=N_{\mathrm{GOY}}$, then $T_{\theta_{1}, \theta_{2}}Z$ also
solves the same equation. Differentiating the transformed
profile with respect to $\theta_{1}$ and $\theta_{2}$ at 
$(\theta_{1}, \theta_{2})=(0, 0)$ gives the directions
$v_{\mathrm{G}, 1}^{c}$ and $v_{\mathrm{G}, 2}^{c}$, respectively.
Since \[
R^{-1}T_{\theta_{1}, \theta_{2}}R = T_{\theta_{2}, -\theta_{1}-\theta_{2}},
\]
we expect $M$ to send the first direction to $-v_{\mathrm{G}, 2}^{c}$ 
and the second to $v_{\mathrm{G}, 1}^{c}-v_{\mathrm{G}, 2}^{c}$.

We now verify this using the variational equation. Define 
$w_{\mathrm{G}, 1}(s)$ by the derivative of $T_{\theta, 0}Z^{(\ast)}(s)$ with 
respect to $\theta$ at $0$, \[
(w_{\mathrm{G}, 1}(s))_{n} := 
\begin{cases}
-i(Z^{(\ast)}(s))_{n} & n \equiv 0 \pmod{3}, \\
i(Z^{(\ast)}(s))_{n} & n \equiv 1 \pmod{3}, \\
0 & n \equiv 2 \pmod{3}.
\end{cases}
\] Then
\begin{align*}
(\partial_{s}w_{\mathrm{G}, 1}(s))_{n} &= 
\begin{cases}
-i\partial_{s}(Z^{(\ast)}(s))_{n} & n \equiv 0 \pmod{3}, \\
i\partial_{s}(Z^{(\ast)}(s))_{n} & n \equiv 1 \pmod{3}, \\
0 & n \equiv 2 \pmod{3}
\end{cases} \\
&=
\begin{cases}
-i(-(Z^{(\ast)}(s))_{n}+(N_{\mathrm{GOY}}(Z^{(\ast)}(s)))_{n}) 
& n \equiv 0 \pmod{3}, \\
i(-(Z^{(\ast)}(s))_{n}+(N_{\mathrm{GOY}}(Z^{(\ast)}(s)))_{n})
& n \equiv 1 \pmod{3}, \\
0 & n \equiv 2 \pmod{3}.
\end{cases}
\end{align*}
Therefore, differentiating \[
N_{\mathrm{GOY}}(T_{\theta, 0}Z^{(\ast)}(s)) 
= T_{\theta, 0}N_{\mathrm{GOY}}(Z^{(\ast)}(s))
\] with respect to $\theta$ at $0$ gives \[
DN_{\mathrm{GOY}}(Z^{(\ast)}(s))w_{\mathrm{G}, 1}(s) 
= w_{\mathrm{G}, 1}(s) + \partial_{s}w_{\mathrm{G}, 1}(s).
\]
Thus $w_{\mathrm{G}, 1}(s)$ solves \[
\partial_{s}w_{\mathrm{G}, 1}(s) 
= (-I + DN_{\mathrm{GOY}}(Z^{(\ast)}(s)))w_{\mathrm{G}, 1}(s), \quad 
w_{\mathrm{G}, 1}(0) = v_{\mathrm{G}, 1}^{c}.
\]
By \eqref{eq:variational-ODE} and uniqueness, we obtain \[
D\Phi_{\mathrm{GOY}}^{s}(Z^{(\ast)}(0))v_{\mathrm{G}, 1}^{c} 
= w_{\mathrm{G}, 1}(s). 
\]
Relative periodicity then gives 
\begin{align*}
(Mv_{\mathrm{G}, 1}^{c})_{n} 
&= (R^{-1}w_{\mathrm{G}, 1}(\delta^{(\ast)}))_{n} \\
&= \begin{cases}
-i(Z^{(\ast)}(\delta^{(\ast)}))_{n+1} & n+1 \equiv 0 \pmod{3}, \\
i(Z^{(\ast)}(\delta^{(\ast)}))_{n+1} & n+1 \equiv 1 \pmod{3}, \\
0 & n+1 \equiv 2 \pmod{3}
\end{cases} \\
&= \begin{cases}
-i(Z^{(\ast)}(0))_{n} & n \equiv 2 \pmod{3}, \\
i(Z^{(\ast)}(0))_{n} & n \equiv 0 \pmod{3}, \\
0 & n \equiv 1 \pmod{3}
\end{cases} 
= -(v_{\mathrm{G}, 2}^{c})_{n}.
\end{align*}
A similar argument shows that 
$Mv_{\mathrm{G}, 2}^{c} = v_{\mathrm{G}, 1}^{c} 
- v_{\mathrm{G}, 2}^{c}$.

For the Sabra model, define
\begin{equation}\label{eq:Sabraphase}
(T_{\theta}Z)_{n} := e^{i\theta \left ( \frac{1+\sqrt{5}}{2}\right )^{n}}
Z_{n} \quad (\theta \in \RR).
\end{equation}
Substitution into the nonlinearity gives
$N_{\mathrm{Sabra}}(T_{\theta}Z)=T_{\theta}N_{\mathrm{Sabra}}(Z)$.
Define $w_{\mathrm{S}}(s)$ by the formal derivative of
$T_{\theta}Z^{(\ast)}(s)$ with respect to $\theta$ at $0$, \[
(w_{\mathrm{S}}(s))_{n} := i\left ( \frac{1+\sqrt{5}}{2}\right )^{n} 
(Z^{(\ast)}(s))_{n}.
\]
One should be careful since the derivative of $T_{\theta}$ with respect to 
$\theta$ is not well-defined on $\ell^{\infty}$; 
$\theta \mapsto T_{\theta}$ is not even continuous in the operator norm. 
Still, $w_{\mathrm{S}}$ belongs to
$C^{1}([0, \delta^{(\ast)}]; \ell^{\infty})$, thanks to the right spatial 
decay of $Z^{(\ast)}$ proved in Theorem \ref{thm:orbit-existence}. 
Then, differentiating each component of \[
N_{\mathrm{Sabra}}(T_{\theta}Z^{(\ast)}(s)) 
= T_{\theta}N_{\mathrm{Sabra}}(Z^{(\ast)}(s))
\] with respect to $\theta$ at $0$ and arguing as in the GOY case gives \[
D\Phi_{\mathrm{Sabra}}^{s}(Z^{(\ast)}(0)) v_{\mathrm{S}}^{u} 
= w_{\mathrm{S}}(s),
\]
and relative periodicity again gives 
$Mv_{\mathrm{S}}^{u} = \frac{1+\sqrt{5}}{2}v_{\mathrm{S}}^{u}$.

These computations show that $V$ is $M$-invariant in either
model. The required linear independence of its generators
will follow from the invertibility condition introduced in the
next subsection.

\subsection{The projection and the contraction criterion}
\label{subsec:monodromy-criterion}

For the computation of the matrix norms, it is convenient to identify 
$\ell^{\infty} = \ell^{\infty}(\ZZ; \CC)$ with
$\ell^{\infty}(\ZZ \times \{1, 2\}; \RR)$. The real and imaginary parts 
of the $n$th component $z_{n}$ correspond to the coordinates $(n, 1)$ 
and $(n, 2)$, respectively: \[
z_{n, 1}:=\mathrm{Re}z_{n}, \quad z_{n, 2}:=\mathrm{Im}z_{n}.
\]
We denote the corresponding supremum norm by
\begin{equation}\label{eq:monodromy-maximum-norm}
\|z\|_{\max} := \sup_{\substack{n \in \ZZ \\ a \in \{1, 2\}}}
|z_{n, a}| = \sup_{n \in \ZZ} \max \{ |\mathrm{Re}z_{n}|, 
|\mathrm{Im}z_{n}| \}.
\end{equation}
This norm is equivalent to the standard norm on
$\ell^{\infty}(\ZZ; \CC)$, since \[
\|z\|_{\max} \le \|z\|_{\ell^{\infty}} \le \sqrt{2} \|z\|_{\max}.
\]
For the rest of this section, finite-dimensional spaces
$\RR^{d}$ carry their maximum norms. Unless otherwise specified,
operator norms are induced by these norms and by
$\|\cdot\|_{\max}$ and its restrictions.

We now choose the projection $\Pi$ used in the proof of Theorem
\ref{thm:monodromy-splitting}.

Let $m=4$ for GOY and $m=3$ for Sabra. Define
$\cV: \RR^{m} \to \ell^{\infty}$ by
\begin{equation}\label{eq:monodromy-mode-map}
\cV\xi := 
\begin{cases}
\xi_{1}v_{\mathrm{r}}^{u} + \xi_{2}v_{\mathrm{r}}^{c}
+\xi_{3}v_{\mathrm{G}, 1}^{c} + \xi_{4}v_{\mathrm{G}, 2}^{c}
& \text{for GOY}, \\
\xi_{1}v_{\mathrm{r}}^{u} + \xi_{2}v_{\mathrm{r}}^{c}
+ \xi_{3}v_{\mathrm{S}}^{u} & \text{for Sabra}
\end{cases} \quad (\xi \in \RR^{m}).
\end{equation}
Thus $\mathrm{ran}\cV = V$. To choose its complement, define the 
coordinate map $\pi: \ell^{\infty} \to \RR^{m}$ by
\begin{equation}\label{eq:monodromy-coordinate-map}
\pi z :=
\begin{cases}
(\mathrm{Re}z_{n_{1}}, \mathrm{Re}z_{n_{2}}, \mathrm{Im}z_{n_{3}}, 
\mathrm{Im}z_{n_{4}}) & \text{for GOY}, \\
(\mathrm{Re}z_{n_{1}}, \mathrm{Re}z_{n_{2}}, \mathrm{Im}z_{n_{3}}) 
& \text{for Sabra},
\end{cases}
\end{equation}
where the sampling indices $n_{i} (1 \le i \le m)$ will be specified 
later. We require $n_{1} \neq n_{2}$, and also $n_{3} \neq n_{4}$ 
for GOY. The two modes $v_{\mathrm{r}}^{u}, v_{\mathrm{r}}^{c}$ are 
real, while $v_{\mathrm{G}, 1}^{c}, v_{\mathrm{G}, 2}^{c}$ and
$v_{\mathrm{S}}^{u}$ are purely imaginary. This is why we choose two 
real coordinates and two imaginary coordinates for GOY, and two real 
coordinates and one imaginary coordinate for Sabra.

Assume that $\pi\cV: \RR^{m} \to \RR^{m}$ is invertible; this will be 
checked later. Note that this invertibility condition implies that the 
generators of $V$ in \eqref{eq:monodromy-V} are linearly independent. 
Define
\begin{equation}\label{eq:monodromy-projection}
\Pi := \cV(\pi\cV)^{-1}\pi.
\end{equation}
Then $\Pi$ is a bounded projection onto $V$, and $\Pi z$ depends only 
on the coordinates selected by $\pi$. Its kernel is $Y=\ker\pi$, so we 
obtain the decomposition \[
\ell^{\infty} = \mathrm{ran}\Pi \oplus \ker\Pi = V \oplus Y.
\]
For this choice of $\Pi$, we define $\widetilde{M}$ as in 
\eqref{eq:modulated-operator}, \[
\widetilde{M} := (I-\Pi)M|_{Y}: Y \to Y.
\]

Since $V$ is $M$-invariant, $(I-\Pi)M\Pi=0$, and hence \[
\widetilde{M}^{k} = (I-\Pi)M^{k}|_{Y} \quad (k \ge 1).
\]
Now observe that to prove \eqref{eq:modulated-operator-contractivity}, 
it suffices to find an integer $p \ge 1$ such that
\begin{equation}\label{eq:monodromy-power-condition}
\|\widetilde{M}^{p}\|<1.
\end{equation}
(Recall that the unspecified norm is the one induced from 
$\|\cdot \|_{\max}$.) Indeed, for every integer $k \ge 1$ and 
$y \in Y$, we have \[
\|\widetilde{M}^{kp}y\|_{\ell^{\infty}} 
\le \sqrt{2} \|\widetilde{M}^{kp}y\|_{\max}
\le \sqrt{2} \|\widetilde{M}^{p}\|_{\max \to \max}^{k}\|y\|_{\max}
\le \sqrt{2} \|\widetilde{M}^{p}\|_{\max \to \max}^{k}
\|y\|_{\ell^{\infty}}.
\]
Thus for sufficiently large $k$, $\widetilde{M}^{kp}$ is contractive 
in the restriction of the standard norm on $\ell^{\infty}(\ZZ; \CC)$ 
to $Y$. 

It therefore remains to check that $\pi\cV$ is invertible and that 
\eqref{eq:monodromy-power-condition} holds for some integer $p \ge 1$.

\subsection{Reduction to finite computations}
\label{subsec:monodromy-estimates}

We now reduce the conditions at the end of Subsection
\ref{subsec:monodromy-criterion} to finite computations. For the 
contraction estimate $\|\widetilde{M}^{p}\| < 1$, we compare
$\widetilde{M}$ with an approximation 
$\overline{\widetilde{M}}: Y \to Y$ and control the error in the operator 
norm. Note that $\|\widetilde{M}^{p}\|$ can be bounded in terms of
$\|\overline{\widetilde{M}}^{j}\|$ $(1 \le j \le p)$ and
$\|\widetilde{M}-\overline{\widetilde{M}}\|$. Indeed, observe that \[
\widetilde{M}^{p}-\overline{\widetilde{M}}^{p}
= \sum_{j=0}^{p-1}
\overline{\widetilde{M}}^{p-1-j}
(\widetilde{M}-\overline{\widetilde{M}})
(\widetilde{M}^{j} - \overline{\widetilde{M}}^{j} 
+ \overline{\widetilde{M}}^{j}).
\]
Thus, given any bound $\|\widetilde{M}-\overline{\widetilde{M}}\|
\le \epsilon_{\widetilde{M}}$, setting
$\epsilon_{\widetilde{M}, 0}:=0$ and defining recursively
\begin{equation}\label{eq:monodromy-power-error}
\epsilon_{\widetilde{M}, p} :=
\epsilon_{\widetilde{M}}\sum_{j=0}^{p-1} \|\overline{\widetilde{M}}^{p-1-j}\|
(\|\overline{\widetilde{M}}^{j}\| + \epsilon_{\widetilde{M}, j})
\quad (p \ge 1)
\end{equation}
gives $\|\widetilde{M}^{p}-\overline{\widetilde{M}}^{p}\|
\le \epsilon_{\widetilde{M}, p}$ by induction.

We therefore seek an approximation $\overline{\widetilde{M}}$ for which
$\|\overline{\widetilde{M}}^{j}\| (1 \le j \le p)$ and
$\|\widetilde{M}-\overline{\widetilde{M}}\|$ can all be bounded by finite
computations. The resulting finite conditions are collected in Remark 
\ref{rmk:monodromy-verification}.

\subsubsection{Construction of the approximate operator}

Recall that \[
\widetilde{M}=(I-\Pi)M|_{Y}, 
\quad M=R^{-1}D\Phi^{\delta^{(\ast)}}(Z^{(\ast)}(0)).
\]
To construct an approximation $\overline{\widetilde{M}}$ to $\widetilde{M}$, 
we approximate the projection $\Pi$ and the flow derivative 
$D\Phi^{\delta^{(\ast)}}(Z^{(\ast)}(0))$ separately.

Since $\Pi = \cV(\pi \cV)^{-1}\pi$, we approximate $\Pi$ by replacing the
modes in \eqref{eq:monodromy-modes} with their numerical approximations. 
For this, we use the numerical approximation
$(\overline{\delta}, \overline{a})$ to $(\delta^{(\ast)}, a^{(\ast)})$ 
from Subsection \ref{subsec:orbit-finite-verification}. Recall that \[
Z^{(\ast)}(\delta^{(\ast)}\theta) = \sum_{k \ge 0}
(a^{(\ast)})^{(k)}\theta^{k} \quad (0 \le \theta \le 1).
\]
Since $\overline{a}^{(k)}=0$ for $k>K$, the corresponding 
approximation is \[
\overline{Z}(\overline{\delta}\theta):= 
\sum_{k=0}^{K}\overline{a}^{(k)} \theta^{k}.
\]
We therefore define the approximate modes
$\overline{v}_{\mathrm{r}}^{u}$, 
$\overline{v}_{\mathrm{r}}^{c}$, 
$\overline{v}_{\mathrm{S}}^{u}$, 
$\overline{v}_{\mathrm{G}, 1}^{c}$, and
$\overline{v}_{\mathrm{G}, 2}^{c}$
by replacing $Z^{(\ast)}(0)$ with $\overline{Z}(0)$ in
\eqref{eq:monodromy-modes}.

Let $\overline{\cV}:\RR^{m} \to \ell^{\infty}$ be defined in the 
same way as in \eqref{eq:monodromy-mode-map}, where the modes are replaced by the 
corresponding approximate modes.
Assuming that $\pi\overline{\cV}: \RR^{m} \to \RR^{m}$ is
invertible, as will be verified later, we take \[
\overline{\Pi} := \overline{\cV}(\pi\overline{\cV})^{-1}\pi
\]
as our approximation to $\Pi$. This $\overline{\Pi}$ is a bounded
projection onto $\mathrm{ran}\overline{\cV}$, with
$\ker\overline{\Pi} = \ker\pi = Y$.

We next approximate the flow derivative. Note that 
$s \mapsto D\Phi^{s}(Z^{(\ast)}(0))$ solves the variational equation
\eqref{eq:variational-ODE}. Setting $s=\delta^{(\ast)}\theta$ gives
\begin{equation} \label{eq:variational-ODE-rescaled}
\frac{\dd}{\dd\theta} D\Phi^{\delta^{(\ast)}\theta}(Z^{(\ast)}(0))=
\delta^{(\ast)}\left ( -I+DN(Z^{(\ast)}(\delta^{(\ast)}\theta)) 
\right ) D\Phi^{\delta^{(\ast)}\theta}(Z^{(\ast)}(0)), \quad
D\Phi^{0}(Z^{(\ast)}(0))=I.
\end{equation}
We approximate this equation by replacing $\delta^{(\ast)}$ with 
$\overline{\delta}$ and $Z^{(\ast)}(\delta^{(\ast)}\theta)$ with
$\overline{Z}(\overline{\delta}\theta)$. Let $\overline{\Psi}$ denote 
the solution of the resulting equation \[
\frac{\dd}{\dd\theta}\overline{\Psi}(\theta) = \overline{\delta}
\left ( -I+DN(\overline{Z}(\overline{\delta}\theta)) \right )
\overline{\Psi}(\theta), \quad \overline{\Psi}(0)=I.
\]
Setting $\overline{\Xi}(\theta)=e^{\overline{\delta}\theta}
\overline{\Psi}(\theta)$ gives
\begin{equation}
\label{eq:approximate-variational-equation}
\frac{\dd}{\dd\theta}\overline{\Xi}(\theta)=\overline{\delta}
DN(\overline{Z}(\overline{\delta}\theta)) \overline{\Xi}(\theta),
\quad \overline{\Xi}(0)=I.
\end{equation}

Now we extend this equation to complex time, allowing $\theta$ to 
be complex. Here, recall that the sequences and operators appearing 
in \eqref{eq:approximate-variational-equation} are written in the real 
coordinates of 
$\ell^{\infty}(\ZZ \times \{1, 2\}; \RR)$. In particular, complex conjugation 
in the definition of $N$ is already represented by a real-linear map on these 
coordinates. We extend the real coordinate matrices
$DN(\overline{a}^{(k)})$ complex-linearly to
$\ell^{\infty}(\ZZ \times \{1, 2\}; \CC)$, equipped with its supremum 
norm, and allow the time variable to be complex.
Using this complexified setup, we claim that the solution 
$\overline{\Xi}$ to \eqref{eq:approximate-variational-equation} can 
be extended to an operator-valued holomorphic function on $\CC$.

Consider the Picard iterates \[
\overline{\Xi}_{0}(\zeta):= I, \quad
\overline{\Xi}_{n+1}(\zeta) := I + \int_{0}^{\zeta} 
\overline{\delta}\left ( \sum_{k=0}^{K}DN(\overline{a}^{(k)})(\zeta')^{k}\right )
\overline{\Xi}_{n}(\zeta')\dd\zeta' \quad (n \ge 0),
\]
where the integral is along the line segment from $0$ to $\zeta$.
Since the coefficient in the integral equation is a polynomial in 
$\zeta'$, each iterate is also a polynomial. For any $r_{C}>0$, set \[
C_{r_{C}} := \sup_{|\zeta| \le r_{C}}
\left \| \overline{\delta} \sum_{k=0}^{K} 
DN(\overline{a}^{(k)})\zeta^{k} \right \| < \infty.
\]
Then induction on $n$ gives \[
\|\overline{\Xi}_{n+1}(\zeta) - \overline{\Xi}_{n}(\zeta)\| \le 
\frac{(C_{r_{C}}|\zeta|)^{n+1}}{(n+1)!} \quad 
(|\zeta| \le r_{C}, n \ge 0).
\]
Thus the iterates converge uniformly in the operator norm on every
closed disk. Therefore the limit
$\lim_{n \to \infty}\overline{\Xi}_{n}$ is an entire function.
Passing to the limit in the integral equation and using uniqueness
shows that this limit extends the solution of
\eqref{eq:approximate-variational-equation}.

Returning to the linearized bilateral relative return map 
$M=R^{-1}D\Phi^{\delta^{(\ast)}}(Z^{(\ast)}(0))$, recall that 
$\overline{\Psi}(1)$ approximates 
$D\Phi^{\delta^{(\ast)}}(Z^{(\ast)}(0))$, and
$\overline{\Psi}(1)=e^{-\overline{\delta}}\overline{\Xi}(1)$. Hence 
the corresponding approximation to $M$ is \[
R^{-1}\overline{\Psi}(1)
=e^{-\overline{\delta}}R^{-1}\overline{\Xi}(1).
\]
Since $\overline{\Xi}$ is entire, it can be expressed as a Taylor 
expansion \[
\overline{\Xi}(\theta) = \sum_{j \ge 0}\overline{\Xi}^{(j)}\theta^{j}.
\]
Replacing $\overline{\Xi}(1)$ by a finite sum, we define
our approximation to $M$ by
\begin{equation}\label{eq:approximate-monodromy}
\overline{M}:=e^{-\overline{\delta}}R^{-1}\sum_{j=0}^{J}
\overline{\Xi}^{(j)},
\end{equation}
where $J\ge 0$ is an integer to be chosen later. 

Combining $\overline{M}$ with the projection $\overline{\Pi}$, we 
obtain the desired approximation to $\widetilde{M}$: \[
\overline{\widetilde{M}}:=(I-\overline{\Pi})\overline{M}|_{Y}.
\]

\subsubsection{Norm bounds for the approximate operator}
\label{subsubsec:monodromy-norm-bounds}

We now bound the norms of $\overline{\widetilde{M}}$ and its powers 
using finite computations. We also obtain bounds for 
$\|\overline{M}\|$ and $\|I-\overline{\Pi}\|$, which will be used later 
to estimate $\|\widetilde{M}-\overline{\widetilde{M}}\|$.

We first show that these approximations are determined by finite
matrices and explicit shifts. Choose integers $I_{-} \le I_{+}$ such
that $[I_{-}, I_{+}]$ contains the spatial support of $\overline{a}$
and the sampling indices $n_{1}, n_{2}, n_{3}$ (and $n_{4}$ for GOY) 
in the definition of $\pi$.
Let $e_{n, a}$ denote the unit vector at $(n, a)$ in
$\ell^{\infty}(\ZZ \times \{1, 2\}; \RR)$. In complex notation, \[
e_{n, 1} = \mathbf{1}_{\{n\}}, \quad e_{n, 2} = i\mathbf{1}_{\{n\}}.
\]
Let $I_{\pi} \subset \ZZ \times \{1, 2\}$ be the set of scalar coordinates
sampled by $\pi$. Thus \[
Y = \{z \in \ell^{\infty}: z_{n, a}=0 \text{ for } (n, a) \in I_{\pi}\}.
\]

For the projection
$\overline{\Pi}=\overline{\cV}(\pi\overline{\cV})^{-1}\pi$, the input 
support is contained in $I_{\pi}$, since $\pi$ samples only these coordinates.
Its output lies in the span of the approximate modes. Since 
$\overline{Z}(0)=\overline{a}^{(0)}$ is a real sequence supported in
$[I_{-}, I_{+}] \times \{1\}$, $N(\overline{Z}(0))$ is supported in
$[I_{-}-1, I_{+}+1] \times \{1\}$, and so are
$\overline{v}_{\mathrm{r}}^{u}$ and $\overline{v}_{\mathrm{r}}^{c}$.
The vectors $\overline{v}_{\mathrm{G}, 1}^{c}$,
$\overline{v}_{\mathrm{G}, 2}^{c}$, and
$\overline{v}_{\mathrm{S}}^{u}$ are supported in
$[I_{-}, I_{+}] \times \{2\}$. Thus $\overline{\Pi}$ has its output 
support in $[I_{-}-1, I_{+}+1] \times \{1, 2\}$ and is represented by 
a finite matrix.

We next consider $\overline{M}$. Comparing Taylor coefficients in
\eqref{eq:approximate-variational-equation} gives
\begin{equation}\label{eq:approximate-variational-coefficients}
\overline{\Xi}^{(0)}=I, \quad \overline{\Xi}^{(j+1)} 
= \frac{\overline{\delta}}{j+1}\sum_{k=0}^{\min\{j, K\}}
DN(\overline{a}^{(k)})\overline{\Xi}^{(j-k)} \quad (j \ge 0),
\end{equation}
where we recall that $\overline{a}^{(k)}=0$ for $k>K$. By the finite 
interaction range of $N$, each $DN(\overline{a}^{(k)})$ has both 
input and output support in $[I_{-}-2, I_{+}+2] \times \{1, 2\}$.
So the recurrence shows that $\overline{\Xi}^{(j)}$ has both
input and output support in
$[I_{-}-2, I_{+}+2] \times \{1, 2\}$ for every $j \ge 1$.
Combining this with \eqref{eq:approximate-monodromy}, we conclude
that $\overline{M}-e^{-\overline{\delta}}R^{-1}$ has input
support in $[I_{-}-2, I_{+}+2] \times \{1, 2\}$ and output support 
in $[I_{-}-3, I_{+}+1] \times \{1, 2\}$. Thus $\overline{M}$ is the 
sum of the explicit shift $e^{-\overline{\delta}}R^{-1}$ and a finite 
matrix.

Therefore, \[
\overline{\widetilde{M}}
=(I-\overline{\Pi})\overline{M}|_{Y}
=\left . \left ( e^{-\overline{\delta}}R^{-1}
-e^{-\overline{\delta}}\overline{\Pi} R^{-1}
+(I-\overline{\Pi})(\overline{M}-e^{-\overline{\delta}}R^{-1})
\right ) \right |_{Y}
\]
is the restriction to $Y$ of the sum of
$\overline{\widetilde{M}}_{\mathrm{shift}} := 
e^{-\overline{\delta}}R^{-1}$ and the finite matrix \[
\overline{\widetilde{M}}_{\mathrm{finite}} := 
-e^{-\overline{\delta}}\overline{\Pi} R^{-1}
+(I-\overline{\Pi})(\overline{M}-e^{-\overline{\delta}}R^{-1}).
\] 
This matrix has its input support in $[I_{-}-2, I_{+}+2] \times \{1, 2\}$ 
and output support in $[I_{-}-3, I_{+}+1] \times \{1, 2\}$. 
The entries of $\overline{\widetilde{M}}_{\mathrm{finite}}$ 
are given by explicit finite formulas.

We now use these finite representations to reduce the required
norm estimates to finite computations.

\begin{lemma}[Finite-row norm formulas]
\label{lem:monodromy-finite-rows}
For every integer $p \ge 1$, we have
\begin{equation}\label{eq:monodromy-finite-row-norms}
\|\overline{\widetilde{M}}^{p}\| = \max \left \{ e^{-p\overline{\delta}},
\max_{\substack{I_{-}-p-2 \le i \le I_{+}+1 \\ a \in \{1, 2\}}}
\sum_{\substack{I_{-}-2 \le j \le I_{+}+p+1 \\ b \in \{1, 2\}, 
(j, b) \not\in I_{\pi}}}
\left | (\overline{\widetilde{M}}^{p}e_{j, b})_{i, a} \right |
\right \}.
\end{equation}
Moreover,
\begin{align}
\|\overline{M}\| &= \max\left \{ e^{-\overline{\delta}}, 
\max_{\substack{I_{-}-3 \le i \le I_{+}+1 \\ a \in \{1, 2\}}}
\sum_{j=I_{-}-2}^{I_{+}+2} \sum_{b=1}^{2} 
\left | (\overline{M}e_{j, b})_{i, a} \right | \right \},
\label{eq:monodromy-approximation-norm} \\
\|I-\overline{\Pi}\| &= \max\left \{ 1, 
\max_{\substack{I_{-}-1 \le i \le I_{+}+1 \\ a \in \{1, 2\}}}
\sum_{(j, b) \in I_{\pi} \cup \{(i, a)\}}
\left | ((I-\overline{\Pi})e_{j, b})_{i, a} \right | \right \}.
\label{eq:monodromy-complement-norm}
\end{align}
In particular, these norms can be evaluated using finite computations.
\end{lemma}

\begin{proof}
Since $Y$ consists of sequences whose coordinates in $I_{\pi}$ 
vanish and every output coordinate of $\overline{\widetilde{M}}^{p}$ 
depends on finitely many input coordinates, we have \[
\|\overline{\widetilde{M}}^{p}\| =
\sup_{\substack{i \in \ZZ \\ a \in \{1, 2\}}}
\sum_{\substack{j \in \ZZ, b \in \{1, 2\} \\ (j, b) \not\in I_{\pi}}}
|(\overline{\widetilde{M}}^{p}e_{j, b})_{i, a}|.
\]
We therefore determine which input coordinates can affect each output
coordinate, and then take the maximum of the resulting row sums.

We first consider $\overline{\widetilde{M}}$. If
$i\not\in [I_{-}-3, I_{+}+1]$, then 
$\overline{\widetilde{M}}_{\mathrm{finite}}$ does not contribute 
to the output coordinate $(i, a)$, so \[
(\overline{\widetilde{M}}z)_{i, a} = e^{-\overline{\delta}}z_{i+1, a} 
\quad (z \in Y).
\]
Since all the sampling indices selected by $\pi$ lie in 
$[I_{-}, I_{+}]$, we have $(i+1, a) \not\in I_{\pi}$. Thus the 
row sum corresponding to the output coordinate $(i, a)$ is 
$e^{-\overline{\delta}}$.

For $i \in [I_{-}-3, I_{+}+1]$, the shift part 
$\overline{\widetilde{M}}_{\mathrm{shift}}$ uses only the input 
coordinate $(i+1, a)$, while the finite part 
$\overline{\widetilde{M}}_{\mathrm{finite}}$ uses only input coordinates 
in $[I_{-}-2, I_{+}+2] \times \{1, 2\}$. Since $i+1$ also lies in 
$[I_{-}-2, I_{+}+2]$, taking the sums over these input coordinates 
gives \eqref{eq:monodromy-finite-row-norms} for $p=1$. 

Similarly, for $\overline{\widetilde{M}}^{p}$,
if $i \not\in [I_{-}-p-2, I_{+}+1]$, then $i, i+1, \ldots, i+p-1$ 
all lie outside $[I_{-}-3, I_{+}+1]$. Therefore 
$\overline{\widetilde{M}}_{\mathrm{finite}}$ does not affect this 
output and hence \[
(\overline{\widetilde{M}}^{p}z)_{i, a}
= e^{-p\overline{\delta}}z_{i+p, a} \quad (z \in Y).
\]
Again, $(i+p, a)\not\in I_{\pi}$, so the corresponding row sum 
is $e^{-p\overline{\delta}}$.

For $i \in [I_{-}-p-2, I_{+}+1]$, we show by induction on $p$ that 
the output coordinate $(\overline{\widetilde{M}}^{p}z)_{i, a}$ depends 
only on input coordinates in $[I_{-}-2, I_{+}+p+1] \times \{1, 2\}$.
The case $p=1$ was proved above.

For $p \ge 2$ and $i \in [I_{-}-p-1, I_{+}+1]$, apply the induction 
hypothesis to $\overline{\widetilde{M}}z\in Y$. The output coordinate
$(\overline{\widetilde{M}}^{p}z)_{i, a}$ then depends only
on the coordinates $(\overline{\widetilde{M}}z)_{j, b}$ with
$(j, b) \in [I_{-}-2, I_{+}+p] \times \{1, 2\}$. Through the shift part, 
$\overline{\widetilde{M}}_{\mathrm{shift}}$, these depend on input 
coordinates of $z$ in $[I_{-}-1, I_{+}+p+1] \times \{1, 2\}$, and
through the finite part $\overline{\widetilde{M}}_{\mathrm{finite}}$, 
these depend on input coordinates in $[I_{-}-2, I_{+}+2] \times \{1, 2\}$.
Both sets are contained in $[I_{-}-2, I_{+}+p+1] \times \{1, 2\}$.
For the remaining output coordinate with $i=I_{-}-p-2$, the one-step 
shift formula above gives \[
(\overline{\widetilde{M}}^{p}z)_{i, a}
=e^{-\overline{\delta}} (\overline{\widetilde{M}}^{p-1}z)_{i+1, a}.
\]
Since $i+1=I_{-}-p-1 \in [I_{-}-p-1, I_{+}+1]$, the induction hypothesis 
gives the required input support bound for this coordinate as well.
Thus, for every output coordinate with $i \in [I_{-}-p-2, I_{+}+1]$, 
only the input coordinates in $[I_{-}-2, I_{+}+p+1] \times \{1, 2\}$ 
can contribute. This proves 
\eqref{eq:monodromy-finite-row-norms}.

For $\overline{M}$, if $i\not\in [I_{-}-3, I_{+}+1]$, the finite part
$\overline{M}-e^{-\overline{\delta}}R^{-1}$ does not contribute to 
the output coordinate $(i, a)$, so \[
(\overline{M}z)_{i, a}
=e^{-\overline{\delta}}z_{i+1, a}.
\]
Thus the corresponding row sum is $e^{-\overline{\delta}}$.
For $i \in [I_{-}-3, I_{+}+1]$, the shift part 
$e^{-\overline{\delta}}R^{-1}$ uses only the input coordinate 
$(i+1, a)$, while the finite part uses only input coordinates in
$[I_{-}-2, I_{+}+2] \times \{1, 2\}$. Since $i+1$ also lies in 
$[I_{-}-2, I_{+}+2]$, taking the corresponding row sum and 
taking the maximum over all $(i, a)$ gives 
\eqref{eq:monodromy-approximation-norm}.

For $I-\overline{\Pi}$, if $i \not\in [I_{-}-1, I_{+}+1]$, the 
projection $\overline{\Pi}$ does not contribute to the output
coordinate $(i, a)$, so \[
((I-\overline{\Pi})z)_{i, a}=z_{i, a}.
\]
The sum of the absolute values of the input coefficients is therefore 
one. For $i \in [I_{-}-1, I_{+}+1]$, the identity uses only the input 
coordinate $(i, a)$, while $\overline{\Pi}$ uses only the input 
coordinates in $I_{\pi}$. Thus this output coordinate depends only on 
input coordinates in $I_{\pi} \cup \{(i, a)\}$. Taking the sums over 
these input coordinates and then the maximum over the output 
coordinates gives \eqref{eq:monodromy-complement-norm}.
\end{proof}

\subsubsection{Approximation error bounds}
\label{subsubsec:monodromy-errors}

We next estimate $\|\widetilde{M}-\overline{\widetilde{M}}\|$. Recall that
$\overline{\widetilde{M}}=(I-\overline{\Pi})\overline{M}|_{Y}$, where
$\overline{\Pi}$ approximates $\Pi$ and $\overline{M}$ is obtained by
approximating the variational equation and truncating the Taylor
series of $\overline{\Xi}$. We bound these errors separately and then
combine them.

\begin{lemma}[Error in the Sabra phase mode]
\label{lem:monodromy-sabra-phase-error}
Let $C_{Q}=\frac{11}{4}$ and $r=10^{-5}$. Assume that there exists 
an integer $j_{0}$ with $3 \le j_{0} \le I_{+}$ such that 
\begin{equation} \label{eq:sabra-tail-condition}
C_{Q} \left ( \frac{1+\sqrt{5}}{2}\right )
\left ( r + \max_{j_{0}-2 \le n \le I_{+}} \sum_{k=0}^{K} 
|\overline{a}_{n}^{(k)}| \right ) \le 1.
\end{equation}
Then,
\begin{equation} \label{eq:sabra-phase-mode-error} 
\begin{aligned}
& \|v_{\mathrm{S}}^{u}-\overline{v}_{\mathrm{S}}^{u}\|_{\ell^{\infty}} 
\le \epsilon_{\mathrm{S}} \\ 
& := \max \Bigg \{ \left ( \frac{1+\sqrt{5}}{2}\right )^{j_{0}-1}r, \\
& \qquad \qquad
C_{Q}\left ( \frac{1+\sqrt{5}}{2}\right )^{j_{0}+1}
\left ( r + \max_{j_{0}-2 \le n \le I_{+}}\sum_{k=0}^{K}
|\overline{a}_{n}^{(k)}| \right )^{2} + 
\max_{j_{0} \le n \le I_{+}} \left ( \frac{1+\sqrt{5}}{2}\right )^{n}
|\overline{a}_{n}^{(0)}| \Bigg \}.
\end{aligned}
\end{equation}
In particular, $\epsilon_{\mathrm{S}}$ can be computed using finite computations.
\end{lemma}

\begin{proof}
Let $R_{j}$ be as in Subsection \ref{subsec:orbit-spatial-decay}, \[
R_{j} = \sup_{0 \le \tau \le \delta^{(\ast)}}
\sup_{n \ge j} |(Z^{(\ast)}(\tau))_{n}|.
\]
By the construction in Subsection
\ref{subsec:orbit-finite-verification}, we have
\[
\sup_{0 \le \theta \le 1} \|Z^{(\ast)}(\delta^{(\ast)}\theta)
-\overline{Z}(\overline{\delta}\theta)\|_{\ell^{\infty}} \le
\|a^{(\ast)}-\overline{a}\|_{\ell_{\rho}^{1}} \le r.
\]
Since $\overline{a}$ has spatial support in $[I_{-}, I_{+}]$, this 
gives
\begin{equation} \label{eq:sabra-tail-envelope}
R_{j_{0}-2} \le 
\sup_{0 \le \theta \le 1} \|Z^{(\ast)}(\delta^{(\ast)}\theta)
-\overline{Z}(\overline{\delta}\theta)\|_{\ell^{\infty}}
+ \sup_{0 \le \theta \le 1}\sup_{n \ge j_{0}-2} 
|(\overline{Z}(\overline{\delta}\theta))_{n}| 
\le r + \max_{j_{0}-2 \le n \le I_{+}} 
\sum_{k=0}^{K}|\overline{a}_{n}^{(k)}|.
\end{equation}
For every $n \ge j_{0}$, applying \eqref{eq:orbit-right-tail-recurrence}
with $j=n+2$ and $j=j_{0}$, and using monotonicity of $R_{j}$, 
we obtain
\begin{align*}
\left ( \frac{1+\sqrt{5}}{2} \right )^{n+2} R_{n+2} 
&\le C_{Q} \left ( \frac{1+\sqrt{5}}{2}\right )^{n+2} R_{n}^{2} \\
&\le C_{Q} \left ( \frac{1+\sqrt{5}}{2}\right )^{2} R_{j_{0}}
\left ( \frac{1+\sqrt{5}}{2}\right )^{n}R_{n} \\
&\le C_{Q}^{2} \left ( \frac{1+\sqrt{5}}{2}\right )^{2} 
R_{j_{0}-2}^{2} \left ( \frac{1+\sqrt{5}}{2}\right )^{n}R_{n} \\
&\le \left ( \frac{1+\sqrt{5}}{2}\right )^{n}R_{n},
\end{align*}
where the last inequality uses \eqref{eq:sabra-tail-envelope} and the assumption 
\eqref{eq:sabra-tail-condition}. Iterating this inequality gives \[
\left ( \frac{1+\sqrt{5}}{2} \right )^{n} R_{n} 
\le \max \left \{ \left ( \frac{1+\sqrt{5}}{2}\right )^{j_{0}}R_{j_{0}},
\left ( \frac{1+\sqrt{5}}{2}\right )^{j_{0}+1}R_{j_{0}+1} \right \}
\le \left ( \frac{1+\sqrt{5}}{2}\right )^{j_{0}+1}R_{j_{0}}
\]
for any $n \ge j_{0}$. Consequently, 
\begin{align*}
\|v_{\mathrm{S}}^{u} - \overline{v}_{\mathrm{S}}^{u}\|_{\ell^{\infty}} 
&= \max \left \{ \sup_{n \le j_{0}-1} 
\left ( \frac{1+\sqrt{5}}{2}\right )^{n}|(Z^{(\ast)}(0))_{n}
-\overline{a}_{n}^{(0)}|,  
\sup_{n \ge j_{0}} \left ( \frac{1+\sqrt{5}}{2}\right )^{n}|(Z^{(\ast)}(0))_{n}
-\overline{a}_{n}^{(0)}| \right \}\\
& \le \max \left \{ \left (\frac{1+\sqrt{5}}{2}\right )^{j_{0}-1}r, 
\left ( \frac{1+\sqrt{5}}{2}\right )^{j_{0}+1}R_{j_{0}} 
+ \max_{j_{0} \le n \le I_{+}} \left ( \frac{1+\sqrt{5}}{2} \right )^{n} 
|\overline{a}_{n}^{(0)}| \right \}.
\end{align*}
Using $R_{j_{0}} \le C_{Q}R_{j_{0}-2}^{2}$ from
\eqref{eq:orbit-right-tail-recurrence} and \eqref{eq:sabra-tail-envelope} proves 
\eqref{eq:sabra-phase-mode-error}.
\end{proof}

\begin{lemma}[Mode and projection errors]
\label{lem:monodromy-projection-error}
Let $C_{Q}=\frac{11}{4}$ and $r=10^{-5}$. In the Sabra case, assume 
\eqref{eq:sabra-tail-condition} and let $\epsilon_{\mathrm{S}}$ be as 
in \eqref{eq:sabra-phase-mode-error}. Then we have
\begin{equation}\label{eq:monodromy-mode-error}
\|\cV-\overline{\cV}\| \le \epsilon_{V} :=
\left [ 1 + 2C_{Q}(2\|\overline{Z}(0)\|_{\ell^{\infty}}+r) \right ]r
+ \begin{cases}
2r & \text{for GOY}, \\
\epsilon_{\mathrm{S}} & \text{for Sabra}.
\end{cases}
\end{equation}
If, in addition, $\pi\overline{\cV}$ is invertible and
\begin{equation}\label{eq:monodromy-Neumann-condition}
\eta := \|(\pi\overline{\cV})^{-1}\|\epsilon_{V} < 1,
\end{equation}
then $\pi\cV$ is invertible and
\begin{equation}\label{eq:monodromy-projection-error}
\|\Pi-\overline{\Pi}\| \le \epsilon_{\Pi} := 
\frac{\eta}{1-\eta}\|I-\overline{\Pi}\|.
\end{equation}
Here $\|I-\overline{\Pi}\|$ is given by
\eqref{eq:monodromy-complement-norm}.
\end{lemma}

\begin{proof}
Since we are using the maximum norm, we have 
\begin{equation} \label{eq:mode-map-column-bound}
\|\cV-\overline{\cV}\| \le \|v_{\mathrm{r}}^{u}-\overline{v}_{\mathrm{r}}^{u}\|
+ \| v_{\mathrm{r}}^{c} - \overline{v}_{\mathrm{r}}^{c}\|      
+\begin{cases}
\|v_{\mathrm{G}, 1}^{c}-\overline{v}_{\mathrm{G}, 1}^{c}\| + 
\|v_{\mathrm{G}, 2}^{c}-\overline{v}_{\mathrm{G}, 2}^{c}\| & 
\text{for GOY}, \\
\|v_{\mathrm{S}}^{u} - \overline{v}_{\mathrm{S}}^{u}\| & 
\text{for Sabra}.
\end{cases}
\end{equation}
We estimate the error coming from each mode. By the construction 
in Subsection \ref{subsec:orbit-finite-verification}, we have \[
\|Z^{(\ast)}(0)-\overline{Z}(0)\|_{\ell^{\infty}} 
\le \|a^{(\ast)}-\overline{a}\|_{\ell_{\rho}^{1}} \le r.
\]
The quadratic formulas for $N$ give
\begin{align}
\|v_{\mathrm{r}}^{u}-\overline{v}_{\mathrm{r}}^{u}\|_{\ell^{\infty}}
&= \|N(Z^{(\ast)}(0))-N(\overline{Z}(0))\|_{\ell^{\infty}} 
\nonumber \\
&\le C_{Q} ( \|Z^{(\ast)}(0)\|_{\ell^{\infty}} 
+ \|\overline{Z}(0)\|_{\ell^{\infty}} )
\|Z^{(\ast)}(0)-\overline{Z}(0)\|_{\ell^{\infty}} \nonumber \\
&\le C_{Q}(2\|\overline{Z}(0)\|_{\ell^{\infty}}+r)r. \label{eq:real-unstable-mode-error}
\end{align}
Consequently,
\begin{align}
\|v_{\mathrm{r}}^{c}-\overline{v}_{\mathrm{r}}^{c}\|_{\ell^{\infty}}
&= \|-Z^{(\ast)}(0)+N(Z^{(\ast)}(0)) + 
\overline{Z}(0)-N(\overline{Z}(0))\|_{\ell^{\infty}} \nonumber \\
&\le \|v_{\mathrm{r}}^{u}-\overline{v}_{\mathrm{r}}^{u}\|_{\ell^{\infty}}
+ \|Z^{(\ast)}(0)-\overline{Z}(0)\|_{\ell^{\infty}} \nonumber \\
&\le \left [ 1+C_{Q}(2\|\overline{Z}(0)\|_{\ell^{\infty}}+r)\right ]r.
\label{eq:real-center-mode-error}
\end{align}
From the definitions of $v_{\mathrm{G}, 1}^{c}$ and 
$v_{\mathrm{G}, 2}^{c}$ given in \eqref{eq:monodromy-modes}, it is 
clear that 
\begin{equation} \label{eq:goy-phase-mode-error}
\|v_{\mathrm{G}, j}^{c}-\overline{v}_{\mathrm{G}, j}^{c}\|
\le \|Z^{(\ast)}(0)-\overline{Z}(0)\| \le r \quad 
(j=1, 2).
\end{equation}
Also, by Lemma \ref{lem:monodromy-sabra-phase-error}, we have
\begin{equation} \label{eq:sabra-mode-error-used}
\|v_{\mathrm{S}}^{u}-\overline{v}_{\mathrm{S}}^{u}\| \le
\epsilon_{\mathrm{S}}.
\end{equation}
Plugging \eqref{eq:real-unstable-mode-error}, \eqref{eq:real-center-mode-error}, \eqref{eq:goy-phase-mode-error}, and \eqref{eq:sabra-mode-error-used}
into \eqref{eq:mode-map-column-bound} gives \eqref{eq:monodromy-mode-error}.

Now suppose that $\pi\overline{\cV}$ is invertible and $\eta<1$.
The identity \[
\pi\cV = \pi \overline{\cV} \left [ I + (\pi\overline{\cV})^{-1} 
\pi (\cV-\overline{\cV}) \right ],
\]
together with $\|\pi\|=1$, shows that $\pi\cV$ is invertible by the 
Neumann series. Moreover, \[
\|(\pi\cV)^{-1}\| \le \frac{\|(\pi\overline{\cV})^{-1}\|}{1-\eta}.
\]
We are left to estimate $\|\Pi-\overline{\Pi}\|$.
From the definitions of $\Pi$ and $\overline{\Pi}$, observe that
\begin{equation*}
\begin{aligned}
\Pi-\overline{\Pi} &= \cV(\pi\cV)^{-1}\pi - \overline{\cV}
(\pi \overline{\cV})^{-1}\pi \\ 
&= \cV(\pi\cV)^{-1}\pi-\overline{\cV}(\pi\overline{\cV})^{-1}\pi\cV
(\pi\cV)^{-1}\pi \\
&= (I-\overline{\Pi})\cV(\pi\cV)^{-1}\pi \\
&= (I-\overline{\Pi})(\cV-\overline{\cV})(\pi\cV)^{-1}\pi,
\end{aligned}
\end{equation*}
where the last equality follows from 
$(I-\overline{\Pi})\overline{\cV}=0$. Taking norms gives \[
\|\Pi-\overline{\Pi}\| \le \|I-\overline{\Pi}\|\epsilon_{V}
\frac{\|(\pi\overline{\cV})^{-1}\|}{1-\eta}=\epsilon_{\Pi}.
\]
This completes the proof of the lemma.
\end{proof}

Next we bound the error coming from truncating the Taylor series of
$\overline{\Xi}$.

\begin{lemma}[Taylor remainder]
\label{lem:monodromy-Taylor-remainder}
For any $r_{C}>1$, we have
\begin{equation}\label{eq:monodromy-Taylor-error}
\left \| \overline{\Xi}(1)-\sum_{j=0}^{J}\overline{\Xi}^{(j)}
\right \| \le \epsilon_{T} := \frac{\displaystyle
\exp\left ( \overline{\delta}\sum_{k=0}^{K}\frac{r_{C}^{k+1}}{k+1}
\|DN(\overline{a}^{(k)})\|\right )}{r_{C}^{J+1}(1-r_{C}^{-1})}.
\end{equation}
Since each $DN(\overline{a}^{(k)})$ has finite input and output support,
$\epsilon_{T}$ can be evaluated using finite computations.
\end{lemma}

\begin{proof}
Recall that $\overline{\Xi}$ is entire as a function taking
values in the space of bounded operators on
$\ell^{\infty}(\ZZ \times \{1, 2\}; \CC)$. Cauchy's integral 
formula gives \[
\|\overline{\Xi}^{(j)}\| \le r_{C}^{-j}\sup_{|\zeta|=r_{C}}
\|\overline{\Xi}(\zeta)\|.
\]
Summing over $j>J$, we obtain \[
\left \| \overline{\Xi}(1)-\sum_{j=0}^{J}\overline{\Xi}^{(j)}
\right \| \le \frac{\displaystyle
\sup_{|\zeta|=r_{C}}\|\overline{\Xi}(\zeta)\|}{r_{C}^{J+1}(1-r_{C}^{-1})}.
\]
It remains to bound this supremum. For any $\phi \in \RR$, 
integrating \eqref{eq:approximate-variational-equation} along the 
segment $\zeta=t e^{i\phi}$, $0 \le t \le r_{C}$, and applying
Gr\"onwall's inequality gives
\begin{equation*}
\begin{aligned}
\|\overline{\Xi}(r_{C}e^{i\phi})\| &\le \exp\left (
\int_{0}^{r_{C}} \left \| \overline{\delta} \sum_{k=0}^{K}
DN(\overline{a}^{(k)})t^{k}e^{ik\phi}\right \|\dd t \right ) \\
&\le \exp\left ( \overline{\delta}\sum_{k=0}^{K}
\|DN(\overline{a}^{(k)})\|\int_{0}^{r_{C}}t^{k}\dd t \right ) \\
&=\exp\left ( \overline{\delta}\sum_{k=0}^{K}
\frac{r_{C}^{k+1}}{k+1}\|DN(\overline{a}^{(k)})\| \right ).
\end{aligned}
\end{equation*}
Taking the supremum over 
$\phi$ and substituting into the preceding remainder estimate gives 
the claim.
\end{proof}

To bound $\|M-\overline{M}\|$, we now combine this Taylor remainder
bound with an estimate for the error in the variational equation.

\begin{lemma}[Linearized return map error] \label{lem:monodromy-operator-error}
We have
\begin{equation}
\label{eq:monodromy-operator-error}
\|M-\overline{M}\| \le \epsilon_{M} := e^{-(\overline{\delta}-r)}
\left ( \epsilon_{\mathrm{var}} e^{\Lambda+\epsilon_{\mathrm{var}}}
+ re^{\Lambda} \right ) + e^{-\overline{\delta}}\epsilon_{T},
\end{equation}
where \[
\Lambda := \overline{\delta} \sum_{k=0}^{K}\|DN(\overline{a}^{(k)})\|,
\quad \epsilon_{\mathrm{var}} := \frac{r}{\overline{\delta}}\Lambda
+2\sqrt{2}C_{Q}(\overline{\delta}+r)r
\]
and $\epsilon_{T}$ is defined by \eqref{eq:monodromy-Taylor-error}.
\end{lemma}

\begin{proof}
Define \[
\Xi^{(\ast)}(\theta) := e^{\delta^{(\ast)}\theta}
D\Phi^{\delta^{(\ast)}\theta}(Z^{(\ast)}(0)).
\]
The variational equation \eqref{eq:variational-ODE-rescaled} gives \[
\frac{\dd}{\dd\theta} \Xi^{(\ast)}(\theta) = \delta^{(\ast)}
DN(Z^{(\ast)}(\delta^{(\ast)}\theta)) \Xi^{(\ast)}(\theta),
\quad \Xi^{(\ast)}(0)=I.
\]
We compare this equation with the approximate equation
\eqref{eq:approximate-variational-equation} on $0 \le \theta \le 1$.
Subtracting the two equations gives
\begin{equation*}
\begin{aligned}
\frac{\dd}{\dd\theta} \left ( \Xi^{(\ast)}(\theta) 
-\overline{\Xi}(\theta)\right ) 
&= \delta^{(\ast)}DN(Z^{(\ast)}(\delta^{(\ast)}\theta))
\left ( \Xi^{(\ast)}(\theta)-\overline{\Xi}(\theta)\right ) \\
& \quad + \left [ \delta^{(\ast)}DN(Z^{(\ast)}(\delta^{(\ast)}\theta))
-\overline{\delta}DN(\overline{Z}(\overline{\delta}\theta))\right ]
\overline{\Xi}(\theta),
\end{aligned}
\end{equation*}
with $\Xi^{(\ast)}(0)-\overline{\Xi}(0)=0$.
Applying Duhamel's formula to this equation, we obtain
\begin{equation} \label{eq:monodromy-variational-difference}
\|\Xi^{(\ast)}(1)-\overline{\Xi}(1)\| \le \int_{0}^{1} 
\|\Xi^{(\ast)}(1)\Xi^{(\ast)}(\theta)^{-1}\| 
\left \| \left [ \delta^{(\ast)}
DN(Z^{(\ast)}(\delta^{(\ast)}\theta))-\overline{\delta}
DN(\overline{Z}(\overline{\delta}\theta))\right ]\overline{\Xi}(\theta)    
\right \| \dd\theta.
\end{equation}

Since the coefficient in the approximate equation
\eqref{eq:approximate-variational-equation} satisfies \[
\sup_{0 \le \theta \le 1} \left \| \overline{\delta}
DN(\overline{Z}(\overline{\delta}\theta))\right \| \le 
\overline{\delta} \sum_{k=0}^{K}\|DN(\overline{a}^{(k)})\| = \Lambda,
\]
by Gr\"onwall's inequality we have
$\|\overline{\Xi}(\theta)\| \le e^{\Lambda\theta}$.

Also, since $N$ is quadratic, we have \[
\|DN(z)-DN(\widetilde{z})\| \le \sqrt{2} \|DN(z)-DN(\widetilde{z})\|_{
\ell^{\infty}(\ZZ; \CC) \to \ell^{\infty}(\ZZ; \CC)} 
\le 2\sqrt{2}C_{Q}\|z-\widetilde{z}\|_{\ell^{\infty}(\ZZ; \CC)},
\]
so it follows that
\begin{equation*}
\begin{aligned}
& \left \| \delta^{(\ast)}DN(Z^{(\ast)}(\delta^{(\ast)}\theta))
-\overline{\delta}DN(\overline{Z}(\overline{\delta}\theta)) \right \| \\
&\le |\delta^{(\ast)}-\overline{\delta}| 
\|DN(\overline{Z}(\overline{\delta}\theta))\| \\
& \quad + \delta^{(\ast)}\left \| DN(Z^{(\ast)}(\delta^{(\ast)}\theta))
- DN(\overline{Z}(\overline{\delta}\theta))\right \| \\
&\le r \frac{\Lambda}{\overline{\delta}} + 2\sqrt{2}C_{Q}(\overline{\delta}+r)
\|Z^{(\ast)}(\delta^{(\ast)}\theta)-\overline{Z}(\overline{\delta}\theta)    
\|_{\ell^{\infty}(\ZZ; \CC)} \le \epsilon_{\mathrm{var}},
\end{aligned}
\end{equation*}
where we have used $|\delta^{(\ast)}-\overline{\delta}| \le r$
and \[
\|Z^{(\ast)}(\delta^{(\ast)}\theta)-\overline{Z}(\overline{\delta}\theta)\|_{
\ell^{\infty}(\ZZ; \CC)}
\le \sum_{k \ge 0}\|(a^{(\ast)})^{(k)}-\overline{a}^{(k)}\|_{\ell^{1}} 
\le \|a^{(\ast)}-\overline{a}\|_{\ell_{\rho}^{1}} \le r,
\]
which follow from the construction in Subsection
\ref{subsec:orbit-finite-verification}.

Since \[
\|\delta^{(\ast)}DN(Z^{(\ast)}(\delta^{(\ast)}\theta))\|      
\le \|\overline{\delta}DN(\overline{Z}(\overline{\delta}\theta)) \|
+ \|\delta^{(\ast)}DN(Z^{(\ast)}(\delta^{(\ast)}\theta))
-\overline{\delta}DN(\overline{Z}(\overline{\delta}\theta)) \|
\le \Lambda+\epsilon_{\mathrm{var}},
\]
applying Gr\"onwall's inequality on $[\theta, 1]$ gives \[
\|\Xi^{(\ast)}(1)\Xi^{(\ast)}(\theta)^{-1}\| 
\le e^{(\Lambda+\epsilon_{\mathrm{var}})(1-\theta)}.
\] Plugging these estimates into \eqref{eq:monodromy-variational-difference} gives \[
\|\Xi^{(\ast)}(1)-\overline{\Xi}(1)\| \le \int_{0}^{1} 
e^{(\Lambda+\epsilon_{\mathrm{var}})(1-\theta)}\epsilon_{\mathrm{var}}
e^{\Lambda\theta}\dd\theta \le \epsilon_{\mathrm{var}}e^{\Lambda 
+\epsilon_{\mathrm{var}}}.
\]

Finally, the definitions of $M$ and $\overline{M}$ give
\begin{equation*}
\begin{aligned}
M-\overline{M} &= R^{-1}D\Phi^{\delta^{(\ast)}}(Z^{(\ast)}(0))
-e^{-\overline{\delta}}R^{-1}\sum_{j=0}^{J}\overline{\Xi}^{(j)} \\
&= R^{-1} \left [ e^{-\delta^{(\ast)}}\Xi^{(\ast)}(1)
- e^{-\overline{\delta}}\sum_{j=0}^{J}\overline{\Xi}^{(j)}\right ] \\
&= R^{-1}\left [ e^{-\delta^{(\ast)}} \left ( \Xi^{(\ast)}(1)
- \overline{\Xi}(1)\right ) + \left ( e^{-\delta^{(\ast)}}
-e^{-\overline{\delta}}\right )\overline{\Xi}(1) + 
e^{-\overline{\delta}} \left (\overline{\Xi}(1)
-\sum_{j=0}^{J}\overline{\Xi}^{(j)}\right )\right ].
\end{aligned}
\end{equation*}
Using $\|R^{-1}\|=1$, the inequalities \[
e^{-\delta^{(\ast)}} \le e^{-(\overline{\delta}-r)}, \quad
|e^{-\delta^{(\ast)}}-e^{-\overline{\delta}}| \le 
r e^{-(\overline{\delta}-r)},
\]
and Lemma \ref{lem:monodromy-Taylor-remainder}, we conclude that \[
\|M-\overline{M}\| \le e^{-(\overline{\delta}-r)} \left ( 
\epsilon_{\mathrm{var}}e^{\Lambda+\epsilon_{\mathrm{var}}}+re^{\Lambda}
\right ) + e^{-\overline{\delta}}\epsilon_{T}=\epsilon_{M}.
\]
\end{proof}

We now estimate $\|\widetilde{M}-\overline{\widetilde{M}}\|$ by combining the
preceding bounds. Assuming that $\pi\overline{\cV}$ is invertible
and $\eta<1$, we have
\begin{equation*}
\begin{aligned}
\widetilde{M}-\overline{\widetilde{M}} 
&= (I-\Pi)M|_{Y}-(I-\overline{\Pi})\overline{M}|_{Y}\\
&= (\overline{\Pi}-\Pi)M|_{Y}+(I-\overline{\Pi})(M-\overline{M})|_{Y}.
\end{aligned}
\end{equation*}
Lemmas \ref{lem:monodromy-projection-error} and
\ref{lem:monodromy-operator-error} therefore give
\begin{equation} \label{eq:monodromy-reduced-error}
\|\widetilde{M}-\overline{\widetilde{M}}\|_{Y \to Y} 
\le \epsilon_{\widetilde{M}} := \epsilon_{\Pi} ( \|\overline{M}\|
+ \epsilon_{M}) + \|I-\overline{\Pi}\|\epsilon_{M}.
\end{equation}
The two norms on the right-hand side are given by Lemma
\ref{lem:monodromy-finite-rows}.

\begin{remark}
\label{rmk:monodromy-verification}
The preceding estimates reduce the verification of the conditions
at the end of Subsection \ref{subsec:monodromy-criterion} to finite 
computations. It remains to specify $I_{-}, I_{+}$, the sampling 
indices $n_{1}, n_{2}, n_{3}$ (and $n_{4}$ for GOY) in the definition 
of $\pi$, and $J, r_{C}, p$ as above. For Sabra, we also specify 
$j_{0}$ used to compute $\epsilon_{\mathrm{S}}$. We then certify 
the following:
\begin{itemize}
\item $\pi\overline{\cV}$ is invertible.
\item In the Sabra case, \eqref{eq:sabra-tail-condition} holds.
\item $\eta<1$, where $\eta$ is defined in
\eqref{eq:monodromy-Neumann-condition}.
\item $\|\overline{\widetilde{M}}^{p}\| + \epsilon_{\widetilde{M}, p}<1$, 
where $\|\overline{\widetilde{M}}^{j}\|$ is evaluated using
\eqref{eq:monodromy-finite-row-norms}, $\epsilon_{\widetilde{M}}$ is given by
\eqref{eq:monodromy-reduced-error}, and $\epsilon_{\widetilde{M}, p}$ is 
obtained from \eqref{eq:monodromy-power-error}.
\end{itemize}
Indeed, the first three conditions imply that $\pi\cV$ is invertible 
by Lemma \ref{lem:monodromy-projection-error}, and the fourth 
condition, together with \eqref{eq:monodromy-power-error} and
\eqref{eq:monodromy-reduced-error}, gives $\|\widetilde{M}^{p}\|<1$.
\end{remark}

\subsection{Computer-assisted verification and completion of the proof}
\label{subsec:monodromy-verification}

In this subsection, we describe the computer-assisted verification
needed to complete the proof of Theorem \ref{thm:monodromy-splitting}.
All code and input data used here are available in the supplementary
material and in the GitHub repository \cite{CAParchive}.

We use the same numerical data and parameters as in Subsection
\ref{subsec:orbit-finite-verification}. In particular, we use the
same $(\overline{\delta}, \overline{a})$, and \[
K=19, \quad r=10^{-5}.
\]
We now specify the additional parameters listed in Remark
\ref{rmk:monodromy-verification}. Since the spatial support of
$\overline{a}$ is contained in $I=[-51, 21]$, we take $I_{-}=-51$
and $I_{+}=21$. We choose $(n_{1}, n_{2}, n_{3}, n_{4})=(-2, 2, 0, 1)$ 
for GOY and $(n_{1}, n_{2}, n_{3})=(0, 2, 2)$ for Sabra in the 
definition of $\pi$, set the Taylor truncation order in
\eqref{eq:approximate-monodromy} to $J=80$, and use $r_{C}=2$
in Lemma \ref{lem:monodromy-Taylor-remainder}.
We take $p=4$ for GOY, and $p=5$ for Sabra. In the Sabra case, we 
also take $j_{0}=9$.

We now verify the conditions listed in Remark
\ref{rmk:monodromy-verification}. The resulting bounds are recorded 
in the following lemma.

\begin{lemma}[Verified bounds]
\label{lem:monodromy-finite-verification}
For the data and parameters specified above, we have 
\begin{equation}\label{eq:monodromy-verified-bounds}
\|\overline{\widetilde{M}}^{p}\|+\epsilon_{\widetilde{M}, p} \le
\begin{cases}
0.824 & \text{for GOY}, \\
0.990 & \text{for Sabra}.
\end{cases}
\end{equation}
Moreover, $\pi\overline{\cV}$ is invertible and $\eta<1$.
In the Sabra case, \eqref{eq:sabra-tail-condition} also holds.
\end{lemma}

\begin{proof}
The script \texttt{verify\_contraction.py} in \texttt{codes/section\_3} 
carries out the following computations. As in Subsection 
\ref{subsec:orbit-finite-verification}, interval computations are 
performed with \texttt{python-flint}. We use a working precision of 
$64$ bits.

The function \texttt{check\_parameters} checks that the sampling 
indices $n_{1}, n_{2}, n_{3}$ (and $n_{4}$ for GOY) belong to 
$[I_{-}, I_{+}]$, that $n_{1} \neq n_{2}$ in both models, and 
$n_{3} \neq n_{4}$ for GOY. It also checks that 
$3 \le j_{0} \le I_{+}$ for the Sabra case.

For the Sabra case, the function \texttt{compute\_phase\_error} verifies 
\eqref{eq:sabra-tail-condition} using interval arithmetic. It then evaluates the
right-hand side of \eqref{eq:sabra-phase-mode-error} and returns an interval enclosure of 
$\epsilon_{\mathrm{S}}$.

The function \texttt{check\_invertibility} evaluates the approximate 
modes $\overline{v}_{\mathrm{r}}^{u}$, $\overline{v}_{\mathrm{r}}^{c}$, 
$\overline{v}_{\mathrm{G}, 1}^{c}$, $\overline{v}_{\mathrm{G}, 2}^{c}$,
$\overline{v}_{\mathrm{S}}^{u}$ and returns interval enclosures 
of $\overline{\cV}$ and $\pi\overline{\cV}$. It also checks that the 
determinant of $\pi\overline{\cV}$ is nonzero for each model.

The function \texttt{compute\_eta} evaluates $\eta$ using 
\eqref{eq:monodromy-Neumann-condition}. It uses the enclosure of 
$\pi\overline{\cV}$ returned by \texttt{check\_invertibility} and 
computes $\epsilon_{V}$ using \eqref{eq:monodromy-mode-error}. The 
function returns an interval enclosure of $\eta$ and checks that 
its upper endpoint is strictly less than one.

The function \texttt{compute\_contraction\_norm} computes a rigorous
upper bound for $\|\overline{\widetilde{M}}^{p}\|      
+\epsilon_{\widetilde{M}, p}$. It first constructs the finite matrices 
$DN(\overline{a}^{(k)})$ and computes interval enclosures of 
$\overline{\Xi}^{(j)}$ using 
\eqref{eq:approximate-variational-coefficients}. It then uses 
\eqref{eq:approximate-monodromy} to construct the finite submatrix 
of $\overline{M}$ needed in the norm estimates. The function also 
uses the returned $\overline{\cV}$ and computes
$(\pi\overline{\cV})^{-1}$ to form $\overline{\Pi}$, and uses 
$\overline{\widetilde{M}}=(I-\overline{\Pi})\overline{M}|_{Y}$ to 
construct the required finite submatrix of 
$\overline{\widetilde{M}}^{j} (1 \le j \le p)$. The norms
$\|\overline{\widetilde{M}}^{j}\|$ are computed using
\eqref{eq:monodromy-finite-row-norms}. It also bounds 
$\|I-\overline{\Pi}\|$ and $\|\overline{M}\|$ using 
\eqref{eq:monodromy-complement-norm} and 
\eqref{eq:monodromy-approximation-norm}. Next, using the enclosure 
of $\eta$ returned by \texttt{compute\_eta}, it evaluates 
\eqref{eq:monodromy-projection-error} to bound $\epsilon_{\Pi}$. 
It evaluates \eqref{eq:monodromy-Taylor-error} and
\eqref{eq:monodromy-operator-error} to bound $\epsilon_{T}$ and
$\epsilon_{M}$, and combines these bounds using
\eqref{eq:monodromy-reduced-error} to bound 
$\epsilon_{\widetilde{M}}$. All these computations use interval 
arithmetic. Finally, it evaluates the recursion 
\eqref{eq:monodromy-power-error}, using $\epsilon_{\widetilde{M}}$ 
as an upper bound for $\|\widetilde{M}-\overline{\widetilde{M}}\|$, 
to compute $\epsilon_{\widetilde{M}, p}$. 
The function returns the upper endpoint of the resulting enclosure
of $\|\overline{\widetilde{M}}^{p}\| + \epsilon_{\widetilde{M}, p}$ 
and checks that this bound is strictly less than $1$.

Table \ref{tab:monodromy-verification} summarizes the contraction 
bound and verification time for each model. The times 
were measured on an Intel Core Ultra 5 125H processor.

\begin{table}[htbp]
\centering
\begin{tabular}{lccc}
\toprule
Model & $p$ &
$\|\overline{\widetilde{M}}^{p}\|+\epsilon_{\widetilde{M},p}$
& Time (s) \\
\midrule
GOY   & $4$ & $0.8236$ & $\approx 45.482$ \\
Sabra & $5$ & $0.9898$ & $\approx 45.811$ \\
\bottomrule
\end{tabular}
\caption{Upper bounds for 
$\|\overline{\widetilde{M}}^{p}\|+\epsilon_{\widetilde{M},p}$ and 
computation times for the verification. The upper bounds are
rounded upward for display. Times are from one recorded run and are 
rounded to the nearest millisecond.}
\label{tab:monodromy-verification}
\end{table}
\end{proof}

\begin{proof}[Proof of Theorem \ref{thm:monodromy-splitting}]
Lemma \ref{lem:monodromy-finite-verification}, together with
Lemma \ref{lem:monodromy-projection-error}, shows that
$\pi\cV$ is invertible.
Consequently, in each model the vectors in
\eqref{eq:monodromy-modes} that generate $V$ are linearly
independent.
The action computed in Subsection \ref{subsec:monodromy-modes}
shows that $V$ is $M$-invariant.
The formula
$\Pi=\cV(\pi\cV)^{-1}\pi$ gives a bounded projection
onto $V$ with kernel $Y$, and its value depends only on the finitely
many coordinates selected by $\pi$.

By \eqref{eq:monodromy-power-error},
\eqref{eq:monodromy-reduced-error}, and Lemma
\ref{lem:monodromy-finite-verification}, we have \[
\|\widetilde{M}^{p}\|<1.
\]
The argument at the end of Subsection
\ref{subsec:monodromy-criterion} therefore gives
\eqref{eq:modulated-operator-contractivity}
for a sufficiently large multiple of this exponent $p$.
This completes the proof.
\end{proof}

\section{Construction and stability of blowup solutions}
\label{sec:boundary}

In this section, we construct blowup solutions to
\eqref{eq:main-renorm-common} using a modulation argument and
complete the proof of Theorem \ref{thm:main}. 

Throughout this section, $N$ denotes either $N_{\mathrm{GOY}}$
or $N_{\mathrm{Sabra}}$.\footnote{The same convention applies to
other model-dependent quantities. We use the subscripts
$\mathrm{GOY}$ and $\mathrm{Sabra}$ when the two models need
to be distinguished explicitly.} Unless a model is specified,
all statements apply to both choices of $N$. We regard
$\ell^{\infty}=\ell^{\infty}(\ZZ; \CC)$ as a Banach space over
$\RR$ and use its standard complex supremum norm. 
Finite-dimensional spaces $\RR^{d}$ carry their maximum norms. 
Unless otherwise specified, operator and multilinear norms are 
induced by these norms and their restrictions.
Let $Z^{(\ast)}$ be a 
bilateral travelling wave constructed in Section \ref{sec:bilateral} 
and let $M$, $\pi$, $\Pi$, $V$, $Y$, and $\widetilde{M}=(I-\Pi)M|_{Y}$ be defined 
as in Section \ref{sec:monodromy}.

Recall that constructing a solution to \eqref{eq:main-renorm-common}
that blows up at time $t=1$ with the desired rate $(1-t)^{-1}$ is
equivalent to constructing a solution $W$ to \eqref{eq:main-ss}
with \eqref{eq:condition-nontriviality}. Note that $Z^{(\ast)}$
solves \eqref{eq:main-ss-Z} and clearly satisfies the bilateral 
analogue of \eqref{eq:condition-nontriviality}. However, restricting 
$Z^{(\ast)}$ to the positive shells does not yield a solution of
\eqref{eq:main-ss}, since the nonpositive shells affect the dynamics
of the positive shells. We remove this influence by introducing
a boundary correction.

Define the boundary correction
$B: \ell^{\infty}(\ZZ; \CC) \to \ell^{\infty}(\ZZ; \CC)$ by \[
(B_{\mathrm{GOY}}(Z))_{1} = -2Z_{0}^{\ast}Z_{-1}^{\ast} 
-\frac{1}{2} Z_{2}^{\ast}Z_{0}^{\ast}, \quad
(B_{\mathrm{GOY}}(Z))_{2}
=-2Z_{1}^{\ast}Z_{0}^{\ast}, \quad
(B_{\mathrm{GOY}}(Z))_{n}
=0\quad(n\in\ZZ\setminus\{1, 2\})
\]
for the GOY model, and \[
(B_{\mathrm{Sabra}}(Z))_{1}
=-2Z_{0}Z_{-1}-\frac{1}{2} Z_{2}Z_{0}^{\ast}, \quad
(B_{\mathrm{Sabra}}(Z))_{2}
=-2Z_{1}Z_{0}, \quad
(B_{\mathrm{Sabra}}(Z))_{n}
=0\quad(n\in\ZZ\setminus\{1, 2\})
\]
for the Sabra model.

Then $B$ cancels the influence of nonpositive shells on the positive
shells. Consequently, the boundary-corrected equation
\begin{equation}\label{eq:boundary-modified}
\partial_{\tau}\widehat{Z}
=-\widehat{Z}+N(\widehat{Z})+B(\widehat{Z})
\end{equation}
gives a closed system for the components
$\{\widehat{Z}_{n}\}_{n\ge1}$, and their dynamics agree exactly
with those of \eqref{eq:main-ss}, the equation on the positive shells 
with zero boundary conditions. Therefore constructing a solution of
\eqref{eq:boundary-modified} and restricting it to the positive
shells gives a solution to \eqref{eq:main-ss}. Let $\Phi_{0}^{s}$
denote the local flow of \eqref{eq:boundary-modified}.

Recall that $Z^{(\ast)}$ solves \eqref{eq:main-ss-Z} exactly, so the
only residual in \eqref{eq:boundary-modified} comes from the boundary
correction $B(Z^{(\ast)})$. Note that the correction $B(Z)$ depends
only on the components $Z_{-1}, Z_{0}, Z_{1}, Z_{2}$. As $Z^{(\ast)}$
advances one shell to the right over each relative period
$\delta^{(\ast)}$, its left spatial decay from Theorem
\ref{thm:orbit-existence} makes these components tend to zero.
Hence the residual $B(Z^{(\ast)}(\tau))$ tends to zero as
$\tau\to\infty$. To make this residual small from the beginning,
we use the translated travelling wave $R^{J}Z^{(\ast)}$ instead of
$Z^{(\ast)}$, with a sufficiently large integer $J$. This wave
still satisfies \eqref{eq:main-ss-Z} exactly, while the boundary
correction $B(R^{J}Z^{(\ast)}(\tau))$ is small for $\tau\ge0$.
Thus $R^{J}Z^{(\ast)}$ provides an approximate solution of
\eqref{eq:boundary-modified}.

We seek an exact solution $\widehat{Z}$ of \eqref{eq:boundary-modified} 
near this approximate solution $R^{J}Z^{(\ast)}$. For $j \ge J$, the reference 
wave $R^{J}Z^{(\ast)}$ takes the value $R^{j}Z^{(\ast)}(0)$ at time 
$\tau = (j-J)\delta^{(\ast)}$. We express the values of $\widehat{Z}$ at these 
times as \[
\widehat{Z}((j-J)\delta^{(\ast)}) = R^{j}\widehat{Z}^{[j]} \quad (j \ge J).
\]
A first approach would be to find a solution whose recentered values 
$\widehat{Z}^{[j]}$ converge to $Z^{(\ast)}(0)$.

Since $R^{j}\widehat{Z}^{[j]}$ and $R^{j+1}\widehat{Z}^{[j+1]}$ are values
of the solution $\widehat{Z}$ to \eqref{eq:boundary-modified} at 
successive times, they must satisfy \[
R^{j+1}\widehat{Z}^{[j+1]} = \Phi_{0}^{\delta^{(\ast)}}(R^{j}\widehat{Z}^{[j]})
\quad (j \ge J).
\] Removing the spatial shifts gives the map from $\widehat{Z}^{[j]}$ to 
$\widehat{Z}^{[j+1]}$:
\begin{equation}\label{eq:boundary-return-maps}
F_{j}:=R^{-(j+1)}\Phi_{0}^{\delta^{(\ast)}}R^{j}
=R^{-1}\Phi_{j}^{\delta^{(\ast)}}\qquad(j\ge0),
\end{equation}
where $\Phi_{j}^{s}:=R^{-j}\Phi_{0}^{s}R^{j}$ is the local flow of
\begin{equation}\label{eq:moving-boundary-equation}
\partial_{\tau}\widehat{Z}
=-\widehat{Z}+N(\widehat{Z})+B_{j}(\widehat{Z}),
\qquad B_{j}(Z):=R^{-j}B(R^{j}Z).
\end{equation}
Thus, if successive applications of
$F_{J}, F_{J+1}, F_{J+2}, \ldots$ to $\widehat{Z}^{[J]}$
give a sequence converging to $Z^{(\ast)}(0)$, then the
corresponding solution of \eqref{eq:boundary-modified}
has initial value $R^{J}\widehat{Z}^{[J]}$ and approaches
$R^{J}Z^{(\ast)}$ at the sampling times $(j-J)\delta^{(\ast)}$.
As discussed earlier, the left decay of $Z^{(\ast)}$
makes the boundary correction small along the reference
wave for large $j$. If we omit this correction, all the
relative return maps $F_{j}$ become the bilateral relative
return map $F=R^{-1}\Phi^{\delta^{(\ast)}}$.
This is why we studied the linearized dynamics of $F$ in 
Section \ref{sec:monodromy}.
A rigorous comparison between $F_{j}$ and $F$ is given in 
Subsection \ref{subsec:boundary-return-maps}.

However, if we start with an arbitrary initial state $\widehat{Z}^{[J]}$ 
near $Z^{(\ast)}(0)$, then we cannot expect the sequence 
$(\widehat{Z}^{[j]})_{j\ge J}$ to converge to $Z^{(\ast)}(0)$, as the 
linearized bilateral relative return map $M=DF(Z^{(\ast)}(0))$ has 
nonstable modes. This behavior is expected, as the symmetries inherent 
in the equation \eqref{eq:boundary-modified} allow us to construct other 
approximate solutions from $R^{J}Z^{(\ast)}$, and therefore a solution 
starting from arbitrary nearby initial data need not converge to a 
specifically chosen reference wave, $R^{J}Z^{(\ast)}$. But as we saw 
in Theorem \ref{thm:monodromy-splitting} and Subsection 
\ref{subsec:monodromy-modes}, the nonstable modes of $M$ all come from 
symmetries. So, instead of simply applying 
$F_{J}, F_{J+1}, F_{J+2}, \ldots$ to $\widehat{Z}^{[J]}$, we insert 
symmetry transformations so that each modulated iterate lies in 
$Z^{(\ast)}(0) + Y$. This process will be called modulation. 
Since the symmetry transformations act in the directions spanning $V$ to 
first order, after omitting the boundary correction, modulation is 
formally approximated by applying $I-\Pi$ to the perturbation. 
By Theorem \ref{thm:monodromy-splitting}, we may choose an integer $p \ge 1$
such that $\widetilde{M}^{p} = ((I-\Pi)M|_{Y})^{p}=(I-\Pi)M^{p}|_{Y}$ is contractive. 
We fix such a $p$ and apply modulation after $p$ steps of iteration, i.e., 
after solving the equation for a time of $p\delta^{(\ast)}$ and shifting 
$p$ shells to the left. Then we may expect the resulting iteration to converge 
to $Z^{(\ast)}(0)$.

To be more precise, for all sufficiently large $j$ and for each $\widetilde{Z}$ near 
$Z^{(\ast)}(0)$, define the modulated return map by
\begin{equation} \label{eq:boundary-normalized-return}
\widetilde{F}_{j}(\widetilde{Z}) := e^{a_{j}(\widetilde{Z})}
T_{\theta_{j}(\widetilde{Z})}R^{-p}\Phi_{j}^{p\delta^{(\ast)}+b_{j}(\widetilde{Z})}
(\widetilde{Z}),
\end{equation}
where the parameters $a_{j}(\widetilde{Z}), \theta_{j}(\widetilde{Z}), 
b_{j}(\widetilde{Z})$ are chosen 
near zero so that
\begin{equation} \label{eq:boundary-normalization-condition}
\widetilde{F}_{j}(\widetilde{Z}) - Z^{(\ast)}(0) \in Y.
\end{equation} (Here, $\theta_{j}$ represents two real parameters for GOY, and 
one real parameter for Sabra.) These parameters are chosen in 
Subsection \ref{subsec:boundary-modulation}.
Now let $(\widetilde{Z}^{[J+kp]})_{k \ge 0}$ be 
a sequence obtained by successive applications of 
$\widetilde{F}_{J}, \widetilde{F}_{J+p}, \widetilde{F}_{J+2p}, \ldots$ to 
$\widetilde{Z}^{[J]}$. As $\widetilde{M}^{p}$ is contractive, we can expect this to 
converge to $Z^{(\ast)}(0)$. This is explained in detail 
in Subsection \ref{subsec:boundary-contraction}. 

In Subsection \ref{subsec:boundary-reconstruction}, we use the sequence 
$(\widetilde{Z}^{[J+kp]})_{k \ge 0}$ to recover a blowup solution with 
initial data $R^{J}\widetilde{Z}^{[J]}$. We also prove the desired blowup 
rate in $\ell^{\infty}$. Finally, we finish the proof of Theorem 
\ref{thm:main} in Subsection \ref{subsec:boundary-main-proof}.

\subsection{Comparison of the moving-boundary flows}
\label{subsec:boundary-return-maps}

As $j$ increases, the correction $B_{j} = R^{-j}BR^{j}$ involves only
components farther into the negative tail. These components are
exponentially small along the reference wave $Z^{(\ast)}$, whose left 
decay comes from Theorem \ref{thm:orbit-existence}.
Thus at $Z^{(\ast)}(0)$ we expect the moving-boundary flow map 
$\Phi_{j}^{s}$, the local flow of 
\eqref{eq:moving-boundary-equation}, to approach the bilateral flow map 
$\Phi^{s}$, the local flow of \eqref{eq:main-ss-Z}.
The following lemma quantifies this convergence.

\begin{lemma}[Comparison of the flows]
\label{lem:boundary-flow-comparison}
Let $p \ge 1$ be a fixed integer and let $\nu_{-}$ be a fixed real number such 
that $1<\nu_{-}<e^{\delta^{(\ast)}}$. 
Then there exist an open ball $U$ centered at 
$Z^{(\ast)}(0) \in \ell^{\infty}$, an integer $J_{0} \ge 2$, 
and a constant $C>0$ such that the following hold.
\begin{itemize}
\item[(i)] For every $0 \le s \le 2p\delta^{(\ast)}$ and $j \ge J_{0}$, the 
maps $\Phi^{s}$ and $\Phi_{j}^{s}$ are well-defined and $C^{2}$ on $U$. 
Moreover, \begin{equation} \label{eq:boundary-flow-derivative-bounds}
\|\partial_{s}^{\alpha}D_{Z}^{\beta}\Phi^{s}(Z)\| \le C, \quad
\|\partial_{s}^{\alpha}D_{Z}^{\beta}\Phi_{j}^{s}(Z)\| \le C \quad 
(Z \in U, \alpha + \beta \le 2).
\end{equation}
\item[(ii)] For every $0 \le s \le 2p\delta^{(\ast)}$, $j \ge J_{0}$, 
and $Z \in U$, we have 
\begin{equation}\label{eq:boundary-flow-near-profile}
\|\Phi_{j}^{s}(Z) - Z^{(\ast)}(s)\| \le C ( \|Z-Z^{(\ast)}(0)\| + \nu_{-}^{-2j}).
\end{equation}
\item[(iii)] For every $0 \le s \le 2p\delta^{(\ast)}$ and $j \ge J_{0}$, we have
\begin{equation}\label{eq:boundary-flow-difference}
\|D\Phi_{j}^{s}(Z^{(\ast)}(0)) - D\Phi^{s}(Z^{(\ast)}(0))\|
\le C\nu_{-}^{-j}.
\end{equation}
\end{itemize}
\end{lemma}

\begin{proof}
In the proof, $C$ will denote a constant not depending on
$j \ge J_{0}$ and $Z \in U$, whose value may change from line
to line.

Since $\Phi^{s}(Z^{(\ast)}(0)) = Z^{(\ast)}(s)$ is defined 
for all $0 \le s \le 2p\delta^{(\ast)}$, the standard Lipschitz-dependence 
estimates give an open ball $U$ centered at $Z^{(\ast)}(0)$ such that
\begin{equation*}
\|\Phi^{s}(Z) - \Phi^{s}(Z^{(\ast)}(0))\| \le C\|Z-Z^{(\ast)}(0)\|
\quad (0 \le s \le 2p\delta^{(\ast)}, Z \in U).
\end{equation*}
This in particular gives a uniform bound on $\Phi^{s}(Z)$.

To prove well-definedness of $\Phi_{j}^{s}$, fix $Z \in U$ and let \[
X(s) = Z^{(\ast)}(s), \quad X_{j}(s) = \Phi_{j}^{s}(Z).
\] Then we have
\begin{equation} \label{eq:boundary-flow-differential-inequality}
\begin{aligned}
\left \| \frac{\dd}{\dd s}(X-X_{j})\right \| &= 
\|-X + X_{j} + N(X) - N(X_{j}) - B_{j}(X_{j})\| \\             
&\le \|X-X_{j}\| + \|N(X)-N(X_{j})\| + \|B_{j}(X) - B_{j}(X_{j})\| + \|B_{j}(X)\|.
\end{aligned}
\end{equation}
As long as $X_{j}$ exists and $\|X-X_{j}\| \le 1$, both trajectories are 
uniformly bounded, and hence
\begin{equation}\label{eq:boundary-nonlinearity-lipschitz}
\|N(X)-N(X_{j})\| + \|B_{j}(X) - B_{j}(X_{j})\| \le C\|X-X_{j}\|.
\end{equation}
Also, the only nonzero components of $B_{j}=R^{-j}BR^{j}$ satisfy \[
|(B_{j}(Z))_{1-j}| \le 2|Z_{-j}Z_{-j-1}| + \frac{1}{2}|Z_{2-j}Z_{-j}|, 
\quad |(B_{j}(Z))_{2-j}| \le 2|Z_{1-j}Z_{-j}|.
\]
Together with the spatial decay in Theorem \ref{thm:orbit-existence} and 
relative periodicity, we have \[
|(X(s))_{n}| = |(Z^{(\ast)}(s))_{n}| \le C\nu_{-}^{n}
\] for all $0 \le s \le 2p\delta^{(\ast)}$ and $n \le 0$, and therefore
\begin{equation}\label{eq:boundary-correction-on-profile}
\|B_{j}(X(s))\| \le C\nu_{-}^{-2j}.
\end{equation}
Plugging \eqref{eq:boundary-nonlinearity-lipschitz} and \eqref{eq:boundary-correction-on-profile} into \eqref{eq:boundary-flow-differential-inequality}, we conclude that as long 
as $X_{j}$ exists and $\|X-X_{j}\| \le 1$, we have \[
\left \| \frac{\dd}{\dd s}(X-X_{j})\right \| \le C\|X-X_{j}\| + C\nu_{-}^{-2j}.
\]
Gr\"onwall's inequality gives
\begin{equation*}
\|X(s)-X_{j}(s)\| \le C(\|Z-Z^{(\ast)}(0)\| + \nu_{-}^{-2j}).
\end{equation*}
After taking $U$ smaller if necessary and taking $J_{0}$ sufficiently large, 
this improves the bootstrap bound to $\|X-X_{j}\| \le \frac{1}{2}$ whenever 
$j \ge J_{0}$. A standard bootstrapping argument then shows that $X_{j}$ 
is well-defined for all $0 \le s \le 2p\delta^{(\ast)}$, and proves 
\eqref{eq:boundary-flow-near-profile}. In particular, $\Phi_{j}^{s}(Z)$ is 
uniformly bounded.

That $\Phi^{s}$ and $\Phi_{j}^{s}$ are $C^{2}$ on $U$ follows from the fact 
that $-I+N$ and $-I+N+B_{j}$ are $C^{2}$ on $\ell^{\infty}$. To bound their 
derivatives and their difference, recall the variational equations:
\begin{equation} \label{eq:boundary-variational-equations}
\begin{aligned}
\frac{\dd}{\dd s}D\Phi^{s}(Z) &= (-I + DN(\Phi^{s}(Z)))D\Phi^{s}(Z), \\
\frac{\dd}{\dd s}D\Phi_{j}^{s}(Z) &= 
(-I + DN(\Phi_{j}^{s}(Z))+DB_{j}(\Phi_{j}^{s}(Z))) D\Phi_{j}^{s}(Z).
\end{aligned}
\end{equation}
As the trajectories $\Phi^{s}(Z)$ and $\Phi_{j}^{s}(Z)$ are uniformly 
bounded, the coefficients $-I+DN$ and $-I+DN+DB_{j}$ of \eqref{eq:boundary-variational-equations} 
are uniformly bounded. Gr\"onwall's inequality gives uniform bounds 
on $D\Phi^{s}(Z)$ and $D\Phi_{j}^{s}(Z)$. Let \[
Y(s) = D\Phi^{s}(Z^{(\ast)}(0)), \quad Y_{j}(s) = D\Phi_{j}^{s}(Z).
\]
Subtracting the two equations in \eqref{eq:boundary-variational-equations}
gives  
\begin{equation*}
\begin{aligned}
\left \| \frac{\dd}{\dd s}(Y-Y_{j}) \right \| &\le 
\|(-I+DN(\Phi^{s}(Z^{(\ast)}(0))))(Y-Y_{j})\| \\
&+ \|(DN(\Phi^{s}(Z^{(\ast)}(0)))-DN(\Phi_{j}^{s}(Z)))Y_{j}\| \\
&+ \|(DB_{j}(Z^{(\ast)}(s))-DB_{j}(\Phi_{j}^{s}(Z)))Y_{j}\| \\
&+ \|DB_{j}(Z^{(\ast)}(s))Y_{j}\| \\
& \le C\|Y-Y_{j}\| + C(\|Z-Z^{(\ast)}(0)\| + \nu_{-}^{-2j}) + C\nu_{-}^{-j}
\end{aligned}
\end{equation*}
using Lipschitz continuity of $DN$ and $DB_{j}$, \eqref{eq:boundary-flow-near-profile}, 
and the bound $\|DB_{j}(Z^{(\ast)}(s))\| \le C\nu_{-}^{-j}$ following 
from the spatial decay of $Z^{(\ast)}$.
Since $Y(0)=Y_{j}(0)=I$, Gr\"onwall's inequality gives \[
\|Y-Y_{j}\| \le C(\|Z-Z^{(\ast)}(0)\| + \nu_{-}^{-j})
\quad (0 \le s \le 2p\delta^{(\ast)}, j \ge J_{0}, Z \in U).
\]
In particular, plugging $Z=Z^{(\ast)}(0)$ above gives 
\eqref{eq:boundary-flow-difference}.

Finally, differentiating \eqref{eq:boundary-variational-equations} once more with respect 
to $Z$ and using uniform bounds on $D^{2}N$, $D^{2}B_{j}$, and the 
first derivatives of the flows, we conclude by Gr\"onwall's inequality that 
$D_{Z}^{2}\Phi^{s}(Z)$ and $D_{Z}^{2}\Phi_{j}^{s}(Z)$ are uniformly 
bounded. The flow equations \eqref{eq:main-ss-Z} and
\eqref{eq:moving-boundary-equation}, together with their derivatives with respect to 
$s$ and $Z$, give the remaining bounds in \eqref{eq:boundary-flow-derivative-bounds}.

This completes the proof of the lemma.
\end{proof}

\subsection{Choosing the modulation parameters}
\label{subsec:boundary-modulation}

We now choose the parameters in \eqref{eq:boundary-normalized-return} 
so that \eqref{eq:boundary-normalization-condition} holds.

\begin{lemma}[Modulation parameters]
\label{lem:boundary-modulation-parameters}
After shrinking $U$ and increasing $J_{0}$ in Lemma 
\ref{lem:boundary-flow-comparison} if necessary, there exist a 
constant $C>0$ and, for every $j \ge J_{0}$, functions 
$a_{j}$, $b_{j}$, and $\theta_{j}$ defined on $U$ such that
\begin{equation}\label{eq:boundary-parameter-normalization}
e^{a_{j}(Z)}T_{\theta_{j}(Z)}R^{-p}\Phi_{j}^{p\delta^{(\ast)}+b_{j}(Z)}(Z) 
\in Z^{(\ast)}(0) + Y \quad (Z \in U, j \ge J_{0}).
\end{equation} 
Here, the codomain of $\theta_{j}$ is $\RR^{2}$ for GOY and $\RR$ for Sabra.
Moreover, we have
\begin{equation}\label{eq:boundary-parameter-bound}
\begin{aligned}
|a_{j}(Z)| + |b_{j}(Z)| + |\theta_{j}(Z)| &\le C(\|Z-Z^{(\ast)}(0)\| + 
\nu_{-}^{-2j}), \\
|a_{j}(Z)| + |b_{j}(Z)| &< p\delta^{(\ast)}
\end{aligned}
\end{equation}
for all $Z \in U$ and $j \ge J_{0}$.
\end{lemma}

\begin{proof}
In the proof, $C$ denotes a constant independent of $j \ge J_{0}$ and 
$Z \in U$, whose value may change from line to line.

Since $Y=\ker\pi$, the required normalization \eqref{eq:boundary-parameter-normalization}
is equivalent to $H_{j}(Z, a_{j}(Z), b_{j}(Z), \theta_{j}(Z))=0$, where 
$H_{j}$ is defined by \[
H_{j}(Z, a, b, \theta):= \pi (e^{a}T_{\theta}R^{-p}\Phi_{j}^{p\delta^{(\ast)}+b}(Z)
-Z^{(\ast)}(0)).
\]

This is well-defined for $Z \in U$ and $|b| < p\delta^{(\ast)}$, by 
Lemma \ref{lem:boundary-flow-comparison}. Moreover, this map is $C^{2}$, as 
$(Z, a, b) \mapsto e^{a}R^{-p}\Phi_{j}^{p\delta^{(\ast)}+b}(Z)$ is $C^{2}$, 
and $\pi$ only reads finitely many coordinates so that $(\theta, W) \mapsto 
\pi T_{\theta}W$ is also $C^{2}$. Together with
\eqref{eq:boundary-flow-derivative-bounds}, this shows that
the first and second derivatives of $H_{j}$ are uniformly bounded in $j \ge J_{0}$, 
$Z \in U$, and $(a, b, \theta)$ in a fixed small neighborhood of zero.

To find $a_{j}(Z), b_{j}(Z), \theta_{j}(Z)$ solving $H_{j}=0$, we use the 
Newton-like iteration
\begin{equation}\label{eq:boundary-newton-map}
(a, b, \theta) \mapsto \mathcal{N}_{Z, j}(a, b, \theta) := (a, b, \theta) - 
[D_{(a, b, \theta)}H_{j}(Z^{(\ast)}(0), 0, 0, 0)]^{-1}H_{j}(Z, a, b, \theta).
\end{equation}

To show that the iteration \eqref{eq:boundary-newton-map} is contractive, note that its derivative 
with respect to $(a, b, \theta)$ is zero at $(Z^{(\ast)}(0), 0, 0, 0)$. 
Thus, the uniform bounds on the second derivatives of $H_{j}$ make the 
iteration \eqref{eq:boundary-newton-map} contractive near $(Z^{(\ast)}(0), 0, 0, 0)$, provided that 
the inverse \[
[D_{(a, b, \theta)}H_{j}(Z^{(\ast)}(0), 0, 0, 0)]^{-1}
\]
is uniformly bounded in $j$.

Observe that
\begin{align}
\partial_{a}H_{j}(Z^{(\ast)}(0), 0, 0, 0) &= \pi (R^{-p}\Phi_{j}^{p\delta^{(\ast)}}
(Z^{(\ast)}(0))) \nonumber \\
& \xrightarrow[j \to \infty]{} 
\pi (R^{-p}\Phi^{p\delta^{(\ast)}}(Z^{(\ast)}(0))) = \pi Z^{(\ast)}(0) 
= \pi (v_{\mathrm{r}}^{u}-v_{\mathrm{r}}^{c}), \label{eq:boundary-amplitude-derivative} \\
\partial_{b}H_{j}(Z^{(\ast)}(0), 0, 0, 0) &= \pi (R^{-p}
(-W+N(W)+B_{j}(W))|_{W=\Phi_{j}^{p\delta^{(\ast)}}(Z^{(\ast)}(0))}) \nonumber \\
&= \pi ( -Z^{(\ast)}(0) + N(Z^{(\ast)}(0)) + B_{j+p}(Z^{(\ast)}(0)))
+ O(\nu_{-}^{-2j}) \nonumber \\
&\xrightarrow[j \to \infty]{} \pi (-Z^{(\ast)}(0) + N(Z^{(\ast)}(0))) 
= \pi v_{\mathrm{r}}^{c}, \label{eq:boundary-time-derivative}
\end{align}
by \eqref{eq:boundary-flow-near-profile}.

For the GOY model, note that 
\begin{equation*}
\begin{aligned}
\partial_{\theta_{1}} H_{j}(Z^{(\ast)}(0), 0, 0, 0) 
&= \partial_{\theta_{1}} 
\pi (T_{\theta_{1}, 0}R^{-p}\Phi_{j}^{p\delta^{(\ast)}}(Z^{(\ast)}(0)))|_{\theta_{1}=0}, \\
\partial_{\theta_{2}} H_{j}(Z^{(\ast)}(0), 0, 0, 0) 
&= \partial_{\theta_{2}} 
\pi (T_{0, \theta_{2}}R^{-p}\Phi_{j}^{p\delta^{(\ast)}}(Z^{(\ast)}(0)))|_{\theta_{2}=0},
\end{aligned}
\end{equation*} 
and differentiating \eqref{eq:GOYphase} gives
\begin{align*}
\partial_{\theta_{1}}(T_{\theta_{1}, 0}W)_{n}|_{\theta_{1}=0} 
&= \begin{cases}
-iW_{n} & n \equiv 0 \pmod{3}, \\
iW_{n} & n \equiv 1 \pmod{3}, \\
0 & n \equiv 2 \pmod{3},
\end{cases} \\
\partial_{\theta_{2}}(T_{0, \theta_{2}}W)_{n}|_{\theta_{2}=0} 
&= \begin{cases}
-iW_{n} & n \equiv 0 \pmod{3}, \\
0 & n \equiv 1 \pmod{3}, \\
iW_{n} & n \equiv 2 \pmod{3}.
\end{cases}
\end{align*}
Since $R^{-p}\Phi_{j}^{p\delta^{(\ast)}}(Z^{(\ast)}(0)) \to Z^{(\ast)}(0)$ in 
$\ell^{\infty}$ as $j \to \infty$, we conclude that
\begin{align}
\partial_{\theta_{1}}H_{j}(Z^{(\ast)}(0), 0, 0, 0) &\xrightarrow[j \to \infty]{} 
\pi v_{\mathrm{G}, 1}^{c}, \label{eq:boundary-goy-phase-derivative-1} \\
\partial_{\theta_{2}}H_{j}(Z^{(\ast)}(0), 0, 0, 0) &\xrightarrow[j \to \infty]{} 
\pi v_{\mathrm{G}, 2}^{c}. \label{eq:boundary-goy-phase-derivative-2}
\end{align}

For the Sabra model, we have \[
\partial_{\theta}H_{j}(Z^{(\ast)}(0), 0, 0, 0) = \partial_{\theta}\pi   
(T_{\theta}R^{-p}\Phi_{j}^{p\delta^{(\ast)}}(Z^{(\ast)}(0)))|_{\theta=0},
\]
and differentiating \eqref{eq:Sabraphase} at each component gives \[
\partial_{\theta}(T_{\theta}W)_{n}|_{\theta =0}
= i \left ( \frac{1+\sqrt{5}}{2}\right )^{n} W_{n}.
\]
Since $\pi$ only reads finitely many coordinates, we similarly conclude that
\begin{equation} \label{eq:boundary-sabra-phase-derivative}
\partial_{\theta}H_{j}(Z^{(\ast)}(0), 0, 0, 0) \xrightarrow[j\to\infty]{} 
\pi v_{\mathrm{S}}^{u}.
\end{equation}

By \eqref{eq:boundary-amplitude-derivative}, \eqref{eq:boundary-time-derivative}, \eqref{eq:boundary-goy-phase-derivative-1}, \eqref{eq:boundary-goy-phase-derivative-2}, and \eqref{eq:boundary-sabra-phase-derivative}, we 
conclude that $D_{(a, b, \theta)}H_{j}(Z^{(\ast)}(0), 0, 0, 0)$ converges in 
operator norm to the map \[
(\alpha, \beta, \gamma) \mapsto \pi \cV(\alpha, \beta-\alpha, \gamma)
\] 
as $j \to \infty$. Since $\pi \cV$ is invertible,
as shown in the proof of Theorem \ref{thm:monodromy-splitting}, this limit 
is an invertible operator. Therefore, after increasing $J_{0}$ if necessary,
we conclude that \begin{equation} \label{eq:boundary-parameter-inverse-bound}
\|[D_{(a, b, \theta)}H_{j}(Z^{(\ast)}(0), 0, 0, 0)]^{-1} \| \le C \quad 
(j \ge J_{0}).
\end{equation}

We now prove that the Newton-like iteration \eqref{eq:boundary-newton-map} converges when 
started at $(0, 0, 0)$, for all large enough $j$. First observe that 
the map \eqref{eq:boundary-newton-map} is contractive near $(Z^{(\ast)}(0), 0, 0, 0)$. 
Indeed, its derivative is
\begin{equation} \label{eq:boundary-newton-derivative}
D_{(a, b, \theta)} \mathcal{N}_{Z, j}(a, b, \theta) 
= I - [D_{(a, b, \theta)}H_{j}(Z^{(\ast)}(0), 0, 0, 0)]^{-1}
D_{(a, b, \theta)}H_{j}(Z, a, b, \theta),
\end{equation}
which is zero at $(Z, a, b, \theta)=(Z^{(\ast)}(0), 0, 0, 0)$. As the 
second derivatives of $H_{j}$ are uniformly bounded in $j$ and 
the inverse is uniformly bounded by \eqref{eq:boundary-parameter-inverse-bound}, for all 
$Z$ near $Z^{(\ast)}(0)$ and $(a, b, \theta)$ near $(0, 0, 0)$, 
the norm of the operator \eqref{eq:boundary-newton-derivative} is at most $\frac{1}{2}$.
Choose $U$ smaller if necessary, and take $0<\epsilon<p\delta^{(\ast)}/2$ 
so that this holds for every $Z \in U$, $j \ge J_{0}$, and 
$\|(a, b, \theta)\| \le \epsilon$. That is,
\begin{equation}\label{eq:boundary-newton-contraction}
\|D_{(a, b, \theta)}\mathcal{N}_{Z, j}(a, b, \theta)\| \le \frac{1}{2}
\quad (Z \in U, j \ge J_{0}, \|(a, b, \theta)\| \le \epsilon).
\end{equation}

If we plug $(a, b, \theta)=(0, 0, 0)$ into \eqref{eq:boundary-newton-map}, we get \[
\mathcal{N}_{Z, j}(0, 0, 0) = - [D_{(a, b, \theta)}H_{j}(Z^{(\ast)}(0), 0, 0, 0)]^{-1}
H_{j}(Z, 0, 0, 0).
\]
By \eqref{eq:boundary-flow-near-profile} and relative periodicity 
of $Z^{(\ast)}$, we have
\begin{equation}\label{eq:boundary-parameter-residual}
\|H_{j}(Z, 0, 0, 0)\| \le C(\|Z-Z^{(\ast)}(0)\| + \nu_{-}^{-2j})
\quad (Z \in U, j \ge J_{0}).
\end{equation}
Together with \eqref{eq:boundary-parameter-inverse-bound}, this gives
\begin{equation}\label{eq:boundary-newton-residual}
\|\mathcal{N}_{Z, j}(0, 0, 0)\| \le C(\|Z-Z^{(\ast)}(0)\| + \nu_{-}^{-2j})
\quad (Z \in U, j \ge J_{0}).
\end{equation}
By \eqref{eq:boundary-newton-residual}, after shrinking $U$ and increasing $J_{0}$ if necessary,
we may assume that \[
\|\mathcal{N}_{Z, j}(0, 0, 0)\| < \frac{\epsilon}{2}.
\]
Then, by \eqref{eq:boundary-newton-contraction}, we have \[
\|\mathcal{N}_{Z, j}(a, b, \theta)\| 
\le \|\mathcal{N}_{Z, j}(0, 0, 0)\| + \frac{1}{2}\|(a, b, \theta)\| \le \epsilon.
\]
Thus, $\mathcal{N}_{Z, j}$ is a contraction from the ball 
$\{(a, b, \theta): \|(a, b, \theta)\| \le \epsilon\}$ to itself.
By the Banach fixed-point theorem, we conclude that there exists 
a unique $(a, b, \theta)$ in this ball 
satisfying $\mathcal{N}_{Z, j}(a, b, \theta)=(a, b, \theta)$, or in other words, \[
H_{j}(Z, a, b, \theta)=0.
\]
Denote this fixed point by $(a_{j}(Z), b_{j}(Z), \theta_{j}(Z))$.
This construction applies to every $Z \in U$ and $j \ge J_{0}$, and 
the fixed point satisfies \[
\|(a_{j}(Z), b_{j}(Z), \theta_{j}(Z))\| \le \epsilon.
\]
Since $\epsilon < p\delta^{(\ast)}/2$, this gives the second inequality of 
\eqref{eq:boundary-parameter-bound}. To prove the first inequality of 
\eqref{eq:boundary-parameter-bound}, note that 
\begin{equation*}
\begin{aligned}
\|(a_{j}(Z), b_{j}(Z), \theta_{j}(Z))\| &= 
\|\mathcal{N}_{Z, j}(a_{j}(Z), b_{j}(Z), \theta_{j}(Z))\| \\    
&\le \|\mathcal{N}_{Z, j}(0, 0, 0)\| 
+ \frac{1}{2}\|(a_{j}(Z), b_{j}(Z), \theta_{j}(Z))\|.
\end{aligned}
\end{equation*}
Therefore, by \eqref{eq:boundary-newton-residual}, we have \[
\|(a_{j}(Z), b_{j}(Z), \theta_{j}(Z))\| \le C 
(\|Z-Z^{(\ast)}(0)\| + \nu_{-}^{-2j}).
\]
After increasing $C$ if necessary, this proves the first
inequality in \eqref{eq:boundary-parameter-bound}.
This completes the proof of the lemma.
\end{proof}

\subsection{Convergence of the modulated iterates}
\label{subsec:boundary-contraction}

Recall that we have defined \[
\widetilde{F}_{j}(Z) = e^{a_{j}(Z)}T_{\theta_{j}(Z)}R^{-p}\Phi_{j}^{p\delta^{(\ast)}
+b_{j}(Z)}(Z).
\]
Let $(\widetilde{Z}^{[J+kp]})_{k \ge 0}$ be a sequence obtained by successively 
applying $\widetilde{F}_{J}, \widetilde{F}_{J+p}, \widetilde{F}_{J+2p}, \ldots$ 
to $\widetilde{Z}^{[J]}$. Our goal is to prove that, for every sufficiently 
large $J$ and every $\widetilde{Z}^{[J]}$ sufficiently close to $Z^{(\ast)}(0)$,
these iterates are well-defined and converge to $Z^{(\ast)}(0)$.
We first use the contractivity of $\widetilde{M}^{p}$ to estimate the size 
of the perturbation after applying the modulated return map 
$\widetilde{F}_{j}$.

\begin{proposition}\label{prop:boundary-normalized-return}
After shrinking $U$ and increasing $J_{0}$ if necessary, there are 
constants $C>0$, $\epsilon_{0} > 0$, and $0<\kappa<1$ such that
\begin{equation}\label{eq:normalized-return-contraction}
\|\widetilde{F}_{j}(Z^{(\ast)}(0)+y)-Z^{(\ast)}(0)\| 
\le \kappa \|y\| + C\nu_{-}^{-2j} \quad (j \ge J_{0}, y \in Y, 
\|y\| \le \epsilon_{0}).
\end{equation}
We also have
\begin{equation}\label{eq:boundary-first-return}
\|\widetilde{F}_{j}(Z) - Z^{(\ast)}(0)\| \le C(\|Z-Z^{(\ast)}(0)\| + 
\nu_{-}^{-2j})
\end{equation}
for every $Z \in U$ and $j \ge J_{0}$.
\end{proposition}

\begin{proof}
In the proof, $C$ denotes a constant independent of $j \ge J_{0}$ and 
$Z \in U$, whose value may change from line to line.

In the Sabra case, the map $(\theta, Z) \mapsto T_{\theta}Z$ is not a $C^{1}$ 
function. We begin by moving the phase transformation $T_{\theta}$ onto 
$Z^{(\ast)}(0)$, which has sufficient spatial decay so that 
$\theta \mapsto T_{\theta}Z^{(\ast)}(0)$ is $C^{2}$. Define \[
E_{j}(Z) := T_{-\theta_{j}(Z)}(\widetilde{F}_{j}(Z) - Z^{(\ast)}(0)) 
= e^{a_{j}(Z)}R^{-p}\Phi_{j}^{p\delta^{(\ast)}+b_{j}(Z)}(Z) 
- T_{-\theta_{j}(Z)}Z^{(\ast)}(0).
\]
Since the phase transformation $T_{\theta}$ preserves the $\ell^{\infty}$ norm, 
it suffices to estimate $\|E_{j}(Z)\|$. 

Next, to study $e^{a_{j}(Z)}R^{-p}\Phi_{j}^{p\delta^{(\ast)}+b_{j}(Z)}(Z)$,
recall that by Lemma \ref{lem:boundary-modulation-parameters}, 
the parameters $a_{j}(Z)$ and $b_{j}(Z)$ are small for $Z \in U$. 
We apply Taylor's formula to \[
(Z, a, b) \mapsto G_{j}(Z, a, b) := e^{a}R^{-p}\Phi_{j}^{p\delta^{(\ast)}+b}(Z).
\]
By Lemma \ref{lem:boundary-flow-comparison}(i), $G_{j}$ is $C^{2}$, and its 
first and second derivatives are uniformly bounded in $j\ge J_{0}$, 
$Z \in U$, and $(a, b)$ in a fixed small neighborhood of zero. 
\begin{align*}
& \|\partial_{a}G_{j}(Z^{(\ast)}(0), 0, 0) - (v_{\mathrm{r}}^{u}-v_{\mathrm{r}}^{c})\| \\
&= \|R^{-p}\Phi_{j}^{p\delta^{(\ast)}}(Z^{(\ast)}(0)) - Z^{(\ast)}(0) \| \\
&\le \|\Phi_{j}^{p\delta^{(\ast)}}(Z^{(\ast)}(0)) - \Phi^{p\delta^{(\ast)}}
(Z^{(\ast)}(0)) \| \\
& \le C \nu_{-}^{-2j} \quad \text{by \eqref{eq:boundary-flow-near-profile}}, \\
& \|\partial_{b}G_{j}(Z^{(\ast)}(0), 0, 0) - v_{\mathrm{r}}^{c}\| \\
&= 
\|R^{-p}(-W+N(W)+B_{j}(W))|_{W=\Phi_{j}^{p\delta^{(\ast)}}(Z^{(\ast)}(0))}
-(-Z^{(\ast)}(0) + N(Z^{(\ast)}(0))) \| \\
& \le \|\Phi_{j}^{p\delta^{(\ast)}}(Z^{(\ast)}(0)) - Z^{(\ast)}(p\delta^{(\ast)})\| + 
\| N(\Phi_{j}^{p\delta^{(\ast)}}(Z^{(\ast)}(0))) - N(Z^{(\ast)}(p\delta^{(\ast)}))\| 
\\
& + \|B_{j}(\Phi_{j}^{p\delta^{(\ast)}}(Z^{(\ast)}(0)))
-B_{j}(Z^{(\ast)}(p\delta^{(\ast)}))\| + \| B_{j}(Z^{(\ast)}(p\delta^{(\ast)}))\| \\
& \le C\nu_{-}^{-2j} \quad \text{by \eqref{eq:boundary-flow-near-profile}, \eqref{eq:boundary-correction-on-profile}}, \\
& \|D_{Z}G_{j}(Z^{(\ast)}(0), 0, 0) - M^{p}\| \\ 
&= \|R^{-p}D\Phi_{j}^{p\delta^{(\ast)}}(Z^{(\ast)}(0)) - (R^{-1}D\Phi^{\delta^{(\ast)}}
(Z^{(\ast)}(0)))^{p}\| \\
&= \|R^{-p}D\Phi_{j}^{p\delta^{(\ast)}}(Z^{(\ast)}(0)) - R^{-p}D\Phi^{p\delta^{(\ast)}}
(Z^{(\ast)}(0))\| \\
&= \|D\Phi_{j}^{p\delta^{(\ast)}}(Z^{(\ast)}(0)) 
- D\Phi^{p\delta^{(\ast)}}(Z^{(\ast)}(0))\| \\
& \le C\nu_{-}^{-j} \quad \text{by \eqref{eq:boundary-flow-difference}}.
\end{align*}
Therefore, by Taylor's formula applied to $G_{j}$, we have
\begin{align}
& \|e^{a}R^{-p}\Phi_{j}^{p\delta^{(\ast)}+b}(Z) - Z^{(\ast)}(0) 
- a(v_{\mathrm{r}}^{u}-v_{\mathrm{r}}^{c}) - bv_{\mathrm{r}}^{c}
- M^{p}(Z-Z^{(\ast)}(0))\| \nonumber \\
& \le C\nu_{-}^{-2j} + C\nu_{-}^{-j}(|a|+|b|+\|Z-Z^{(\ast)}(0)\|) 
+ C(|a|^{2}+|b|^{2}+\|Z-Z^{(\ast)}(0)\|^{2}) \label{eq:boundary-flow-taylor-remainder}
\end{align}
for all $Z \in U$, $j \ge J_{0}$, and $(a, b)$ in a fixed sufficiently small 
neighborhood of zero.

Next, the phase transformations satisfy
\begin{equation}\label{eq:boundary-profile-phase-remainder}
\begin{cases}
\|T_{-\theta_{1}, -\theta_{2}}Z^{(\ast)}(0)-Z^{(\ast)}(0)
+\theta_{1}v_{\mathrm{G}, 1}^{c} + \theta_{2}v_{\mathrm{G}, 2}^{c}\| \le C|\theta|^{2},
& \text{for GOY}, \\
\|T_{-\theta}Z^{(\ast)}(0) - Z^{(\ast)}(0) + \theta v_{\mathrm{S}}^{u}\| \le 
C|\theta|^{2}, & \text{for Sabra}.
\end{cases}
\end{equation}
For the GOY model, this follows directly from differentiating the 
$C^{2}$ function $\theta \mapsto T_{\theta}Z$ given in \eqref{eq:GOYphase}. 
For the Sabra model, the right spatial decay in Theorem \ref{thm:orbit-existence} 
gives \[
\sum_{n\in \ZZ} \left ( \frac{1+\sqrt{5}}{2}\right )^{2n} |(Z^{(\ast)}(0))_{n}| < \infty,
\]
so the estimate follows from differentiating the $C^{2}$ function 
$\theta \mapsto T_{\theta}Z^{(\ast)}(0)$ given in \eqref{eq:Sabraphase}.

We now substitute $a_{j}(Z)$ and $b_{j}(Z)$ for $a$ and $b$
in \eqref{eq:boundary-flow-taylor-remainder}.
The estimate \eqref{eq:boundary-parameter-bound} ensures
that these values lie in the required neighborhood of zero,
after shrinking $U$ and increasing $J_{0}$ if necessary.
Combining \eqref{eq:boundary-flow-taylor-remainder},
\eqref{eq:boundary-profile-phase-remainder}, and
\eqref{eq:boundary-parameter-bound}, we obtain
\begin{equation}\label{eq:boundary-comparison-expansion}
\begin{aligned}
& \|E_{j}(Z) - M^{p}(Z-Z^{(\ast)}(0)) - \cV(a_{j}(Z), b_{j}(Z)-a_{j}(Z), 
\theta_{j}(Z)) \| \\ 
& \le C\nu_{-}^{-2j} + C\nu_{-}^{-j}(|a_{j}(Z)| + |b_{j}(Z)| + \|Z-Z^{(\ast)}(0)\|) \\
& + C(|a_{j}(Z)|^{2} + |b_{j}(Z)|^{2} + \|Z-Z^{(\ast)}(0)\|^{2})
+ C|\theta_{j}(Z)|^{2} \\
& \le C ( \|Z-Z^{(\ast)}(0)\|^{2} + \nu_{-}^{-2j})
\end{aligned}
\end{equation}
for every $Z \in U$ and $j \ge J_{0}$.

Since $M^{p}$ is bounded and $\cV$ is bounded, 
\eqref{eq:boundary-comparison-expansion} and \eqref{eq:boundary-parameter-bound} 
give \[
\|E_{j}(Z)\| \le C\|Z-Z^{(\ast)}(0)\| + C(\|Z-Z^{(\ast)}(0)\|+\nu_{-}^{-2j}) 
+ C(\|Z-Z^{(\ast)}(0)\|^{2}+\nu_{-}^{-2j})
\]
for every $Z \in U$ and $j \ge J_{0}$. This proves \eqref{eq:boundary-first-return}.

Now assume that $Z=Z^{(\ast)}(0)+y$ with $y \in Y$. By the
definition of $E_{j}(Z)$ and \eqref{eq:boundary-parameter-normalization}, we have \[
\pi T_{\theta_{j}(Z)}E_{j}(Z) = \pi (\widetilde{F}_{j}(Z)-Z^{(\ast)}(0)) = 0.
\]
Since $\Pi=\cV(\pi \cV)^{-1}\pi$, we also have 
$\Pi T_{\theta_{j}(Z)}E_{j}(Z)=0$. Therefore, we have
\begin{align}
\|E_{j}(Z)\| & = \|E_{j}(Z) - \Pi T_{\theta_{j}(Z)}E_{j}(Z)\| \nonumber \\
& \le \|(I-\Pi)E_{j}(Z)\| + \|\Pi(I-T_{\theta_{j}(Z)})E_{j}(Z)\| \nonumber \\
& \le \|(I-\Pi)(E_{j}(Z) - M^{p}(Z-Z^{(\ast)}(0))-\cV
(a_{j}(Z), b_{j}(Z)-a_{j}(Z), \theta_{j}(Z)))\| \nonumber \\
& + \|(I-\Pi)M^{p}(Z-Z^{(\ast)}(0))\| + \|\Pi(I-T_{\theta_{j}(Z)})E_{j}(Z)\|,
\label{eq:boundary-projection-estimate}
\end{align}
where the last inequality follows from $(I-\Pi)\cV=0$.
Recall that $V$ is $M$-invariant, so that $(I-\Pi)M\Pi = (I-\Pi)\Pi M =0$. 
Therefore, 
\begin{equation} \label{eq:boundary-projected-power}
(I-\Pi)M^{p}|_{Y} = ((I-\Pi)M|_{Y})^{p} = \widetilde{M}^{p}.
\end{equation} 
Also, since $\pi$ reads only finitely many input coordinates, 
the map $\theta \mapsto \pi (I-T_{\theta})$ is smooth in operator norm 
and vanishes at $\theta=0$. Thus, by $\Pi = \cV(\pi \cV)^{-1}\pi$ 
and \eqref{eq:boundary-parameter-bound}, we have
\begin{equation} \label{eq:boundary-projected-phase-bound}
\|\Pi(I-T_{\theta_{j}(Z)})\| \le C|\theta_{j}(Z)| \le 
C(\|Z-Z^{(\ast)}(0)\|+\nu_{-}^{-2j}).
\end{equation} 
Using \eqref{eq:boundary-comparison-expansion}, \eqref{eq:boundary-projected-power}, 
\eqref{eq:boundary-projected-phase-bound},
and \eqref{eq:boundary-first-return} in \eqref{eq:boundary-projection-estimate}, we obtain \[
\|E_{j}(Z)\| \le C(\|Z-Z^{(\ast)}(0)\|^{2}+\nu_{-}^{-2j}) 
+ \|\widetilde{M}^{p}\|_{Y \to Y}\|Z-Z^{(\ast)}(0)\|
+ C(\|Z-Z^{(\ast)}(0)\| + \nu_{-}^{-2j})^{2}.
\]
Choose $\kappa$ such that $\|\widetilde{M}^{p}\|_{Y \to Y} < \kappa < 1$,
which is possible by Theorem \ref{thm:monodromy-splitting}. By taking 
$\epsilon_{0}$ sufficiently small, we conclude that \[
\|E_{j}(Z)\| \le C (\|y\|^{2} + \nu_{-}^{-2j}) + \|\widetilde{M}^{p}\|_{Y \to Y}
\|y\| \le \kappa \|y\| + C\nu_{-}^{-2j}.
\]
This proves \eqref{eq:normalized-return-contraction} and completes the 
proof of the proposition.
\end{proof}

The estimate \eqref{eq:boundary-first-return} permits the initial iterate 
$\widetilde{Z}^{[J]}$ to be any sufficiently small $\ell^{\infty}$
perturbation of $Z^{(\ast)}(0)$. The normalization condition 
\eqref{eq:boundary-parameter-normalization} ensures that the first 
modulated iterate $\widetilde{Z}^{[J+p]}$ lies in $Z^{(\ast)}(0)+Y$. 
We can then apply \eqref{eq:normalized-return-contraction} at every 
subsequent step. This gives the following corollary.

\begin{corollary}\label{cor:boundary-normalized-iteration}
There exists a constant $\epsilon > 0$ such that, for every sufficiently
large $J$ and every initial state $\widetilde{Z}^{[J]}$ satisfying \[
\|\widetilde{Z}^{[J]} - Z^{(\ast)}(0)\| < \epsilon,
\]
the iteration
\begin{equation}\label{eq:normalized-iteration}
\widetilde{Z}^{[J+(k+1)p]} := \widetilde{F}_{J+kp}(\widetilde{Z}^{[J+kp]})
\quad (k \ge 0)
\end{equation}
is well-defined for all $k$ and remains in $U$, and converges to 
$Z^{(\ast)}(0)$ in $\ell^{\infty}$.
\end{corollary}

\begin{proof}
Let $C>0$, $\epsilon_{0}>0$, and $0<\kappa<1$ be as in Proposition 
\ref{prop:boundary-normalized-return}. Choose $\epsilon>0$ small enough
and $J$ large enough so that \[
C(\epsilon + \nu_{-}^{-2J}) \le \frac{1}{2}\epsilon_{0}, 
\quad \frac{1}{2}\kappa \epsilon_{0} + C \nu_{-}^{-2J} \le \frac{1}{2}\epsilon_{0},
\quad \epsilon < \epsilon_{0}.
\]
Then Proposition \ref{prop:boundary-normalized-return} implies by induction that \[
\|\widetilde{Z}^{[J+kp]} - Z^{(\ast)}(0)\| \le \frac{1}{2}\epsilon_{0} 
\quad (k \ge 1),
\]
so that we can use \eqref{eq:normalized-return-contraction} for all $k \ge 1$. 
Iterating \eqref{eq:normalized-return-contraction} gives
\begin{equation*}
\begin{aligned}
\|\widetilde{Z}^{[J+kp]} - Z^{(\ast)}(0)\| & \le 
\frac{1}{2}\epsilon_{0} \kappa^{k-1} + \sum_{m=1}^{k-1} C\nu_{-}^{-2(J+mp)}
\kappa^{k-1-m} \\
& \le \frac{1}{2}\epsilon_{0}\kappa^{k-1} + C\nu_{-}^{-2J}(k-1)
\max\{\kappa, \nu_{-}^{-2p}\}^{k-1} \xrightarrow[k \to \infty]{} 0.
\end{aligned}
\end{equation*}
This completes the proof.
\end{proof}

\subsection{Reconstruction in the original variables}
\label{subsec:boundary-reconstruction}

We now turn the sequence in \eqref{eq:normalized-iteration}
into a continuous solution of $\partial_{t}z=N(z)+B(z)$ and
restrict it to the positive shells to get a blowup solution 
of $\partial_{t}w = N^{+}(w)$.

\begin{proposition} \label{prop:boundary-reconstruction}
Let $\epsilon > 0$ and $J$ be as in Corollary 
\ref{cor:boundary-normalized-iteration}. Then for any 
initial data $z(0) \in \ell^{\infty}$ satisfying \[
\|z(0) - R^{J}Z^{(\ast)}(0)\| < \epsilon,
\]
there exists a finite time $T>0$ such that the solution 
$z \in C_{t}^{1}([0, T); \ell^{\infty})$ of $\partial_{t}z = N(z)+B(z)$
satisfies \begin{equation} \label{eq:boundary-reconstructed-blowup-rate}
0 < \liminf_{t \to T-}(T-t)\|z(t)\|_{\ell^{\infty}} 
\le \limsup_{t \to T-}(T-t)\|z(t)\|_{\ell^{\infty}} < \infty.
\end{equation}
Moreover, if we define $w=(z_{n})_{n \ge 1}$ to be the restriction of 
$z$ to positive shell indices, then 
$w \in C_{t}^{1}([0, T); \ell^{\infty}(\ZZ_{>0}; \CC))$ solves 
\eqref{eq:main-renorm-common} and satisfies \[
0 < \liminf_{t \to T-}(T-t)\|w(t)\|_{\ell^{\infty}(\ZZ_{>0}; \CC)} 
\le \limsup_{t \to T-}(T-t)\|w(t)\|_{\ell^{\infty}(\ZZ_{>0}; \CC)} < \infty.
\] 
Finally, there exists a function 
$v \in L^{\infty}([0, T); \ell^{1}(\ZZ_{>0}; \CC))$ such that
\begin{equation} \label{eq:boundary-rescaled-profile-approximation}
\|(T-t)w(t)-v(t)\|_{\ell^{\infty}(\ZZ_{>0}; \CC)}
\end{equation}
converges to zero as $t$ approaches $T-$.
\end{proposition}

\begin{proof}
Let $\widetilde{Z}^{[J]} = R^{-J}z(0)$. Define $(\widetilde{Z}^{[J+kp]})_{k \ge 0}$ 
as in Corollary \ref{cor:boundary-normalized-iteration}. Since \[
\|\widetilde{Z}^{[J]} - Z^{(\ast)}(0)\| = \|z(0) - R^{J}Z^{(\ast)}(0)\| < \epsilon,
\]
the corollary shows that the modulated iterates $(\widetilde{Z}^{[J+kp]})_{k \ge 0}$ 
converge to $Z^{(\ast)}(0)$.

For brevity, write \[
a_{J+kp}:=a_{J+kp}(\widetilde{Z}^{[J+kp]}), \quad
b_{J+kp}:=b_{J+kp}(\widetilde{Z}^{[J+kp]}), \quad
\theta_{J+kp}:=\theta_{J+kp}(\widetilde{Z}^{[J+kp]}).
\]
By \eqref{eq:boundary-parameter-bound} and the convergence of
$(\widetilde{Z}^{[J+kp]})_{k \ge 0}$, these parameters tend 
to zero as $k \to \infty$. 

Then $\widetilde{Z}^{[J+(k+1)p]} = \widetilde{F}_{J+kp}(\widetilde{Z}^{[J+kp]})$ 
can be rewritten as \[
\Phi_{J+kp}^{p\delta^{(\ast)}+b_{J+kp}}(\widetilde{Z}^{[J+kp]})
= e^{-a_{J+kp}}R^{p}T_{-\theta_{J+kp}}\widetilde{Z}^{[J+(k+1)p]}.
\]
Returning to the original variables gives a solution $\widetilde{z}^{[J+kp]}$ of 
$\partial_{t}\widetilde{z}^{[J+kp]} = N(\widetilde{z}^{[J+kp]})
+B(\widetilde{z}^{[J+kp]})$,
\begin{equation}\label{eq:boundary-original-segment}
\widetilde{z}^{[J+kp]}(1-e^{-s}) 
= e^{s}R^{J+kp}\Phi_{J+kp}^{s}(\widetilde{Z}^{[J+kp]})
\quad (0 \le s \le p\delta^{(\ast)} + b_{J+kp}),
\end{equation}
where the solution starts from $R^{J+kp}\widetilde{Z}^{[J+kp]}$ at time $t=0$
and reaches \[
e^{p\delta^{(\ast)}+b_{J+kp}-a_{J+kp}}R^{J+(k+1)p}T_{-\theta_{J+kp}}
\widetilde{Z}^{[J+(k+1)p]}
\]
at time $t=1-e^{-p\delta^{(\ast)}-b_{J+kp}}$. 

To match this endpoint of $\widetilde{z}^{[J+kp]}$ with the initial value of 
$\widetilde{z}^{[J+(k+1)p]}$, we need to apply \[
S_{J+(k+1)p} := R^{J+(k+1)p}T_{-\theta_{J+kp}}R^{-(J+(k+1)p)}
\]
and multiply its amplitude by \[
A_{J+(k+1)p}:=e^{p\delta^{(\ast)}+b_{J+kp}-a_{J+kp}}.
\]
Note that applying $S_{J+(k+1)p}$ to a 
solution of $\partial_{t}z = N(z)+B(z)$ still gives a solution to the same 
equation, and $z(t) \mapsto Az(At)$ also gives a solution to the same equation. 
Thus, if we define 
\begin{equation}\label{eq:boundary-physical-reconstruction}
z^{[J+kp]}(t) = \left ( \prod_{m=1}^{k}A_{J+mp}\right ) 
\left ( S_{J+p} \circ S_{J+2p} \circ \cdots \circ S_{J+kp}\right )
\widetilde{z}^{[J+kp]}\left ( \left ( \prod_{m=1}^{k}A_{J+mp}\right ) t\right ),
\end{equation}
then this gives a solution to $\partial_{t}z = N(z)+B(z)$ on the time interval \[
\left [ 0, \left ( \prod_{m=1}^{k}A_{J+mp}\right )^{-1} 
\left ( 1-e^{-p\delta^{(\ast)}-b_{J+kp}}\right )\right ],
\]
such that the endpoint of $z^{[J+kp]}$ matches the initial value of 
$z^{[J+(k+1)p]}$.

The bound \eqref{eq:boundary-parameter-bound} ensures that all these 
time intervals have positive lengths.
Therefore, gluing all the constructed solutions $z^{[J+kp]}$ gives a solution $z$
to $\partial_{t}z = N(z)+B(z)$ with initial data \[
z(0) = z^{[J]}(0) = \widetilde{z}^{[J]}(0) = R^{J}\widetilde{Z}^{[J]},
\]
on a time interval $[0, T)$ where \[
T = \sum_{k=0}^{\infty} \left ( \prod_{m=1}^{k}A_{J+mp}\right )^{-1} 
\left ( 1-e^{-p\delta^{(\ast)}-b_{J+kp}}\right ).
\]

Since the parameters $a_{J+kp}$ and $b_{J+kp}$ converge to zero as $k \to \infty$, 
we have \[
e^{\frac{1}{2}p\delta^{(\ast)}} < A_{J+(k+1)p} < e^{\frac{3}{2}p\delta^{(\ast)}}
\]
for all large enough $k$, and therefore \[
\prod_{m=1}^{k}A_{J+mp} >  e^{\frac{1}{4}kp\delta^{(\ast)}}
\] for all large enough $k$. Thus all but finitely many terms in the series 
defining $T$ are bounded by the corresponding terms of the 
convergent geometric series \[
\sum_{k=0}^{\infty} e^{-\frac{1}{4}kp\delta^{(\ast)}}.
\]
Therefore $T<\infty$.

Restricting the solution $z$ to positive shell indices, we obtain a solution 
$w=(z_{n})_{n \ge 1}$ to \eqref{eq:main-renorm-common}.

Now we prove the $\ell^{\infty}$ blowup rate for $z$ and $w$. 
As $\|w\|_{\ell^{\infty}(\ZZ_{>0}; \CC)} \le \|z\|_{\ell^{\infty}}$, it 
is enough to prove that \[
c(T-t)^{-1} \le \|w(t)\|_{\ell^{\infty}(\ZZ_{>0}; \CC)} 
\le \|z(t)\|_{\ell^{\infty}(\ZZ; \CC)} \le C(T-t)^{-1}
\]
for all $t$ sufficiently close to $T$. Define $w^{[J+kp]}$ and 
$\widetilde{w}^{[J+kp]}$ to be restrictions of $z^{[J+kp]}$
and $\widetilde{z}^{[J+kp]}$ to positive shell indices.

Let
\begin{equation}\label{eq:boundary-physical-times}
t_{J+kp}:= \sum_{r=0}^{k-1} \left ( \prod_{m=1}^{r}A_{J+mp}\right )^{-1} 
\left ( 1 - e^{-p\delta^{(\ast)}-b_{J+rp}}\right ).
\end{equation}
By construction of the solution $z$, we have \[
z(t)=z^{[J+kp]}(t-t_{J+kp}) \quad (t \in [t_{J+kp}, t_{J+(k+1)p}]).
\]
Therefore it is enough to prove that
\begin{equation*}
\begin{aligned}
c(T-t)^{-1} &\le \|w^{[J+kp]}(t-t_{J+kp})\|_{\ell^{\infty}(\ZZ_{>0}; \CC)} \\
&\le \|z^{[J+kp]}(t-t_{J+kp})\|_{\ell^{\infty}(\ZZ; \CC)} \le C(T-t)^{-1}
\quad (t \in [t_{J+kp}, t_{J+(k+1)p}])
\end{aligned}
\end{equation*}
for all large $k$ and for some constants $c, C>0$ independent of $k$.

Note that for large enough $k$,
\begin{align*}
T-t_{J+kp} &= \sum_{r=k}^{\infty} \left ( \prod_{m=1}^{r} A_{J+mp} \right )^{-1} 
(1-e^{-p\delta^{(\ast)}-b_{J+rp}}) \\
&= \left ( \prod_{m=1}^{k} A_{J+mp}\right )^{-1} 
\sum_{r=k}^{\infty} \left ( \prod_{m=k+1}^{r} A_{J+mp}\right )^{-1}
(1-e^{-p\delta^{(\ast)}-b_{J+rp}}) \\
& \le  \left ( \prod_{m=1}^{k} A_{J+mp}\right )^{-1} 
\sum_{r=k}^{\infty} e^{-\frac{1}{2}p\delta^{(\ast)}(r-k)}(1-e^{-2p\delta^{(\ast)}}) \\
& \le C \left ( \prod_{m=1}^{k} A_{J+mp}\right )^{-1}, \\
T-t_{J+kp} &= \left ( \prod_{m=1}^{k} A_{J+mp}\right )^{-1} 
\sum_{r=k}^{\infty} \left ( \prod_{m=k+1}^{r} A_{J+mp}\right )^{-1}
(1-e^{-p\delta^{(\ast)}-b_{J+rp}}) \\
& \ge \left ( \prod_{m=1}^{k} A_{J+mp}\right )^{-1}
\sum_{r=k}^{\infty} e^{-\frac{3}{2}p\delta^{(\ast)}(r-k)}
(1-e^{-\frac{1}{2}p\delta^{(\ast)}}) \\
& \ge c \left ( \prod_{m=1}^{k} A_{J+mp}\right )^{-1}.
\end{align*}
As we have $T-t_{J+(k+1)p} \le T-t \le T-t_{J+kp}$ for 
$t \in [t_{J+kp}, t_{J+(k+1)p}]$, we conclude that 
\begin{equation} \label{eq:boundary-remaining-time}
c \prod_{m=1}^{k}A_{J+mp} \le (T-t)^{-1} \le C\prod_{m=1}^{k}A_{J+mp}
\quad (t \in [t_{J+kp}, t_{J+(k+1)p}]).
\end{equation}
Here, we have absorbed one term $A_{J+(k+1)p}$ into the constants.
Thus, we are left to prove that
\begin{equation*}
\begin{aligned}
c \prod_{m=1}^{k}A_{J+mp} & \le 
\|w^{[J+kp]}(t-t_{J+kp})\|_{\ell^{\infty}(\ZZ_{>0}; \CC)} \\ 
&\le \|z^{[J+kp]}(t-t_{J+kp})\|_{\ell^{\infty}(\ZZ; \CC)} 
\le C\prod_{m=1}^{k}A_{J+mp}
\quad (t \in [t_{J+kp}, t_{J+(k+1)p}]).
\end{aligned}
\end{equation*}

Note that in the definition \eqref{eq:boundary-physical-reconstruction}
of $z^{[J+kp]}$, $S_{J+mp}$ preserves the modulus of each component.
Therefore, \[
\|z^{[J+kp]}(t)\|_{\ell^{\infty}(\ZZ; \CC)} 
= \left ( \prod_{m=1}^{k}A_{J+mp}\right ) \left \| \widetilde{z}^{[J+kp]} 
\left ( \left ( \prod_{m=1}^{k}A_{J+mp}\right ) t \right )\right \|_{\ell^{\infty}
(\ZZ; \CC)},
\]
and the same equality holds for the restriction to positive shell indices.
Therefore it remains to prove that
\begin{equation} \label{eq:boundary-original-segment-bounds}
c \le \left \| \widetilde{w}^{[J+kp]}(1-e^{-s}) \right \|_{\ell^{\infty}
(\ZZ_{>0}; \CC)}
\le \left \| \widetilde{z}^{[J+kp]}(1-e^{-s}) \right \|_{\ell^{\infty}
(\ZZ; \CC)} \le C
\end{equation}
for all large enough $k$, where the self-similar time $s$ ranges over 
$[0, p\delta^{(\ast)}+b_{J+kp}]$. In the definition \eqref{eq:boundary-original-segment} of 
$\widetilde{z}^{[J+kp]}$, we may ignore the $e^{s}$ term as it can be 
absorbed into the constants. Now the upper bound in \eqref{eq:boundary-original-segment-bounds} follows as 
$\Phi_{J+kp}^{s}$ is uniformly bounded near $Z^{(\ast)}(0)$ and 
$\widetilde{Z}^{[J+kp]}$ is converging to $Z^{(\ast)}(0)$. To prove the 
lower bound, observe that 
\begin{align}
& \|((R^{J+kp}\Phi_{J+kp}^{s}(\widetilde{Z}^{[J+kp]}))_{n})_{n \ge 1}\|_{
\ell^{\infty}(\ZZ_{>0}; \CC)} \nonumber \\
&= \sup_{n \ge 1-J-kp} \left |\left (\Phi_{J+kp}^{s}(\widetilde{Z}^{[J+kp]})
\right )_{n}\right | \nonumber \\
&\ge \sup_{n \ge 1-J-kp} \left |(Z^{(\ast)}(s))_{n}\right |
- \|Z^{(\ast)}(s) - \Phi_{J+kp}^{s}(\widetilde{Z}^{[J+kp]})\| \nonumber \\
& \ge \|Z^{(\ast)}(s)\| - \sup_{n \le -J-kp} \left |(Z^{(\ast)}(s))_{n}\right | 
- \|Z^{(\ast)}(s) - \Phi_{J+kp}^{s}(\widetilde{Z}^{[J+kp]})\|.
\label{eq:boundary-positive-shell-lower-bound}
\end{align}
Here, $\|Z^{(\ast)}(s)\|$ can be bounded below by a positive number, as the 
bilateral travelling wave $Z^{(\ast)}$ is continuous and is nonvanishing 
by uniqueness. The second term on the right-hand side of \eqref{eq:boundary-positive-shell-lower-bound}, \[
\sup_{n \le -J-kp} \left | (Z^{(\ast)}(s))_{n}\right |,
\]
converges to zero as $k \to \infty$, by the left spatial decay from Theorem 
\ref{thm:orbit-existence}. The last term of \eqref{eq:boundary-positive-shell-lower-bound} can be bounded using 
\eqref{eq:boundary-flow-near-profile}, so that \[
 \|Z^{(\ast)}(s) - \Phi_{J+kp}^{s}(\widetilde{Z}^{[J+kp]})\|
 \le C(\|\widetilde{Z}^{[J+kp]}-Z^{(\ast)}(0)\| + \nu_{-}^{-2(J+kp)})
 \xrightarrow[k \to \infty]{} 0.
\]
Since both error terms tend uniformly to zero, \eqref{eq:boundary-positive-shell-lower-bound} proves the lower 
bound in \eqref{eq:boundary-original-segment-bounds}. This completes the proof of
the $\ell^{\infty}$ blowup rate for both $z$ and $w$.

It remains to construct $v \in L^{\infty}([0, T); \ell^{1}
(\ZZ_{>0}; \CC))$ such that 
\eqref{eq:boundary-rescaled-profile-approximation} converges to zero as 
$t \to T-$. We will define 
$\widehat{v}$ in the same way as we defined $z$, except that we use 
$Z^{(\ast)}(s)$ instead of $\Phi_{J+kp}^{s}(\widetilde{Z}^{[J+kp]})$ in 
\eqref{eq:boundary-original-segment}. That is, define \[
\widetilde{v}^{[J+kp]}(1-e^{-s}) = e^{s}R^{J+kp}Z^{(\ast)}(s) 
\quad (0 \le s < p\delta^{(\ast)}+b_{J+kp}),
\]
and define \[
\widehat{v}(t+t_{J+kp}) = \left ( \prod_{m=1}^{k}A_{J+mp}\right )    
(S_{J+p} \circ S_{J+2p} \circ \cdots \circ S_{J+kp}) \widetilde{v}^{[J+kp]} 
\left ( \left ( \prod_{m=1}^{k} A_{J+mp} \right ) t \right ).
\]
Then this gives a function $\widehat{v}$ defined on each time interval 
$[t_{J+kp}, t_{J+(k+1)p})$, and hence on $[0, T)$. Now define \[
v(t) = (T-t)((\widehat{v}(t))_{n})_{n \ge 1} \quad (0 \le t < T).
\]

To prove that $v \in L^{\infty}_{t}([0, T); \ell^{1}(\ZZ_{>0}; \CC))$, 
observe that $S_{j}$ does not change the $\ell^{1}(\ZZ_{>0}; \CC)$ norm 
of a function, and thus by \eqref{eq:boundary-remaining-time} it is 
enough to check that each $\widetilde{v}^{[J+kp]}$ is bounded in 
$\ell^{1}(\ZZ; \CC)$ uniformly in $k$. And this is true since 
$Z^{(\ast)}(s)$ has exponential spatial decay by Theorem 
\ref{thm:orbit-existence}.

To prove that \eqref{eq:boundary-rescaled-profile-approximation} converges 
to zero, note that for $t_{J+kp} \le t < t_{J+(k+1)p}$ we have
\begin{align*}
&\|(T-t)w(t)-v(t)\|_{\ell^{\infty}(\ZZ_{>0}; \CC)} \\
&\le (T-t)\|z(t)-\widehat{v}(t)\|_{
  \ell^{\infty}(\ZZ; \CC)} \\ 
&\le C \left \|\widetilde{z}^{[J+kp]}\left ( \left ( \prod_{m=1}^{k}A_{J+mp}\right )
(t-t_{J+kp})  
\right ) - \widetilde{v}^{[J+kp]}\left (\left (\prod_{m=1}^{k}A_{J+mp}\right )
(t-t_{J+kp})\right )   
\right \|_{\ell^{\infty}(\ZZ; \CC)} \\
&\le C\|\Phi_{J+kp}^{s}(\widetilde{Z}^{[J+kp]}) - Z^{(\ast)}(s)\|_{\ell^{\infty}} \\
&\le C(\|\widetilde{Z}^{[J+kp]}-Z^{(\ast)}(0)\| + \nu_{-}^{-2(J+kp)}),
\end{align*}
for some $0 \le s < p\delta^{(\ast)}+b_{J+kp}$, where we used 
\eqref{eq:boundary-remaining-time} in the second inequality, and 
\eqref{eq:boundary-flow-near-profile} 
in the last inequality. As $t \to T-$, $k \to \infty$, and thus 
$\|(T-t)w(t)-v(t)\|_{\ell^{\infty}(\ZZ_{>0}; \CC)}$ converges to zero. This completes the proof of the proposition.
\end{proof}

\subsection{Completion of the proof of Theorem \ref{thm:main}}
\label{subsec:boundary-main-proof}

\begin{proof}[Proof of Theorem \ref{thm:main}]
Let $\epsilon>0$ and $J$ be as in Proposition \ref{prop:boundary-reconstruction}.
Since $R^{J}Z^{(\ast)}(0)$ has right spatial decay by Theorem 
\ref{thm:orbit-existence},
we may truncate it to obtain real-valued finitely supported initial 
data $\omega(0)$ such that \[
\|\omega(0) - ((R^{J}Z^{(\ast)}(0))_{n})_{n \ge 1}\|_{
\ell^{\infty}(\ZZ_{>0}; \CC)} < \frac{1}{2}\epsilon.
\]

Now let $\widetilde{\omega}(0) \in \ell^{\infty}(\ZZ_{>0}; \CC)$ satisfy \[
\|\widetilde{\omega}(0) - \omega(0)\|_{\ell^{\infty}(\ZZ_{>0}; \CC)} 
< \frac{1}{2}\epsilon.
\]
Then by the triangle inequality, we have \[
\|\widetilde{\omega}(0) - ((R^{J}Z^{(\ast)}(0))_{n})_{n \ge 1}\|_{
\ell^{\infty}(\ZZ_{>0}; \CC)} < \epsilon.
\]
Then, after extending $\widetilde{\omega}(0)$ to nonpositive shells by 
$R^{J}Z^{(\ast)}(0)$, we can apply Proposition 
\ref{prop:boundary-reconstruction} to show that the strong solution 
$\widetilde{\omega} \in C_{t}^{1}([0, T); \ell^{\infty}(\ZZ_{>0}; \CC))$ 
to \eqref{eq:main-renorm-common} blows up at finite time $T$ with the desired $\ell^{\infty}$ 
blowup rate.

Now assume that $\widetilde{\omega}(0) \in \ell^{q}(\ZZ_{>0}; \CC)$ 
for some $1 \le q < \infty$. The continuation criterion 
\eqref{eq:continuation-criterion} and uniqueness show that the strong 
solution 
$\widetilde{\omega} \in C_{t}^{1}([0, T); \ell^{\infty}(\ZZ_{>0}; \CC))$ is 
also in $C_{t}^{1}([0, T); \ell^{q}(\ZZ_{>0}; \CC))$. Since 
$\|\widetilde{\omega}\|_{\ell^{q}(\ZZ_{>0}; \CC)} \ge 
\|\widetilde{\omega}\|_{\ell^{\infty}(\ZZ_{>0}; \CC)}$, 
we obtain the lower bound for the $\ell^{q}$ blowup rate. We are left to prove 
the upper bound.

Let $v \in L_{t}^{\infty}([0, T); \ell^{1}(\ZZ_{>0}; \CC))$ be given as in 
Proposition \ref{prop:boundary-reconstruction}, corresponding to 
$\widetilde{\omega}$. Observe that 
\begin{align*}
& (T-t)^{2}\|N^{+}(\widetilde{\omega}(t))\|_{\ell^{q}(\ZZ_{>0}; \CC)} \\ 
&= \|N^{+}((T-t)\widetilde{\omega}(t)-v(t) + v(t))\|_{\ell^{q}(\ZZ_{>0}; \CC)} \\
& \le C \|(T-t)\widetilde{\omega}(t)-v(t)\|_{\ell^{\infty}(\ZZ_{>0}; \CC)}
\|(T-t)\widetilde{\omega}(t)\|_{\ell^{q}(\ZZ_{>0}; \CC)} 
+ C \|v(t)\|_{\ell^{q}(\ZZ_{>0}; \CC)}\|(T-t)\widetilde{\omega}(t)\|_{\ell^{\infty}(\ZZ_{>0}; \CC)} \\
& \le \frac{1}{2}(T-t)\|\widetilde{\omega}(t)\|_{\ell^{q}(\ZZ_{>0}; \CC)} + C
\end{align*}
for $t$ sufficiently close to $T$, 
where $C$ is a constant independent of $t$, but may change from line to line, 
and may depend on $\widetilde{\omega}(0)$.
Here, the last inequality follows since 
$(T-t)\widetilde{\omega}(t)-v(t)$ converges to zero in 
$\ell^{\infty}(\ZZ_{>0}; \CC)$ norm, 
$\|v(t)\|_{\ell^{q}(\ZZ_{>0}; \CC)} 
\le \|v(t)\|_{\ell^{1}(\ZZ_{>0}; \CC)}$ has a uniform bound,
and $\widetilde{\omega}$ has the blowup rate $(T-t)^{-1}$ in 
$\ell^{\infty}(\ZZ_{>0}; \CC)$. 

Therefore, we have \[
\left \|\partial_{t}\widetilde{\omega}(t)\right \|_{\ell^{q}(\ZZ_{>0}; \CC)} 
\le \frac{1}{2}\frac{1}{T-t}\|\widetilde{\omega}(t)\|_{\ell^{q}(\ZZ_{>0}; \CC)} 
+ \frac{C}{(T-t)^{2}}
\]
for $t$ sufficiently close to $T$, say for $t>t_{0}$.
By Gr\"onwall's inequality,
\begin{equation*}
\begin{aligned}
\|\widetilde{\omega}(t)\|_{\ell^{q}(\ZZ_{>0}; \CC)} 
&\le \exp\left ( \int_{t_{0}}^{t}\frac{1}{2(T-s)}\dd s \right ) 
\|\widetilde{\omega}(t_{0})\|_{\ell^{q}(\ZZ_{>0}; \CC)}
+ \int_{t_{0}}^{t}\exp\left ( \int_{s}^{t} \frac{1}{2(T-r)}\dd r\right )
\frac{C}{(T-s)^{2}}\dd s \\
&= \left ( \frac{T-t_{0}}{T-t} \right )^{\frac{1}{2}}\|\widetilde{\omega}(t_{0})\|_{
  \ell^{q}(\ZZ_{>0}; \CC)} + \int_{t_{0}}^{t} \left ( \frac{T-s}{T-t} \right )^{\frac{1}{2}}
\frac{C}{(T-s)^{2}}\dd s \\
&\le \left ( \frac{T-t_{0}}{T-t} \right )^{\frac{1}{2}}
\|\widetilde{\omega}(t_{0})\|_{\ell^{q}(\ZZ_{>0}; \CC)}
+ \frac{C}{T-t} \\
& \le \frac{C}{T-t}.
\end{aligned}
\end{equation*}
This completes the proof of the blowup rate in $\ell^{q}$.

Finally, let $T_{\ast}$ be the blowup time for the finitely supported 
initial data $\omega(0)$. The scaling $\omega(t) \mapsto T_{\ast}\omega(T_{\ast}t)$ 
sends its solution to one blowing up at time $1$ and sends the above 
neighborhood onto the ball of radius $T_{\ast}\epsilon/2$ centered at 
$T_{\ast}\omega(0)$. Each perturbed blowup time becomes $T/T_{\ast}$, 
and all the asserted blowup rates are preserved. This completes the 
proof of Theorem \ref{thm:main}.
\end{proof}

\appendix

\section{Construction of the numerical data}
\label{app:numerical-approximation}

This appendix describes the construction of the numerical data
$(\overline{\delta}, \overline{a})$ used in Section \ref{sec:bilateral}. 
All code used in these computations is available both in the 
supplementary material and in the GitHub repository \cite{CAParchive}.

Recall that we seek an approximation to a pair $(\delta, a)$ satisfying 
\eqref{eq:orbit-nonlinear-matching-condition} and 
\eqref{eq:orbit-nonlinear-Taylor-condition}. Once $\delta$ and 
$a^{(0)}$ are given, 
\eqref{eq:orbit-nonlinear-Taylor-condition} uniquely determines 
$a^{(k)}$ for all $k \ge 1$, so it remains only to check
\eqref{eq:orbit-nonlinear-matching-condition}. Since $a^{(0)}$ 
represents $Z(0)$ and \eqref{eq:orbit-nonlinear-matching-condition}
expresses the relative periodicity condition
\begin{equation}\label{eq:orbit-relative-periodicity}
Z(\delta)=RZ(0),
\end{equation} 
we are led to search for approximations to an initial profile $Z(0)$
and a relative period $\delta$ satisfying 
\eqref{eq:orbit-relative-periodicity}.

Undoing the renormalization in Subsection \ref{subsec:renorm},
we return to the original variable \[
z(t)=\frac{1}{1-t} Z(-\log(1-t)), \quad t=1-e^{-\tau}.
\] 
Then $z$ solves $\partial_{t}z = N_{\mathrm{real}}(z)$, and
\eqref{eq:orbit-relative-periodicity} becomes \[
z(1-e^{-\delta}) = e^{\delta}Rz(0).
\]
This suggests looking for a solution $z$, two times $t_{l}<t_{r}$,
and an amplitude growth factor $A>1$ such that
\begin{equation}\label{eq:orbit-original-relative-periodicity}
z(t_{r}) \approx ARz(t_{l}).
\end{equation}

We solve $\partial_{t}z=N_{\mathrm{real}}(z)$ from the simple initial data \[
z_{n}(0) = \begin{cases}
1 & n=0, 1, \\
0 & \text{otherwise},
\end{cases}
\]
in the hope that the solution approaches the self-similar form we seek.
We use the fourth-order Runge-Kutta method with an adaptive timestep, 
where we truncate the system to a finite interval $I$ of shell 
indices $n$.

During this computation, we search for $t_{l}, t_{r}, A$ satisfying 
\eqref{eq:orbit-original-relative-periodicity} as follows. Starting 
at $t_{0}=0$, we record successive times $t_{j+1}$ at which the 
weighted center defined by \[
\frac{\sum_{n \in I}n|z_{n}|^{2}}{\sum_{n \in I} |z_{n}|^{2}}
\]
has moved approximately one shell to the right from its value at time $t_{j}$.
For each $j=0, 1, 2, \ldots$, we set \[
A= \frac{\|z(t_{j+1})\|_{\ell^{2}}}{\|z(t_{j})\|_{\ell^{2}}}
\]
and compare $z(t_{j})$ with $A^{-1}R^{-1}z(t_{j+1})$ in the
$\ell^{2}$ norm, relative to $\|z(t_{j})\|_{\ell^{2}}$.
We repeat this comparison until the profiles are 
sufficiently close. We denote the accepted pair $(t_{j}, t_{j+1})$ by 
$(t_{l}, t_{r})$, and take $A$ to be the corresponding amplitude ratio.

After shifting the time origin to $t_{l}$, we obtain an approximate 
solution $z$ of $\partial_{t}z=N_{\mathrm{real}}(z)$ such that \[
z(T) \approx ARz(0),
\] where $T=t_{r}-t_{l}$. We now rescale $z$ so that its estimated 
blowup time is $1$. Over a time interval of length $T$, the profile 
is approximately amplified by a factor of $A$ and shifted one shell 
to the right. Since $N_{\mathrm{real}}$ is quadratic and commutes with 
$R$, the same amplification
and shift are expected to occur after a further time $T/A$. The estimated 
blowup time is therefore \[
T+\frac{T}{A}+\frac{T}{A^{2}}+\cdots = \frac{AT}{A-1}.
\] Hence the rescaled function \[
t \mapsto \frac{AT}{A-1}z\left (\frac{AT}{A-1}t \right )
\] gives an approximate solution of $\partial_{t}z=N_{\mathrm{real}}(z)$ with an 
estimated blowup time of $1$. Thus, by setting \[
Z(0) = \frac{AT}{A-1}z(0), \quad \delta = \log A,
\] we obtain an approximate solution to \eqref{eq:orbit-relative-periodicity},
and hence approximations to $a^{(0)}$ and $\delta$. For the data 
used in Section \ref{sec:bilateral}, we relabel the shell indices 
by an integer shift so that the weighted center of $a^{(0)}$ lies in 
$[0, 1)$, and denote the shifted interval by $I$. Together with 
finitely many coefficients $a^{(k)}$, say for $1 \le k \le K$, computed 
using \eqref{eq:orbit-nonlinear-Taylor-condition}, with coefficients
outside $I$ set to zero, these give an approximation 
$(\overline{\delta}, \overline{a})$.

The computations described above are implemented in 
\texttt{approximation.py}. The function \texttt{generate\_approximation} 
returns $I, K, \overline{\delta}$, and $\overline{a}$. These data 
are stored in \texttt{approximation.npz}. These computations, 
including saving the data, took approximately $25$ seconds 
on an Intel Core Ultra 5 125H processor.

\bibliographystyle{amsplain}
\bibliography{shell_references_audited}

\end{document}